\documentclass[12pt,a4paper,reqno]{amsart}
\usepackage[T1]{fontenc}
\usepackage[utf8]{inputenc}
\usepackage{amsmath,amssymb,amsthm,mathtools}
\usepackage{geometry}
\usepackage{hyperref}
\usepackage{mathrsfs}
\usepackage{graphicx} 
\usepackage{subcaption} 
\usepackage{tikz}
\usetikzlibrary{arrows.meta, calc, math}

\hypersetup{
	colorlinks=true,
	linkcolor=blue,
	citecolor=blue,
	urlcolor=blue,
	pdftitle={Time-optimal travel driven by two vector fields on a Riemannian manifold},
	pdfauthor={Erasmo Caponio, Dario Corona, Ebtsam H. Taha}
}
\numberwithin{equation}{section}

\newtheorem{theorem}{Theorem}[section]
\newtheorem{proposition}[theorem]{Proposition}
\newtheorem{lemma}[theorem]{Lemma}
\newtheorem{corollary}[theorem]{Corollary}
\theoremstyle{definition}

\theoremstyle{remark}
\newtheorem{remark}[theorem]{Remark}

\newcommand{\R}{\mathbb{R}}

\title[Time-optimal travel driven by two vector fields]{Time-optimal travel driven by two vector fields on a Riemannian manifold}

\author[E. Caponio]{Erasmo Caponio}
\address{Dipartimento di Meccanica, Matematica e Management, Politecnico di Bari, Bari, Italy}
\email{erasmo.caponio@poliba.it}

\author[D. Corona]{Dario Corona}
\address{School of Science and Technology, University of Camerino, Camerino, Italy}
\email{dario.corona@unicam.it}

\author[E. H. Taha]{Ebtsam H. Taha}
\address{Department of Mathematics, Faculty of Science, Cairo University, Giza, Egypt}
\email{ebtsam.taha@sci.cu.edu.eg}


\date{}

\begin{document}

\begin{abstract}
We study minimum-time travel on a Riemannian manifold $(M,h)$ for a
kinematic model in which finitely many drift fields contribute additively
to the velocity, after independent rescaling of their components parallel
and orthogonal to a chosen unit direction $u$. The velocity is then
\[
	\mathcal Y_x(u)
	=\bigl(\eta(x)+h_x(\widehat Z_x,u)\bigr)u+\widehat W_x,
	\qquad u\in UT_xM,
\]
for two effective vector fields $\widehat Z$, $\widehat W$ and a
self-propulsion speed $\eta\geq0$. Let $F$ be the gauge of the radial hull
$\mathcal B_x=[0,1]\cdot\mathcal S_x$ of the velocity set
$\mathcal S_x=\mathcal Y_x(UT_xM)$.
Where $\eta>0$, we characterize when $F$ is a Finsler metric and give
two-sided estimates and representation formulae for it; where $\eta=0$ and
$\widehat Z\neq0$, we recover the Randers, Kropina and strong-wind regimes
of Zermelo navigation. Examples include Matsumoto's slope metric,
traction-modified slope models, and alternative models for wildfire
propagation and microswimmer navigation.
Since $F$ need not be a Finsler metric, we treat the minimum-time problem
with tools from differential inclusions. For distinct endpoints, the inclusions defined by $\mathcal S_x$ and
$\mathcal B_x$ give the same infimum travel time, and so does the convex
hull of $\mathcal B_x$ under local controllability at the terminal point.
When $F$ is finite, the infimum is attained if $(M,h)$ is complete,
$\mathcal Y$ is bounded and each $\mathcal B_x$ is convex.
When $F$ is a Finsler metric on $\{\eta>0\}$, the set $\{\eta=0\}$ cannot
be reached from $\{\eta>0\}$ in finite time and, on $\{\eta>0\}$, the
minimum time is the Finsler distance; if moreover $(M,h)$ is complete
and $\eta$ is bounded, $F$ is forward and backward complete and the
minimum time is attained by geodesics.
We also compute the support Hamiltonian directly from the velocity law,
for every $\eta\geq0$.
\end{abstract}

\maketitle

\tableofcontents

\section{Introduction}

Many models of navigation describe an agent that moves with its own speed
along a chosen direction, while external fields act on its velocity.
A walker on a mountain slope is slowed down uphill and sped up downhill by
gravity~\cite{Matsumoto1989,ParkinsonEtAl2018}; a ship or an aircraft is
transported by a current or a wind, as in Zermelo's navigation
problem~\cite{BaoRoblesShen2004,CaponioJavaloyesSanchez2024}; the front of
a wildfire spreads faster uphill and downwind~\cite{Richards1990,
Anderson1983,Finney1998,JavaloyesPendasSanchez2023};
a swimming microorganism is advected by the surrounding flow and settles
under gravity~\cite{PedleyKessler1987,VentrellaEtAl2023}, and its optimal
navigation has been widely studied~\cite{ColabreseEtAl2017,QiuEtAl2022,
DaddiMoussaIderEtAl2021}. In several of these models, the external field
does not simply add to the velocity, since its components along and across
the direction of motion act differently. For instance, in Matsumoto's slope
model the component of gravity orthogonal to the direction of motion is
cancelled by the grip of the walker on the ground, and only the parallel
one affects the motion, whereas the settling velocity of a spheroidal
swimmer depends on the angle between its axis and the vertical.

In all these situations, time-optimal motion is naturally described by a
field of admissible velocity sets. The background Riemannian metric
determines lengths and orthogonality, but the velocities that can be
realized in one unit of time may depend on the chosen direction and on the
external fields. When the admissible velocity sets are smooth, strongly
convex, and contain the origin in their interiors, their fiberwise gauges
define a Finsler metric; however, there is no reason for these properties
to hold in general. In this paper, we take the kinematic law, rather than
a metric prescribed in advance, as the primary datum.

\subsection{The velocity law}
Let $(M,h)$ be a Riemannian manifold of dimension at least $2$,
let $\alpha_x(y)=\sqrt{h_x(y,y)}$,
and let
$u\in UT_xM=\{v\in T_xM:\alpha_x(v)=1\}$
represent the chosen control direction.
We consider a non-negative self-propulsion speed $\eta$ and
finitely many drift fields whose components parallel and orthogonal to $u$
may be rescaled independently. The model is kinematic: the rescaled
components of the fields contribute additively to the velocity. Their
combined effect can always be encoded by two effective smooth vector
fields $\widehat Z$ and $\widehat W$, and the velocity law takes the form
\begin{equation}\label{eq:generalForm}
	\mathcal Y_x(u)
	=\bigl(\eta(x)+h_x(\widehat Z_x,u)\bigr)u+\widehat W_x.
\end{equation}
The field $\widehat Z$ produces a direction-dependent radial deformation,
whereas $\widehat W$ translates the resulting velocity set. The reduction
to these two fields is algebraic, but it separates two geometrically
different effects and remains meaningful even when $\eta$ vanishes.

For each $x\in M$, the maximal velocities form the set
$\mathcal S_x:=\mathcal Y_x(UT_xM)$. We call them maximal because they
are not scaled down, although a ray issuing from $0_x$ may meet
$\mathcal S_x$ more than once, in which case only its farthest
intersection is the largest velocity in that direction. We assume that
the agent can also
move more slowly than a maximal velocity in the same direction, possibly
coming to rest. For instance, a walker on a slope can brake, and in
wildfire models the burned region does not shrink in time, which requires
the zero velocity to be admissible~\cite{Bressan2007}. Thus, the
admissible velocities form the radial set
\begin{equation}\label{admissible}
	\mathcal B_x
	:=\operatorname{rad}\mathcal S_x
	=\{\lambda\mathcal Y_x(u):
	u\in UT_xM,\ 0\leq\lambda\leq1\}.
\end{equation}
Other conventions are possible: for instance, in Zermelo's problem only
the self-propulsion can be reduced, while the drift remains fully active.
However, $\mathcal B_x$ is the largest set of velocities that do not
exceed, along their own direction, a maximal one, and every choice of
admissible velocities either in $\mathcal S_x$ or  $\mathcal B_x$ yields the same infimum travel time
between distinct points (see Theorem~\ref{maximalspeed}). The choice only
affects the possibility of remaining at rest, and hence the reachable
sets (see Subsection~\ref{subsec:additional-drift}). This is a feature of
autonomous laws: when the velocity law depends on time, waiting may
shorten the travel time (see Remark~\ref{rem:time-dependent}).
The gauge function 
\[
F(y):=\inf\{\lambda>0:y\in\lambda\mathcal B_x\}
\]
of $\mathcal B_x$ measures the time needed to realize the displacement
$y\in T_xM$, and a natural question  is when $F$ is a Finsler metric.

\subsection{Finsler regimes}
On the region
\[
M_+:=\{x\in M:\eta(x)>0\},
\]
we introduce the normalized fields
\[
Z:=\frac{\widehat Z}{\eta},
\qquad
W:=\frac{\widehat W}{\eta}.
\]
The untranslated radial deformation is governed by the Matsumoto metric
associated with $Z$. Proposition~\ref{prop:FinslerM-conditions} gives the
sharp strong-convexity condition $\alpha(Z)<1/2$; beyond this threshold,
the Matsumoto gauge is no longer a Finsler metric and, for
$\alpha(Z)>1/2$, it is not even convex. For $\alpha(Z)<1/2$,
Proposition~\ref{generalcase} shows that $\mathcal B_x$ is the unit ball
of a Finsler norm $F$ precisely when the opposite translation vector
$-W_x$ belongs to the interior of the untranslated velocity ball. Both
conditions are also necessary for $F$ to be a Finsler metric on $M_+$
(Corollary~\ref{cor:F-Finsler}). The
resulting metric is generally defined only implicitly. Nevertheless, we
obtain semi-explicit fiberwise representation formulae, both when $Z_x$
and $W_x$ are collinear and when they are linearly independent (see
Theorems~\ref{thm:single} and~\ref{thm:two} in Appendix~\ref{appendix}).

The kinematic formulation also covers points at which self-propulsion
vanishes. If $\eta(x)=0$ and $\widehat Z_x=0_x$, then the admissible
velocity set is simply the segment $[0_x,\widehat W_x]$. If instead
$\widehat Z_x\neq0_x$, the maximal velocities form a translated Riemannian
sphere, and $\mathcal B_x$ is its radial hull.
Proposition~\ref{zeroeta-classification} classifies the corresponding
gauge according to whether the origin lies inside, on, or outside the
enclosed ball. These alternatives produce, respectively, Randers, Kropina,
and conic Minkowski gauges, and they agree with the mild-, critical-, and
strong-wind alternatives of Zermelo navigation, in the framework of wind
Finslerian structures~\cite{CaponioJavaloyesSanchez2024}. In the
strong-wind case, however, the radial hull is larger than the translated
velocity ball, and its gauge represents the least-time branch associated
with the farthest intersection of each admissible ray with $\mathcal S_x$.

\subsection{Examples}
Besides Matsumoto's slope metric and Zermelo navigation mentioned above 	(see also Subsections~\ref{subsec:Matsumoto} and~\ref{subsec:zermelo}), several models fit within this framework.
\begin{itemize}
\item The one-field case contains the traction-modified slope models,
	including the slippery-slope and superwind regimes,
	studied by Aldea and Kopacz~\cite{AldeaKopacz2025a,AldeaKopacz2026}
	(see Subsection~\ref{subsec:onevector}).
	\item The general model with two fields yields an alternative kinematic model
	for anisotropic wildfire propagation,
	which combines the Matsumoto-type slope term of~\cite{JavaloyesPendasSanchez2023}
	with a wind acting as an additive translational drift,
	as in \cite{Markvorsen2016,Dehkordi2022},
	(see Subsection~\ref{wildfireJPS}).
	\item A time-dependent law with two noncollinear fields arises
	from a free-orientation relaxation of the settling model
	for spheroidal microswimmers in~\cite{VentrellaEtAl2023}.
	This relaxation replaces the dynamically constrained body orientation
	by a freely selectable direction,
	and leads to a natural minimum-time navigation problem,
	as considered, for different hydrodynamic velocity laws,
	in~\cite{DaddiMoussaIderEtAl2021}
	(see Subsection~\ref{subsec:microswimmers}).
\end{itemize}

\subsection{The minimum-time problem}
The main motivation for the velocity law~\eqref{eq:generalForm} is, however,
the time-optimal problem,
which lies at the core of almost all the applications mentioned above.
Since $F$ need not be a Finsler metric,
time-minimizing trajectories cannot in general be described as Finsler geodesics.
For $x_1,x_2\in M$ and $T>0$, let us consider the class of admissible curves
\[
	\Gamma_T(x_1,x_2)
	:=\bigl\{\gamma\in AC([0,T],M):
	\gamma(0)=x_1,\ \gamma(T)=x_2,\
	\dot\gamma(t)\in\mathcal B_{\gamma(t)}\ \text{for a.e.}\ t\in[0,T]\bigr\},
\]
and the minimum-time problem
\begin{equation}\label{TT}\tag{TT}
	\mathcal T(x_1,x_2):=\inf\bigl\{T>0:\Gamma_T(x_1,x_2)\neq\emptyset\bigr\},
\end{equation}
with the convention that $\mathcal T(x_1,x_2)=+\infty$
if $\Gamma_T(x_1,x_2)=\emptyset$ for every $T>0$.
This formulation,
which falls within the theory of differential inclusions~\cite{AubinCellina}
and of minimum-time control~\cite{BardiCapuzzoDolcetta1997,CannarsaSinestrari2004},
does not require $\mathcal B_x$ to be convex, strongly convex,
or even to contain the origin in its interior.
A closely related framework was introduced by Bressan~\cite{Bressan2007},
who modelled the region burned by a wildfire
as the reachable set of a differential inclusion $\dot x\in \mathcal{F}(x)$,
with compact and convex values containing the origin,
and has since been developed in the study of fire confinement
problems~\cite{BressanEtAl2008,BressanDeLellis2009,BianchiniZizza2025}.
In our setting, the admissible velocity sets always contain the origin,
but they need not be convex.

A central part of the paper concerns the relation
among the control problems generated by the maximal velocities $\mathcal S_x$,
their radial hull $\mathcal B_x$,
and their convex hull $\widehat{\mathcal B}_x:=\operatorname{co}\mathcal B_x$.
Denoting by $\mathcal T_{\mathcal S}$ and $\widehat{\mathcal T}$
the infimum travel times of the differential inclusions
generated by $\mathcal S_x$ and by $\widehat{\mathcal B}_x$, respectively,
we obtain the following results.
\begin{enumerate}
	\item For distinct endpoints, $\mathcal T_{\mathcal S}(x_1,x_2)=\mathcal T(x_1,x_2)$
	(Theorem~\ref{maximalspeed}).
	The proof is based on a measurable reparametrization
	of a $\mathcal B$-admissible trajectory,
	and requires neither convexity nor Finsler regularity.
	On the diagonal, the same conclusion holds
	under the condition $0_x\in\mathcal S_x$.
	\item If $0_{x_2}\in\operatorname{int}\widehat{\mathcal B}_{x_2}$,
	then $\widehat{\mathcal T}(x_1,x_2)=\mathcal T(x_1,x_2)$
	(Theorem~\ref{convexification}).
	The Filippov--Wa\.zewski relaxation theorem gives uniform approximation
	by trajectories of the original inclusion,
	but does not in general preserve the terminal point. 
	Condition $0_{x_2}\in\operatorname{int}\widehat{\mathcal B}_{x_2}$ gives 
	Petrov's condition which allows one to control the terminal error
	in a time tending to zero.
	\item If $(M,h)$ is complete, $\mathcal S_x$ is uniformly bounded with respect to $x\in M$,
	and every $\mathcal B_x$ is convex,
	then every finite infimum travel time between distinct endpoints
	is attained by a maximal-velocity trajectory
	(Proposition~\ref{nonfinsler-existence}).
	Existence results for minimum-time problems whose prescribed velocity sets are not convex, such as the brachistochrone problem,  were obtained in~\cite{CellinaFerrieroMarchini2006}. In our setting, the sets $\mathcal S_x$ may be nonconvex, and the main assumption of~\cite{CellinaFerrieroMarchini2006} reduces to the convexity of the radial hulls $\mathcal B_x$; the optimal trajectory still has velocities in $\mathcal S_x$.
	For the law~\eqref{eq:generalForm}, $\mathcal B_x$ is convex
	wherever $\eta(x)=0$ or $\alpha(\widehat Z_x)\leq\eta(x)/2$,
	whatever $\widehat W_x$ is (Corollary~\ref{cor:model-existence}).
\end{enumerate}
In particular, for $x_1\neq x_2$ and under the assumption $0_{x_2}\in\operatorname{int}\widehat{\mathcal B}_{x_2}$,
the three choices of admissible velocities give the same infimum travel time:
\[
	\mathcal T_{\mathcal S}(x_1,x_2)
	=\mathcal T(x_1,x_2)
	=\widehat{\mathcal T}(x_1,x_2).
\]
The results of Subsection~\ref{timeoptimatraj} also hold when, at each point $x$, the control directions are restricted to a nonempty compact set $K_x\subseteq UT_xM$, with $\mathcal S_x:=\mathcal Y_x(K_x)$ and $\mathcal B_x$ its radial hull, provided that, locally, $K_x=\Phi_x(\Theta)$ for a fixed compact set $\Theta\subset\mathbb R^m$ and a map $(x,\theta)\mapsto\Phi_x(\theta)\in UT_xM$ which, in local trivializations, is continuous and Lipschitz in $x$ uniformly in $\theta$. Indeed, the proofs of Theorem~\ref{maximalspeed} and Proposition~\ref{nonfinsler-existence} only use that $x\mapsto\mathcal S_x$ has nonempty compact values and is locally Lipschitz with respect to the Hausdorff distance, as in Lemma~\ref{hausdorff}, whose proof carries over with $\Theta$ in place of $\mathbb S^{n-1}$. The proof of Theorem~\ref{convexification} also uses a representation of the dynamics as a control system with a fixed compact control set, which can then be provided by $[0,1]\times\Theta$.

When  the gauge function $F$ is a Finsler metric on $M_+$,
the time-optimal problem can instead be studied
with the tools of Finsler geometry,
provided that the behaviour near the zero set of $\eta$ is controlled.
Proposition~\ref{2bounds} gives two-sided estimates
comparing $F$ with $\alpha/\eta$,
and Theorem~\ref{bfcompleteness} shows that,
if $(M,h)$ is complete and $\eta$ is bounded from above,
then $(M_+,F)$ is both forward and backward complete.
Even without these global assumptions, $\{\eta=0\}$ lies at infinite travel time from $M_+$
(Proposition~\ref{unreachable}) and,
for points in the same connected component of $M_+$,
the minimum time agrees with the Finsler distance. Under the assumptions of Theorem~\ref{bfcompleteness}, it is attained, for distinct endpoints,
by a unit-speed minimizing geodesic (Corollary~\ref{cor:geodesics}).

\subsection{The support Hamiltonian}
On the dual side, we compute the support Hamiltonian
$H(p):=\max_{z\in\mathcal B_x}p(z)$, $p\in T_x^*M$,
directly from the velocity law, for every $\eta\geq0$
(Proposition~\ref{support-general}).
Its relevance comes from the eikonal equation. Let $\tau_{x_0}:=\mathcal T(x_0,\cdot)$ be the minimum-time function from $x_0$. Then $H(d\tau_{x_0})=1$ at every point $x\neq x_0$ reachable from $x_0$ where $\tau_{x_0}$ is differentiable. Moreover, if $0_{x_0}\in\operatorname{int}\widehat{\mathcal B}_{x_0}$, the set of points reachable from $x_0$ is open, $\tau_{x_0}$ is continuous on it, and $\tau_{x_0}$ is a viscosity solution of $H(d\tau_{x_0})=1$ on this set minus $\{x_0\}$ (see, e.g.,~\cite{BardiCapuzzoDolcetta1997,CannarsaSinestrari2004}).
Since a linear functional has the same maximum on a set and on its convex hull,
$H$ does not distinguish $\mathcal B_x$ from $\widehat{\mathcal B}_x$,
consistently with Theorem~\ref{convexification}.
In the regular Finsler regime, $H$ is the dual metric $F^*$,
and its computation reduces to a one-dimensional maximization
(see Proposition~\ref{prop:Hscalar}).
We also consider time-dependent data,
and we compare the radial-hull convention
with the translated-ball convention for additive drift
(Subsection~\ref{subsec:additional-drift}).

The eikonal and level-set formulations are also the basis of the numerical computation of arrival times, for instance through level-set methods~\cite{OsherSethian1988}, ordered upwind methods for anisotropic static Hamilton--Jacobi equations~\cite{SethianVladimirsky2003}, and Hamiltonian fast marching for anisotropic and non-holonomic eikonal equations~\cite{MirebeauPortegies2019}; methods of this kind have been applied to optimal navigation on steep terrain~\cite{ParkinsonEtAl2018}. Since our formulae express $H$ directly in terms of $\eta$, $\widehat Z$ and $\widehat W$, they can be used in schemes that admit  a  convex Hamiltonian also beyond the regular Finsler regime, provided that $0_x\in\operatorname{int}\widehat{\mathcal B}_x$ at every point of the computational domain, that is, $H(p)>0$ for every $p\neq0$ (see Lemma~\ref{separation}). In this case, by Theorem~\ref{convexification}, the computed arrival time is the infimum travel time of the original problem, even where $\mathcal B_x$ is not convex. 
\subsection{Overview}
The paper is organized as follows.
Section~\ref{sec general results} collects general facts
about admissible velocity sets, their radial hulls, and their gauge functions,
including the gauge of the radial hull of a translated Riemannian sphere.
Section~\ref{model} derives the velocity law~\eqref{eq:generalForm}
from finitely many drift fields
and studies the fiberwise geometry of the admissible velocity sets,
in the regimes of positive and vanishing self-propulsion.
Section~\ref{examples} discusses the models that fit within our framework:
Matsumoto's slope metric, Zermelo navigation,
the traction-modified slope models,
wildfire propagation, and microswimmer navigation.
Section~\ref{sec:time-travel} is devoted to the minimum-time problem~\eqref{TT}:
we compare the maximal, radial, and convexified admissible velocity sets,
we prove the existence of time-optimal trajectories,
and we study the completeness of $(M_+,F)$ in the regular Finsler regime.
Section~\ref{hamiltonian} develops the support-Hamiltonian formulation,
including its scalar representation in the regular regime,
time-dependent laws, and additive drifts.
Finally, Appendix~\ref{appendix} derives the fiberwise representation
formulae for $F$ in the one-field and two-field cases.

\section{Admissible velocities and gauge functions}\label{sec general results}

In this section, we recall and prove general results
about admissible velocity sets and their gauge functions.
This means that, within this section,
we do not assume the admissible velocity sets
to be defined by a specific velocity law
(e.g., as in~\eqref{admissible}),
but rather consider general subsets of the tangent space.
All results in this section are pointwise,
so we fix a point $x \in M$
and we omit the dependence on $x$
for all sets and functions.

Throughout the paper, $(M,h)$ is a smooth Riemannian manifold
of dimension $n\geq2$, $\alpha$ is the induced norm,
$\alpha_x(y):=\sqrt{h_x(y,y)}$, and
$UT_xM:=\{v\in T_xM:\alpha_x(v)=1\}$.

Let us first fix some terminology.
For the standard notions of convex hull
and gauge we refer, for instance,
to \cite[\S 1.6--1.7]{Schneider2014}.
A subset $K\subset T_xM$ is \emph{radial} (with respect to $0$)
if $tz\in K$ whenever $z\in K$ and $t\in[0,1]$;
the \emph{radial hull} of a subset $S\subset T_xM$ is
\[
	\operatorname{rad}S:=\{tz:z\in S,\ 0\leq t\leq1\},
\]
i.e. the smallest radial set containing $S$.
Notice that $\operatorname{rad}S$ is compact whenever $S$ is,
being the image of the compact set $[0,1]\times S$
under the continuous map $(t,z)\mapsto tz$.
The sets $\operatorname{cone}S$ and $\operatorname{co}S$ are
respectively the {\em generated cone}  of $S$
\[
	\operatorname{cone}S:=\{tz:z\in S,\  t\geq 0\},
\]
and  the convex hull of $S$.
For any $K \subset T_xM$,
its \emph{gauge function} $F\colon T_xM \to [0,+\infty]$
is defined as 
\begin{equation}
	\label{eq:def-F}
	F(y):=
	\inf\bigl\{\lambda>0:
	y\in\lambda K\bigr\}\in[0,+\infty].
\end{equation}
We call
\[
	\Lambda:=\{y\in T_xM:F(y)<+\infty\}
\]
the \emph{finiteness cone} of $F$ at $x$.
Note that no convexity is assumed anywhere in this section,
and that $F$ is allowed to take the value $+\infty$.
All the following lemmas are well-known;
we isolate and prove them for the reader's convenience.

\begin{lemma}\label{radialgauge}
	Let $K \subset T_xM$ be nonempty, compact and radial.
	Let $F$ be its gauge function.
	Then $0\in K$, $F$ is positively homogeneous,
	$F(0)=0$, $F>0$ on $T_xM\setminus\{0\}$,
	and $\Lambda=\operatorname{cone}(K)$.
	Moreover, 
	\begin{enumerate}
		\item[(i)] for every $y\in\Lambda\setminus\{0\}$
			the infimum defining $F(y)$ is attained;
			more precisely, $1/F(y)$ is the largest $s>0$
			such that $s y\in K$,
			and therefore
			\[
		\frac{y}{F(y)}\in\partial K;
		\]
	\item[(ii)] $K=\bigl\{y\in T_xM:F(y)\leq1\bigr\}$.
\end{enumerate}
\end{lemma}
\begin{proof}
As $0\in K$, we immediately get $F(0)=0$.
Let $y\neq0$.
If $F(y)=0$, there are arbitrarily small $\lambda>0$
with $y/\lambda\in K$, contradicting the boundedness of $K$;
hence $F(y)>0$.
Both the identity $\Lambda=\operatorname{cone}(K)$
and the positive homogeneity of $F$ are a restatement
of the definition of the gauge function,
since $y\in\Lambda$
if and only if $y\in\lambda K$ for some $\lambda>0$. 
	
(i) Let $y\in\Lambda\setminus\{0\}$ and let $(\lambda_j)$
be a sequence with $y\in\lambda_jK$ and $\lambda_j\to F(y)>0$.
Then $y/\lambda_j\in K$ and $y/\lambda_j\to y/F(y)$,
so $y/F(y)\in K$ because $K$ is closed,
hence the infimum is attained.
Since $sy\in K$ if and only if $1/s\geq F(y)$,
the largest such $s$ is $1/F(y)$;
in particular $(1+\varepsilon)y/F(y)\notin K$
for every $\varepsilon>0$,
so $y/F(y)$ is a point of $K$ which
is a limit of points of the complement of $K$,
that is, $y/F(y)\in\partial K$.
	
(ii) If $y\in K$, then $y\in1\cdot K$ and $F(y)\leq1$.
Conversely, let $F(y)\leq1$.
If $y=0$, then $y\in K$.
Otherwise $y\in\Lambda\setminus\{0\}$, so $y/F(y)\in K$
by (i), and radiality gives $y=F(y)\bigl(y/F(y)\bigr)\in K$,
since $0<F(y)\leq1$.
\end{proof}

Observe that only the boundedness of $K$ is used for $F>0$,
only its closedness for the attainment in (i),
and only its radiality for the nontrivial inclusion in (ii).
Note also that $\partial K$ in general 
strictly contains  $\{y\in T_xM:F(y)=1\}$:
if $0\in\partial K$,
the boundary of $K$ contains points $z$ with $F(z)<1$.

\begin{remark}[Gauge functions and Okubo's technique]
The procedure commonly referred to as \emph{Okubo's technique} 	(see, for instance, \cite[Exercise~1.2.8(a)]{bcs}) consists in starting 	from an equation
$Q(y)=0$, 	 and determining its gauge $F$ from
\begin{equation}\label{Okubo-equation}
Q\left(\frac{y}{F(y)}\right)=0.
\end{equation}
This is precisely what Lemma~\ref{radialgauge}-(i) establishes. Indeed, if $K\subset T_xM$ is nonempty, 	compact and radial and $F$ is its gauge, then
\[
\frac{y}{F(y)}\in\partial K
\]
for every $y\in\operatorname{cone}(K)\setminus\{0\}$. Consequently, 	whenever $\partial K$ is described by the equation $Q=0$, 	\eqref{Okubo-equation} follows immediately. 	
 A qualification is needed only if one attempts to determine $F(y)$ by 	solving~\eqref{Okubo-equation} without specifying which solution must be 	chosen.
 If a ray meets the prescribed hypersurface more than once, the equation may have several positive solutions.
 The gauge is determined by 	the farthest intersection along the ray, equivalently by the smallest 	positive solution.
 Proposition~\ref{minimalgauge} below formalizes precisely this selection rule for translated compact hypersurfaces.
\end{remark}

We now specialize Lemma~\ref{radialgauge}
to the case in which the compact radial set is
the radial hull of a translate
of the boundary of another compact radial set,
which is the structure of several admissible
velocity sets considered in this paper. 
\begin{proposition}	\label{minimalgauge}
	Let $D\subset T_xM$ be nonempty,
	compact and radial and let $W\in T_xM$.
	Let $K :=\operatorname{rad}(\partial D+W)$ and
	let $F$ and $\Lambda$ be its gauge function
	and its finiteness cone, respectively.
	Then, for every $y\in\Lambda\setminus\{0\}$,
	\begin{equation}\label{gaugeboundary}
		\frac{y}{F(y)}\in\bigl(\partial D+W\bigr)\cap\partial K.
	\end{equation}
	Moreover, let $G$ be the gauge function of $D$ and assume that 
	\begin{equation}\label{noflat}
		G^{-1}(1)=\partial D\setminus\{0\}.
	\end{equation}
	Then, for every $y\in\Lambda\setminus\{0\}$ such that
	$y/F(y)\neq W$, the number $t_0:=F(y)$
	is the smallest positive solution of
	\begin{equation}\label{implicit}
		G(y-tW)=t ;
	\end{equation}
	if, moreover, every ray originating from $0$
	meets $\partial D+W$ at most once, then it is the unique one.
\end{proposition}
\begin{proof}
	The set $\partial D+W$ is compact, so $K$,
	being the image of $[0,1]\times(\partial D+W)$
	under the continuous map $(t,z)\mapsto tz$, is compact;
	Lemma~\ref{radialgauge} therefore applies and gives,
	in particular, $\Lambda
	=\operatorname{cone}(K)=\operatorname{cone}(\partial D+W)$.
	Fix $y\in\Lambda\setminus\{0\}$.
	By Lemma~\ref{radialgauge}-(i),
	$1/F(y)$ is the largest $s>0$ with $sy\in K$,
	and $y/F(y)\in\partial K$.
	As $K:=\operatorname{rad}(\partial D+W)$,
	$sy\in K$ if and only if $s'y\in\partial D+W$
	for some $s'\geq s$.
	The largest of such $s$ must be equal to $s'$,
	otherwise any  $\bar s\in \big(1/F(y), s'\big)$ would give 
	$\bar sy\in K$, a contradiction.
	Thus, $y/F(y)\in \partial D+W$ as well,
	proving \eqref{gaugeboundary}.

	Assume now~\eqref{noflat} and $y/F(y)\neq W$.
	By~\eqref{gaugeboundary}, the vector $z:=y/F(y)-W$
	lies in $\partial D\setminus\{0\}$,
	so that $G(z)=1$.
	Since $z=\bigl(y-F(y)W\bigr)/F(y)$ and $F(y)>0$,
	positive homogeneity gives $G\bigl(y-F(y)W\bigr)=F(y)$,
	that is, $t_0=F(y)$ solves~\eqref{implicit}. 
	If  $t>0$ is any other solution of~\eqref{implicit},
	then positive homogeneity gives $G(y/t-W)=1$,
	thus  $y/t-W\in\partial D$, i.e.  $y/t\in\partial D+W\subseteq K$.
	By the definition of the gauge, $F(y)\leq t$,
	so that $t_0$ is the smallest positive solution.
	The uniqueness is immediate, since distinct positive solutions correspond to distinct intersections of the ray
	$\mathbb R_{>0}y$ with $\partial D+W$.
\end{proof}

\begin{remark}\label{degenerateray}
	The condition $y/F(y)\neq W$ cannot be excluded in general. Indeed, if $y/F(y)=W$ and $t_0=F(y)$, then $y-t_0W=0$, and hence $G(y-t_0W)=G(0)=0\neq t_0$. Thus $t_0$ does not solve~\eqref{implicit}. Notice that $y/F(y)=W$ can occur only if $W\neq0$, $y\in\mathbb R_{>0}W$, and $0\in\partial D$.	Indeed, if $y/F(y)=W$, then $y=F(y)W$. Since $y\neq0$ and
	$F(y)>0$, this implies that $W\neq0$ and
	$y\in\mathbb R_{>0}W$. Moreover, by~\eqref{gaugeboundary}, $W=y/F(y)\in\partial D+W$, and hence $0\in\partial D$. Consequently, the  equality $y/F(y)=W$ cannot occur when $0\in\operatorname{int}D$.
\end{remark}

When the origin belongs to a compact convex set,
the radial hull of its boundary returns the set itself;
in general one still recovers the radial hull of the set.
This elementary fact is what links the definition
of the set of admissible velocities
as a radial hull in \eqref{admissible} to the translated body \eqref{ballidentity}
(see Remark~\ref{identityremark}).
\begin{lemma}\label{radialhull}
Let $C\subset T_xM$ be nonempty and compact. Then:
\begin{enumerate}
\item[(i)] $C\subseteq\operatorname{rad}\partial C$, and consequently
$\operatorname{rad}\partial C=\operatorname{rad}C$;
\item[(ii)] if $C$ is in addition convex, then $\operatorname{rad}\partial C=C$
if and only if $0\in C$.
\end{enumerate}
\end{lemma}
\begin{proof}
(i)  As $\partial C\neq \emptyset$ then  $0\in \operatorname{rad}\partial C$;
let then $y\in C\setminus\{0\}$.
The number 	$t^*:=\sup\{t\geq1:ty\in C\}$ is finite because $C$ is bounded, and it is attained 	because $C$ is closed; since $ty\notin C$ for $t>t^*$, the point $t^*y$ lies in 	$\partial C$, so that $y=\frac{1}{t^*}(t^*y)\in\operatorname{rad}\partial C$.
This 	proves $C\subseteq\operatorname{rad}\partial C$.
Applying $\operatorname{rad}$ and 	using its monotonicity and idempotence gives 	$\operatorname{rad}C\subseteq\operatorname{rad}\partial C$, while the reverse 	inclusion is immediate from $\partial C\subseteq C$.
	
(ii) Assume $0\in C$. Then $C\subseteq\operatorname{rad}\partial C$ by (i),
and the 	reverse inclusion follows from convexity: $\lambda z=\lambda z+(1-\lambda)0\in C$ 	for every $z\in\partial C$ and $\lambda\in[0,1]$. Conversely, if $\operatorname{rad}\partial C=C$ then $0\in\operatorname{rad}\partial C=C$.
\end{proof}

\begin{remark}
	Convexity cannot be removed from~(ii),
	as a compact set may contain $0$ in its
	interior and fail to be radial with respect to it,
	whereas $\operatorname{rad}\partial C$
	is radial by construction.
\end{remark}

We now specialize Proposition~\ref{minimalgauge}
to the case where $D$ is the closed $h$-ball of radius $r$
centred at the origin. In this case, the gauge of the radial hull
of the translated sphere $\partial D+W$ is a conic Minkowski norm
(see~\cite{JavSan11}) on the interior of its finiteness cone,
and it is a Minkowski norm on the whole $T_xM$ only when the origin
lies in the interior of the translated ball $D+W$.

\begin{lemma}\label{spheregauge}
	Let $W\in T_xM$, $r>0$ and
	$D:=\{v\in T_xM:\alpha(v)\leq r\}$,
	and let us set
	\[
		\mathcal C := D + W = \big\{v \in T_xM: \alpha(v - W) \le r	\big\},
		\quad\text{and}\quad
		K:=\operatorname{rad}(\partial \mathcal C).
	\]
	Then $K=\operatorname{co}\bigl(\{0\}\cup \mathcal C\bigr)$;
	in particular $K$ is compact and convex and,
	denoting by $F$ its gauge function,
	$K=\{y\in T_xM:F(y)\leq1\}$.
	Moreover, setting $\Delta:=\alpha^2(W)-r^2$, we have
	\begin{enumerate}
		\item[(i)] $F$ is finite exactly on $\Lambda:=\operatorname{cone}(\mathcal C)$,
			on which
			\begin{equation}\label{Fwind}
				F(y)=\frac{\alpha^2(y)}{h(W,y)+\sqrt{h(W,y)^2-\Delta\,\alpha^2(y)}},			\qquad y\in\Lambda\setminus\{0\};
			\end{equation}
		\item[(ii)] $\Lambda=T_xM$ if $\Delta<0$;
			$\Lambda=\{y\in T_xM:h(W,y)>0\}\cup\{0\}$ if $\Delta=0$;
			and
			$\Lambda=\{y\in T_xM:h(W,y)\geq\sqrt{\Delta}\,\alpha(y)\}$ if $\Delta>0$, in which case $\Lambda$ is the closed cone of half-aperture $\arcsin\bigl(r/\alpha(W)\bigr)$ about $W$;
		\item[(iii)] $0\in\operatorname{int}K$ if and only if $\Delta<0$;
			otherwise $0\in\partial K$.
			In all cases $F$ is smooth on $\operatorname{int}\Lambda\setminus\{0\}$ with strongly convex indicatrix there, hence it is a conic Minkowski norm on $\operatorname{int}\Lambda$.
	\end{enumerate}
\end{lemma}
\begin{proof}
By definition of $K$ and by
Lemma~\ref{radialhull}-(i),
$\operatorname{rad}\mathcal C=\operatorname{rad}\partial \mathcal C=K$.
Since $\mathcal C$ is convex and compact, 
$\operatorname{rad}\mathcal C=\operatorname{co}\bigl(\{0\}\cup \mathcal C\bigr)$,
the union of the segments joining $0$ to $\mathcal C$.
Thus $K=\operatorname{co}\bigl(\{0\}\cup \mathcal C\bigr)$ and it is compact,
radial, convex and Lemma~\ref{radialgauge}-(ii) gives $K=\{y\in T_xM:F(y)\leq1\}$.
	
Apply now Proposition~\ref{minimalgauge} with $D=\{v\in T_xM:\alpha(v)\leq r\}$,
whose gauge function is $G=\alpha/r$,
and with translation vector $W$:
then $\partial D+W=\partial \mathcal C$.
As $r>0$, we have $\partial D=\{\alpha=r\}=\{G=1\}$,
so that~\eqref{noflat} holds, and $0\in\operatorname{int}D$,
hence by Remark~\ref{degenerateray}
the exceptional case $y/F(y)=W$ does not occur.
Hence, $\Lambda=\operatorname{cone}(\mathcal C)$ and,
for every $y\in\Lambda\setminus\{0\}$,
$F(y)$ is the smallest positive solution of $\alpha(y-tW )=t r$.
Both sides being non-negative, this equation is equivalent to
\begin{equation}\label{raysphere}
\Delta\,t^2-2h(W,y)\,t+\alpha^2(y)=0 .
\end{equation}
If $\Delta=0$, then~\eqref{raysphere} is linear with the single solution
$t=\alpha^2(y)/\bigl(2h(W,y)\bigr)$,
which is positive if and only if $h(W,y)>0$.
If $\Delta\neq0$, the 	two roots of~\eqref{raysphere} have product
$\alpha^2(y)/\Delta$ and sum $2h(W,y)/\Delta$:
for	$\Delta<0$ the product is negative,
so there is exactly one positive root and $\Lambda=T_xM$.
For	$\Delta>0$ the roots are real and positive if and only if $h(W,y)\geq\sqrt{\Delta}\,\alpha(y)$.
Dividing by $\alpha(W)\alpha(y)>0$, this reads
$\cos\theta\geq\sqrt{\Delta}/\alpha(W )$,
where $\theta$ 	denotes the $h$-angle between $W$ and $y$;
as $\cos$ is decreasing on $[0,\pi]$, the
inequality 	describes the closed cone about $W$ of half-aperture
$\arccos\bigl(\sqrt{\Delta}/\alpha(W)\bigr)=\arcsin\bigl(r/\alpha(W)\bigr)$.
In all three cases the 	smallest positive solution is
$\bigl(h(W,y)-\sqrt{h(W,y)^2-\Delta\alpha^2(y)}\bigr)/\Delta$ when
$\Delta\neq0$, which upon 	rationalization becomes~\eqref{Fwind}, the same
expression giving the solution also when 	$\Delta=0$.
This proves (i) and (ii).
	
As for (iii), if $\Delta<0$ then 
$0\in\operatorname{int}\mathcal C\subseteq\operatorname{int}K$; if $\Delta\geq0$ then 
$\Lambda\neq T_xM$ by (ii),
whereas $0\in\operatorname{int}K$ would give
$\Lambda=\operatorname{cone}(K)=T_xM$, so $0\in\partial K$. 	
Let $y\in\operatorname{int}\Lambda\setminus\{0\}$. By (ii) the quantity
$h(W,y)^2-\Delta\,\alpha^2(y)$ is positive, so both the square root and the
denominator 	in~\eqref{Fwind} are smooth and nonvanishing at $y$ and $F$ is
smooth there.
Moreover, 	\eqref{gaugeboundary} gives
$y/F(y)\in\partial \mathcal C$, by Remark~\ref{degenerateray} the farthest
point in which the ray $\mathbb R_{>0}y$ meets $\partial \mathcal C$. Hence, the
indicatrix in $\operatorname{int}\Lambda$
is exactly the set $\{v\in\partial \mathcal C\cap\operatorname{int}\Lambda:\ \alpha^2(v)>\Delta\}$.
Indeed, if $v\in\partial \mathcal C\cap\operatorname{int}\Lambda$
one has $v\neq0$ and $t=1$ solves~\eqref{raysphere};
if $\Delta=0$ this is the only root and $\alpha^2(v)>0=\Delta$,
while if $\Delta\neq0$ the other root is
$\alpha^2(v)/\Delta$, distinct from $1$
because the discriminant is positive.
If $\Delta < 0$ this root is negative,
so $t = 1$ is the smallest positive root and 
$\alpha^2(v) > \Delta$ holds trivially;
if $\Delta > 0$ it is positive,
so $t = 1$ is the smallest positive root 
if and only if $\alpha^2(v)/\Delta > 1$,
i.e., $\alpha^2(v) > \Delta$.
\end{proof}	
\begin{remark}
When $\Delta>0$ the tangency locus $\partial \mathcal C\cap\partial\Lambda$ is
exactly 	$\{v\in\partial \mathcal C:\alpha^2(v)=\Delta\}$, so $\partial K$
splits as the disjoint union 	of the indicatrix
$\{v\in\partial \mathcal C:\alpha^2(v)>\Delta\}$ 
and 	$\operatorname{rad}\bigl(\partial \mathcal C\cap\partial\Lambda\bigr)$,
the union of the segments 	joining $0$ to the tangency locus.
Along the latter $K$ is convex but not strongly convex,
being ruled by those segments: this is why $\operatorname{int}\Lambda$ cannot
be enlarged in 	Lemma~\ref{spheregauge}-(iii).
When $\Delta=0$ the situation is different:
there $K=\mathcal C$ is a round ball and $\partial K$ is strongly
convex everywhere,
the only obstruction to $F$ being a Minkowski norm on all
of $T_xM$ being that $0\in\partial K$,
whence $\Lambda$ is a half-space.
\end{remark}

Finally, we recall that the \emph{support function} of a nonempty compact
set $S\subset T_xM$ is the positively homogeneous function
$p\in T_x^*M\mapsto\max_{z\in S}p(z)$, and that the support functions of $S$
and of $\operatorname{co}S$ coincide.
We record, in the form in which we shall use it, the following
elementary separation statement.
\begin{lemma}\label{separation}
Let $S\subset T_xM$ be nonempty and compact. Then:
\begin{enumerate}
\item[(i)] $\max_{z\in S}p(z)\geq0$ for every $p\in T^*_xM$ if and only if $0\in\operatorname{co}S$;
\item[(ii)] $\max_{z\in S}p(z)>0$ for every $p\in T^*_xM\setminus\{0\}$ if and only if $0\in\operatorname{int}\operatorname{co}S$.
\end{enumerate}
\end{lemma}
\begin{proof}
Since $\max_{z\in S}p(z)=\max_{z\in\operatorname{co}S}p(z)$, we may assume that $S$ is convex and compact.

(i) If $0\in S$, then $\max_{z\in S}p(z)\geq p(0)=0$ for every $p$. If $0\notin S$, then, $S$ being compact and convex, the point $0$ can be strictly separated from $S$: there is $p$ with $p(z)<0=p(0)$ for every $z\in S$, so that $\max_{z\in S}p(z)<0$.

(ii) If $0\in\operatorname{int}S$ and $p\neq0$, choose $\xi$ with $p(\xi)>0$ and $\varepsilon>0$ with $\varepsilon\xi\in S$; then $\max_{z\in S}p(z)\geq\varepsilon p(\xi)>0$. If $0\notin\operatorname{int}S$, then $0$ lies outside $S$ or on its boundary, so there is a supporting hyperplane, that is a $p\neq0$ with $p(z)\leq0$ for every $z\in S$, and $\max_{z\in S}p(z)\leq0$.
\end{proof}

\section{A general model of motion under drift fields}\label{model}

In this section, we derive the velocity law~\eqref{eq:generalForm} from a
kinematic model with finitely many drift fields, and we then study the
corresponding admissible velocity sets in the regimes of positive and
vanishing self-propulsion.

\begin{figure}[t]
	\centering
	\begin{tikzpicture}[scale=1.5]

		\def\a{2}
		\def\b{1}
		\def\angle{30}
		
		\def\etaC{1}
		\def\m{1.6}
		\def\s{1.4}
		
		\draw[->,very thin,gray!60] (-2.65,0) -- (3.65,0);
		\draw[->,very thin,gray!60] (0,-1.40) -- (0,2.2);
		
		\def\Vx{1}
		\def\Vy{0.2}
		\coordinate (V) at (\Vx,\Vy);

		\def\uangle{25}
		
		\pgfmathsetmacro{\upx}{\a*cos(\uangle)}
		\pgfmathsetmacro{\upy}{\b*sin(\uangle)}
		
		\pgfmathsetmacro{\ux}{
			\upx*cos(\angle)-\upy*sin(\angle)}
		\pgfmathsetmacro{\uy}{
			\upx*sin(\angle)+\upy*cos(\angle)}
		\coordinate (u) at (\ux,\uy);

		\pgfmathsetmacro{\Vpx}{
			\Vx*cos(\angle)+\Vy*sin(\angle)}
		\pgfmathsetmacro{\Vpy}{
			-\Vx*sin(\angle)+\Vy*cos(\angle)}
		
		\pgfmathsetmacro{\hVu}{
			\Vpx*\upx/(\a*\a)
			+
			\Vpy*\upy/(\b*\b)}
		
		\pgfmathsetmacro{\Vparallelx}{\hVu*\ux}
		\pgfmathsetmacro{\Vparallely}{\hVu*\uy}
		\coordinate (Vparallel) at
		(\Vparallelx,\Vparallely);
		
		\pgfmathsetmacro{\Vperpx}{\Vx-\Vparallelx}
		\pgfmathsetmacro{\Vperpy}{\Vy-\Vparallely}
		\coordinate (Vperp) at
		(\Vperpx,\Vperpy);

		\pgfmathsetmacro{\muVparallelx}{\m*\Vparallelx}
		\pgfmathsetmacro{\muVparallely}{\m*\Vparallely}
		
		\pgfmathsetmacro{\sigmaVperpx}{\s*\Vperpx}
		\pgfmathsetmacro{\sigmaVperpy}{\s*\Vperpy}
		
		\pgfmathsetmacro{\Yux}{
			\etaC*\ux+\muVparallelx+\sigmaVperpx}
		\pgfmathsetmacro{\Yuy}{
			\etaC*\uy+\muVparallely+\sigmaVperpy}
		\coordinate (Yu) at (\Yux,\Yuy);

		\draw[rotate=\angle,blue,thick]
		(0,0) ellipse (\a cm and \b cm);

		\draw[black,dashed]
		(0,0) -- (u);
		
		\filldraw[black] (u) circle (1pt);
		\node[above right] at (u) {$u$};

		\draw[->,black,thick]
		(0,0) -- (V)
		node[pos=1,right] {$V$};
		
		\draw[gray!55,dashed]
		(0,0) -- (Vparallel) -- (V) -- (Vperp) -- cycle;
		
		\node[gray!70,left, xshift=9.5, yshift=7.5]
		at ({0.52*\Vparallelx},
		{0.52*\Vparallely})
		{$V^\parallel$};
		
		\node[gray!70,below]
		at ({0.52*\Vperpx},
		{0.52*\Vperpy})
		{$V^\perp$};

		\draw[->,black,thick]
		(0,0) --
		(\muVparallelx,\muVparallely)
		node[pos=1.3,left]
		{$\mu V^\parallel$};
		
		\draw[->,black,thick]
		(0,0) --
		(\sigmaVperpx,\sigmaVperpy)
		node[below, xshift=0, yshift=1]
		{$\sigma V^\perp$};

		\draw[->,black,very thick]
		(0,0) -- (Yu)
		node[pos=1,above right] {$\mathcal Y(u)$};
		
	\end{tikzpicture}
	
\caption{The velocity $\mathcal Y(u)$ when there is only one  field $V$.
		The ellipse represents the $h_x$-unit sphere $UT_xM$. For $u\in UT_xM$, the vector $V$ is decomposed as
		$V=V^\parallel+V^\perp$, where
		$V^\parallel=h_x(V,u)u$ and $V^\perp$ is $h_x$-orthogonal to $u$.
		The velocity is
		$\mathcal Y(u)=\eta u+\mu V^\parallel+\sigma V^\perp$.
		In the plot, $\eta=1$, $\mu=1.6$, and $\sigma=1.4$.}
\label{fig:constructionYu}
\end{figure}
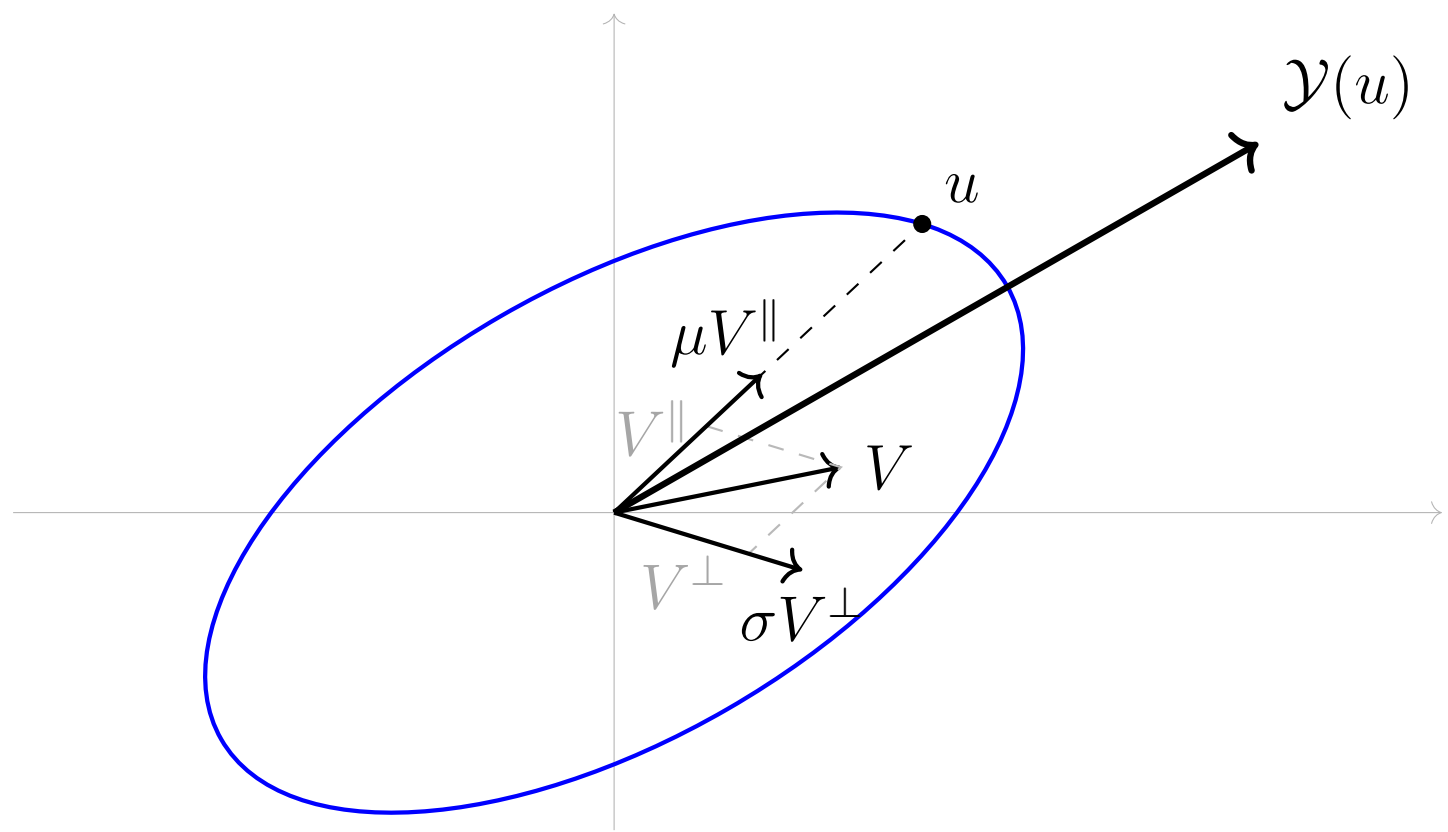

We start with a single smooth vector field $V \in \mathfrak X(M)$
that acts through smooth coefficients $\mu,\sigma\in C^\infty(M)$
which multiply respectively its components parallel
and orthogonal to the chosen control direction.
In other words, at any given point $x \in M$
and for a given vector $u \in UT_xM$,
$V$ is decomposed and deformed along $u$ and its component
$h$-orthogonal to $u$ (see Figure~\ref{fig:constructionYu}):
\[
\mu(x)\, h(V(x),u) u, 	\qquad \text{and} \qquad \sigma(x)\, (V(x) - h(V(x),u)u).
\]
From now on, for the sake of simplicity, we remove the dependence on the point $x$ when no confusion may arise.
Assuming that these two contributions add to the self-propelled
velocity $\eta u$,
the resulting velocity with respect to the control direction
$u$ is given by
\begin{equation}\label{eq:Y-oneV}
	\mathcal{Y}(u)
	= \big(\eta + \mu\, h(V,u)\big) u	
	+ \sigma \big(V - h(V,u)u \big),
\end{equation}
where $\eta\colon M\to[0,+\infty)$ is a smooth function
representing the maximal speed that the system can achieve
without any external contribution.
Defining the two smooth vector fields
\begin{equation}\label{eq:hat-fields-one}
	\widehat Z \coloneqq (\mu-\sigma)V,	\qquad	\widehat W \coloneqq \sigma V,
\end{equation}
the velocity law~\eqref{eq:Y-oneV}
takes the form~\eqref{eq:generalForm}.

For \emph{any finite number} of fields acting on the system,
their actions can always be encoded in two vector fields 
$\widehat Z,\widehat W\in\mathfrak X(M)$ as follows.
Let $(V_i)_{i=1,\ldots,N}$ be a finite family
of smooth vector fields on $M$,
and let $\mu_i,\sigma_i\in C^\infty(M)$
be their respective tangential and orthogonal coefficients
with respect to the control direction $u$.
The resulting velocity is
\[
	\mathcal{Y}(u)=	\bigg(\eta + \sum_{i = 1}^N \mu_i\, h(V_i,u)\bigg)u	+ \sum_{i = 1}^{N}
	\bigg(\sigma_i \Big(V_i - h(V_i,u)u\Big)\bigg).
\]
A simple algebraic computation shows that,
by setting
\begin{equation}\label{eq:hat-fields}
	\widehat Z:=\sum_{i=1}^{N}(\mu_i-\sigma_i)V_i,	\qquad	\widehat W:=\sum_{i=1}^{N}\sigma_iV_i,
\end{equation}
the resulting velocity is again given by~\eqref{eq:generalForm}.

\begin{remark}
The main difference between one drift field and two or more drift fields is that, in the former case, $\widehat Z$ and $\widehat W$ are collinear, since they are both multiples of the original field $V$, whereas in the latter case they need not have the same direction. 
\end{remark}

\subsection{Positive self-propulsion}
\label{subsec:positive-eta}
Let $M_+:=\{x\in M:\eta(x)>0\}$.
On $M_+$ we may define the normalized deformation and translation fields
\begin{equation}\label{eq:ZW-formula}
Z:=\frac{\widehat Z}{\eta}
=\frac1\eta\sum_{i=1}^{N}(\mu_i-\sigma_i)V_i,
\qquad
W:=\frac{\widehat W}{\eta}
=\frac1\eta\sum_{i=1}^{N}\sigma_iV_i,
\end{equation}
so that the fundamental law~\eqref{eq:generalForm}
is equivalently written on $M_+$ as
\begin{equation}\label{eq:normalizedForm}
	\mathcal Y(u)=\eta\Bigl(\bigl(1+h(Z,u)\bigr)u+W\Bigr).
\end{equation}
As the discussion that follows is pointwise,
let us fix $x \in M_+$ and omit the dependence on it.
Consider the term
\[
	\mathcal{X}(u) \coloneqq \big(1 + h(Z,u)\big) u,
\]
so that the subset of $T_xM$ generated by $\mathcal{Y}/\eta$
is a translation of the set generated by $\mathcal{X}$,
namely
\[
	\Big\{\frac{1}{\eta}\mathcal{Y}(u): u \in UT_xM\Big\}
	= \big\{\mathcal{X}(u): u \in UT_xM\big\} + W.
\]
Consider the radial set
\begin{equation*}
\mathcal{M} = \operatorname{rad}\big(\mathcal X(UT_xM)\big),
\end{equation*}
and let $F_{\mathcal{M}}\colon T_xM \to [0,+\infty]$
be its gauge function, namely
\begin{equation}
	\label{eq:def-FM}
	F_{\mathcal{M}}(y)
	:= \inf\big\{
		\lambda > 0:
		y \in \lambda \mathcal{M}
	\big\}.
\end{equation}
Let us define
\begin{equation}
	\label{brho}
	b:=\alpha(Z),
	\qquad\rho(u):=1+h(Z,u),
	\qquad u\in UT_xM.
\end{equation}
Since $|h(Z,u)|\leq b$,
we have $\rho\geq1-b$ on $UT_xM$;
this implies that if $b<1$,
then $\rho>0$, so $\mathcal X(u)=\rho(u)u$
is a radial graph over $UT_xM$ and
\[
	\mathcal M=\bigl\{su:u\in UT_xM,\ 0\leq s\leq\rho(u)\bigr\},
\]
so that
\begin{equation}\label{radialgraph}
	\partial\mathcal M
	=\bigl\{\rho(u)u:u\in UT_xM\bigr\}=\mathcal X(UT_xM),
\end{equation}
and every ray issuing from $0$
meets $\partial\mathcal M$ exactly once.
In particular
$\{y:\alpha(y)\leq1-b\}\subseteq\mathcal M$,
so $0\in\operatorname{int}\mathcal M$.
If instead $b\geq1$, then $\mathcal X(-Z/b)=(b-1)Z/b$ and $\mathcal X(Z/b)=(1+b)Z/b$ lie on the closed ray of $Z$. Since $\mathcal X(u)$ is a multiple of $u$, these are the only points of $\mathcal X(UT_xM)$ on the line spanned by $Z$; hence $\mathcal M$ meets the ray issuing from $0$ in the direction of $-Z$ only at $0$.

Thus, $0\in\operatorname{int}\mathcal M$
if and only if $b<1$, and in this case,
from Lemma~\ref{radialgauge},
we have
$\Lambda=\operatorname{cone}(\mathcal M)=T_xM$;
moreover for $y\in T_xM\setminus\{0\}$,
set $\widehat y:=\frac{y}{\alpha(y)}$, then
\[
y\in\lambda\mathcal M 	\quad\Longleftrightarrow\quad	\alpha(y)	\leq\lambda \rho(\widehat y)
	\quad\Longleftrightarrow\quad	\lambda\geq	\frac{\alpha(y)}{\rho(\widehat y)}.
\]
Consequently,
\begin{equation}	
	\label{eq:almostalphabeta}
F_{\mathcal M}(y)	=	\frac{\alpha(y)}	{1+h(Z,\widehat y)}	=	\frac{\alpha^2(y)}	{\alpha(y)+h(Z,y)},
\end{equation}
and, by homogeneity, we have $F_{\mathcal M}(0) = 0$. The condition $b <1$
guarantees that $F_{\mathcal{M}}$ is well-defined and non-negative on $T_xM$.

In the  two-dimensional case,  the indicatrix of $F_{\mathcal M}$ is a   lima\c{c}on (see \cite{Matsumoto1989}). By  using a general result about strong convexity of $(\alpha,\beta)$-metrics (cf.~\cite[Lemma~1.1.2]{ChernShen2005}),  it is not difficult to see that $F_{\mathcal M}$ has  strongly convex indicatrix if and only if $b<1/2$. 
 We include a short direct proof for the reader's convenience. 
 The argument distinguishes the condition $b<1$, which guarantees the fact that~\eqref{eq:almostalphabeta} is well-defined  and non-negative on $TM_{+}$ from the sharp condition $b<1/2$ for strong convexity, see Figure~\ref{fig:confronto_z}. Moreover, it makes explicit
the degeneracy occurring when  $b=1/2$.
\begin{figure}[t]
	\centering
	\pgfmathsetmacro{\globalScale}{0.5} 
	\begin{subfigure}{0.3\textwidth}
		\centering
		\pgfmathsetmacro{\Aval}{0.2}
		\pgfmathsetmacro{\Bval}{0.2}
		\pgfmathsetmacro{\zx}{0.45}
		\pgfmathsetmacro{\zy}{0.9}
		\pgfmathsetmacro{\bA}{sqrt(\Aval*\zx*\zx + \Bval*\zy*\zy)}
		\begin{tikzpicture}[scale=\globalScale]
			\pgfmathsetmacro{\sqrtA}{sqrt(\Aval)}
			\pgfmathsetmacro{\sqrtB}{sqrt(\Bval)}
			\draw[->, gray] (-3,0) -- (3,0);
			\draw[->, gray] (0,-2) -- (0,3);
			\draw[blue, thick, domain=0:360, samples=200, variable=\t]
			plot ({cos(\t)/\sqrtA}, {sin(\t)/\sqrtB});
			\draw[red, thick, domain=0:360, samples=400, variable=\t]
			plot ({(\sqrtA*\zx*cos(\t) + \sqrtB*\zy*sin(\t) + 1)*cos(\t)/\sqrtA},
			{(\sqrtA*\zx*cos(\t) + \sqrtB*\zy*sin(\t) + 1)*sin(\t)/\sqrtB});
			\draw[->, thick] (0,0) -- (\zx, \zy) node[right]{$Z$};
		\end{tikzpicture}
		\caption{$b\approx
			\protect\pgfmathprintnumber[precision=2]{\bA}$;
			$\partial\mathcal{M}$ is strongly convex.}
	\end{subfigure}
	\hfill
	\begin{subfigure}{0.3\textwidth}
		\centering
		\pgfmathsetmacro{\Aval}{0.2}
		\pgfmathsetmacro{\Bval}{0.2}
		\pgfmathsetmacro{\zx}{0.8}
		\pgfmathsetmacro{\zy}{1.6}
		\pgfmathsetmacro{\bB}{sqrt(\Aval*\zx*\zx + \Bval*\zy*\zy)}
		\begin{tikzpicture}[scale=\globalScale]
			\pgfmathsetmacro{\sqrtA}{sqrt(\Aval)}
			\pgfmathsetmacro{\sqrtB}{sqrt(\Bval)}
			\draw[->, gray] (-3,0) -- (3,0);
			\draw[->, gray] (0,-2) -- (0,3);
			\draw[blue, thick, domain=0:360, samples=200, variable=\t]
			plot ({cos(\t)/\sqrtA}, {sin(\t)/\sqrtB});
			\draw[red, thick, domain=0:360, samples=400, variable=\t]
			plot ({(\sqrtA*\zx*cos(\t) + \sqrtB*\zy*sin(\t) + 1)*cos(\t)/\sqrtA},
			{(\sqrtA*\zx*cos(\t) + \sqrtB*\zy*sin(\t) + 1)*sin(\t)/\sqrtB});
			\draw[->, thick] (0,0) -- (\zx, \zy) node[below right]{$Z$};
		\end{tikzpicture}
		\caption{$b\approx
			\protect\pgfmathprintnumber[precision=2]{\bB}
			$; $\partial\mathcal{M}$ is 
			not convex.}
	\end{subfigure}
	\hfill
	\begin{subfigure}{0.3\textwidth}
		\centering
		\pgfmathsetmacro{\Aval}{0.2}
		\pgfmathsetmacro{\Bval}{0.2}
		\pgfmathsetmacro{\zx}{1.6}
		\pgfmathsetmacro{\zy}{3.2}
		\pgfmathsetmacro{\bC}{sqrt(\Aval*\zx*\zx + \Bval*\zy*\zy)}
		\begin{tikzpicture}[scale=\globalScale]
			\pgfmathsetmacro{\sqrtA}{sqrt(\Aval)}
			\pgfmathsetmacro{\sqrtB}{sqrt(\Bval)}
			\draw[->, gray] (-3,0) -- (3,0);
			\draw[->, gray] (0,-2) -- (0,3);
			\draw[blue, thick, domain=0:360, samples=200, variable=\t]
			plot ({cos(\t)/\sqrtA}, {sin(\t)/\sqrtB});
			\draw[red, thick, domain=0:360, samples=400, variable=\t]
			plot ({(\sqrtA*\zx*cos(\t) + \sqrtB*\zy*sin(\t) + 1)*cos(\t)/\sqrtA},
			{(\sqrtA*\zx*cos(\t) + \sqrtB*\zy*sin(\t) + 1)*sin(\t)/\sqrtB});
			\draw[->, thick] (0,0) -- (\zx, \zy) node[right]{$Z$};
		\end{tikzpicture}
		\caption{$b\approx
			\protect\pgfmathprintnumber[precision=2]{\bC}$; $\partial \mathcal{M}$ is not equal to  $\mathcal X(UT_xM)$.}
	\end{subfigure}
	
\caption{Comparison between the set $\mathcal{X}(UT_xM)$ (red) and the  unit sphere  (blue) in $h$-orthogonal coordinates in $T_xM$ for different  $b$.
The vector $Z$ whose norm is $b$ is		represented as a black arrow.}
\label{fig:confronto_z}
\end{figure}
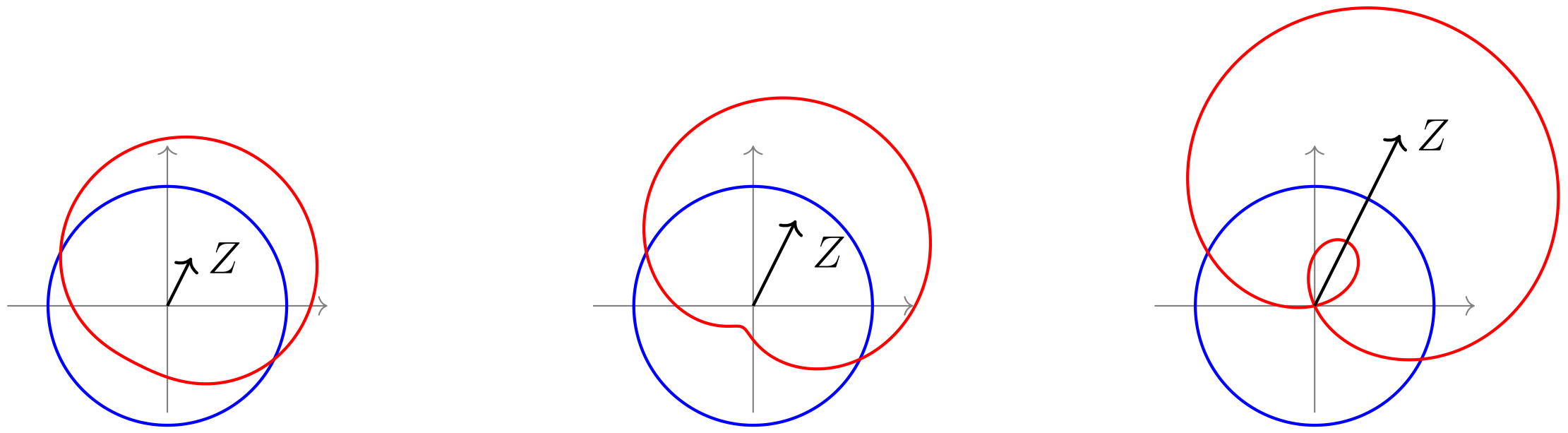
\begin{proposition}
	\label{prop:FinslerM-conditions}
	The gauge function $F_{\mathcal M}\colon T_xM \to [0,+\infty]$
	is a Minkowski norm if and only if $b <\frac{1}{2}$.
\end{proposition}
\begin{proof}
	If $b\geq1$, then, as observed above, $\mathcal M$ misses an entire ray
	issuing from the origin, so that $F_{\mathcal M}$ takes the value
	$+\infty$ and is not a Minkowski norm. Hence, we can assume that $b<1$.
	After choosing $h$-orthonormal linear coordinates
	on	$T_xM$, we can identify $(T_xM,h)$
	with the Euclidean space $(\mathbb R^n,\langle\cdot,\cdot\rangle)$.
	We write
	\[
		|y|:=\sqrt{\langle y,y\rangle},	\qquad	z:=Z,	\qquad	b:=|z|.
	\]
	In these coordinates,
	$\partial\mathcal M	=\{	\bigl(1+\langle z,u\rangle\bigr)u:	|u|=1	\}$.
	Let us consider the function
	\[
		G\colon\mathbb R^n\setminus\{0\}\to \mathbb R,
		\qquad	G(y):=|y|^2-|y|-\langle z,y\rangle.
	\]
	Writing $y=ru$, with $r>0$ and $|u|=1$,
	we have $G(ru)	=	r\big(r-\rho(u)\big)$,
	where $\rho(u):= 1+\langle z,u\rangle$.
	Thus, $\partial\mathcal M	=	G^{-1}(0)$
	while  $G<0$ on the interior side of $\partial\mathcal M$.
	At	$y=\rho(u)u\in\partial\mathcal M$,
	we have 
	\begin{equation}\label{nablaG}
		\nabla G(y)	=	2y-\frac{y}{|y|}-z	=	(2\rho(u)-1)u-z.
	\end{equation}
	Moreover, $	\langle\nabla G(y),u\rangle=2\rho(u)-1-\langle z,u\rangle=\rho(u)>0$.
	Thus, $\nabla G(y)\neq0$ and points toward the exterior of	$\mathcal M$.
	In particular,  $\partial\mathcal M$ is a smooth hypersurface in $T_x M$.

	Let then $\nu(y):=\nabla G(y)/|\nabla G(y)|$	be the outward unit normal. 
	For $\xi\in T_y\partial\mathcal M$, differentiating the normalized gradient gives
	\[
		D_\xi\nu=
		\frac{D^2G_y[\xi]}{|\nabla G(y)|}-\frac{\left\langle D^2G_y[\xi],\nabla G(y)\right\rangle}{	|\nabla G(y)|^3
		}\nabla G(y),
	\]
	where $D^2G_y[\xi]$ denotes the vector obtained by applying the
	Euclidean Hessian of $G$ to $\xi$.
	Consequently,
	\[
		\operatorname{II}_y(\xi,\xi)
		:=	\langle D_\xi\nu,\xi\rangle=
		\frac{D^2G_y(\xi,\xi)}{|\nabla G(y)|}
		-	\frac{\left\langle	D^2G_y[\xi],\nabla G(y)	\right\rangle
		\left\langle\nabla G(y),\xi\right\rangle}{	|\nabla G(y)|^3	}.
	\]
	Since $\partial\mathcal M=G^{-1}(0)$ and $\xi\in T_y\partial\mathcal M$,
	we have $\left\langle\nabla G(y),\xi\right\rangle=0$.
	Therefore,
	\[
		\operatorname{II}_y(\xi,\xi)=\frac{D^2G_y(\xi,\xi)}{|\nabla G(y)|},\qquad\xi\in T_y\partial\mathcal M.
	\]
	From \eqref{nablaG}, at $y=\rho(u)u$ we get
	\begin{align*}
		D^2G_y(\xi,\xi)	&=2|\xi|^2	-	\frac{|\xi|^2-\langle u,\xi\rangle^2}{\rho(u)}=\left(2-\frac{1}{\rho(u)}\right)|\xi|^2	
	+\frac{1}{\rho(u)}\langle u,\xi\rangle^2.\end{align*}
	If $b<1/2$, then $\rho(u)=1+\langle z,u\rangle	\geq1-b>\frac12$	and hence
	$D^2G_y(\xi,\xi)>0$ for every $\xi\neq0$.	In particular, $\operatorname{II}_y$ is positive definite on
	$T_y\partial\mathcal M$, i.e. $\partial\mathcal M$ is strongly convex.

	Conversely, assume that $1/2\leq b<1$.
	Let us set $u_0:=-z/b$,	$r_0:=1-b$ and $y_0:=r_0u_0$.
	Then $y_0\in\partial\mathcal M$ and 	$\nabla G(y_0)=	(2r_0-1)u_0-z=r_0u_0$.
	Since $\dim M\geq2$,
	we can choose a nonzero	$\xi\in u_0^\perp=T_{y_0}\partial\mathcal M$.
	We then have
	$D^2G_{y_0}(\xi,\xi)=		\left(2-\frac1{r_0}\right)|\xi|^2=	\frac{1-2b}{1-b}|\xi|^2$.
	This quantity vanishes when $b=1/2$ and is negative when $b>1/2$.
	Hence, $\partial\mathcal M$ is not strongly convex. 

	Summing up, $\mathcal M$ is a radial, compact set,
	containing $0\in T_x M$ in its interior if and only if $b<1$;
	moreover   $\partial\mathcal M$ is a smooth, strongly convex,
	hypersurface in $T_xM$ if and only if $b<\frac12$.
	This is equivalent to its gauge function 	$F_{\mathcal M}$
	being a real-valued,
	non-negative, positively homogeneous  continuous function on $T_x M$,
	smooth on $T_x M\setminus\{0\}$, such that $D_y^2\left(	\frac12F_{\mathcal M}^2	\right)(y)$ is  positive definite for every $y\neq0$
	(see, e.g. \cite[Prop. 2.3-(v)]{JavSan11}).
	Therefore, $F_{\mathcal{M}}$ is a Minkowski norm if and only if $b < \frac12$.
\end{proof}
\begin{remark}[Pure pursuit navigation]
In~\cite{BidRaf09}, Bidabad and Rafie-Rad consider, in connection with pure
pursuit navigation on a Riemannian manifold, the velocity law $u\mapsto u+W_T$,
where a pursuer steering along $u\in UT_xM$ tracks a target
whose velocity $w$ acts only through its projection on the line of sight,
so that $W_T=h(u,w)u$.
The resulting velocity $\bigl(1+h(w,u)\bigr)u$ is~\eqref{eq:normalizedForm}
with $\eta=1$, $Z=w$ and $W=0$,
and the corresponding gauge is the Matsumoto metric $F_{\mathcal M}$.
The results in~\cite{BidRaf09} are stated under the assumption $b<1$.
However, by Proposition~\ref{prop:FinslerM-conditions},
$F_{\mathcal M}$ is a Minkowski norm if and only if $b<1/2$.
For $b=1/2$ the set $\mathcal M$ is still convex, being a limit of convex sets,
but it is not strongly convex; for $1/2<b<1$ it is not convex,
so that Proposition~\ref{nonfinsler-existence} does not apply
and time-optimal paths need not be geodesics of a Finsler metric.
Nevertheless, $0\in\operatorname{int}\mathcal M$ in the whole range $b<1$,
and Theorems~\ref{maximalspeed} and~\ref{convexification} still apply
(see Section~\ref{sec:time-travel}).
\end{remark}

Let us now consider also the translation 
given by the vector field $W$,
so for a given point $x \in M_+$
we study the admissible velocity set 
\[
	\mathcal{B} :=
\operatorname{rad}\big(\mathcal Y(UT_xM)\big),
\]
where $\mathcal{Y}$ is given as in~\eqref{eq:normalizedForm},
and its gauge function
\[
	F\colon T_xM \to [0,+\infty],
	\qquad F(y) :=
	\inf\big\{\lambda > 0: y \in \lambda\mathcal{B}\big\}.
\]

\begin{proposition}
	\label{generalcase}
	Assume that $b<1/2$. 
	Then the gauge function $F$ of $\mathcal B \subset T_xM$
	is finite valued on $T_xM$ if and only if	
	\begin{equation} \label{translation}
		F_{\mathcal M}(-W) < 1.
	\end{equation}
	If the latter condition holds,
	then $F\colon T_xM \to [0,+\infty)$
	is a Minkowski norm.
\end{proposition}

\begin{proof}
	Since $b<1/2$,
	by Proposition~\ref{prop:FinslerM-conditions}
	$F_{\mathcal M}$ is a Minkowski norm on $T_xM$
	and $\mathcal M$ is its smooth,
	strongly convex, closed unit ball in $T_x M$. 
	Since~\eqref{translation} is equivalent to
	$-W\in\operatorname{int}\mathcal M$,
	the translated strongly convex set $\mathcal M+W$
	contains $0$ in its interior.
	By Lemma~\ref{radialhull}-(ii),
	it therefore coincides with the radial hull of its  boundary.
	As
	\[
		\partial(\mathcal M+W)=\partial\mathcal M+W=\mathcal X(UT_xM)+W,
	\]
	and  $\mathcal Y(u)=\eta \big(\mathcal X(u)+W\big)$, from \eqref{admissible}, we then get
	\begin{equation}\label{ballidentity}
		\mathcal B=\eta\bigl(\mathcal M+W\bigr).
	\end{equation}
	Thus, $\mathcal B$ is a compact, smooth, strongly convex domain containing the	origin in its interior.
	Consequently, its gauge	$F$ is finite and convex on $T_xM$ and, in fact, is a Minkowski norm.

	Conversely, suppose that $F$ is finite-valued on $T_xM$,
	that is $\Lambda=T_xM$,
	and set $\mathcal{S}:=\mathcal Y(UT_xM)=\eta\bigl(\partial\mathcal M+W\bigr)$,
	where the second equality holds by~\eqref{radialgraph}.
	By Lemma~\ref{radialgauge},
	$T_xM=\Lambda=\operatorname{cone}(\mathcal B)=\operatorname{cone}(\mathcal{S})$,
	so that every nonzero covector $p \in T^*_xM$
	satisfies $\max_{z\in \mathcal{S}}p(z)>0$:
	it suffices to pick $y$ with $p(y)>0$
	and $t>0$ with $ty\in \mathcal{S}$.
	Lemma~\ref{separation}-(ii) then gives
	$0\in\operatorname{int}\operatorname{co}\mathcal{S}$;
	since $n\geq2$, the convex hull of the boundary
	of the convex body $\mathcal M+W$ is the body itself,
	so that  $0\in\operatorname{int}(\mathcal M+W)$,
	which is~\eqref{translation}.
\end{proof}

\begin{remark}\label{identityremark}
If $b<1/2$,
by Lemma~\ref{radialhull}-(ii)
we have that \eqref{ballidentity} holds
if and only if
\[
	0\in\mathcal M+W,
	\qquad\text{that is,}\qquad
	F_{\mathcal M}(-W)\leq1 .
\]
In particular,~\eqref{ballidentity} survives in
the critical case $F_{\mathcal M}(-W)=1$,
where $0\in\partial\mathcal B$.
The corresponding gauge is then finite only on
the cone generated by $\mathcal B$
and defines a singular conic Finsler structure
(in the sense of \cite{JavSan11})
on the interior of that cone;
in particular this metric is Kropina
when we have also $Z=0$
(see Subsection~\ref{subsec:zermelo}).
We stress that~\eqref{admissible},
and not~\eqref{ballidentity},
is the definition of the admissible velocity set:
outside the present hypotheses the two sets differ,
even though the infimum travel time can be the same (see Section~\ref{sec:time-travel}).
\end{remark}

\begin{remark}
	Fix $x\in M_+$ with $b<1$ and apply Proposition~\ref{minimalgauge}
    with $D:=\eta\mathcal M$ and translation vector $\widehat W$.
	The gauge function of $D$ is
	$G(v)=F_{\mathcal M}(v)/\eta$,
	and~\eqref{radialgraph} gives
	\[
		\partial D=\bigl\{\eta\rho(u)u:u\in UT_xM\bigr\},
		\qquad	\partial D+\widehat W=\mathcal Y(UT_xM),
	\]
	so that $K=\mathcal B$.
	Moreover, for every $u\in UT_xM$,
	\[
		G\bigl(\eta\rho(u)u\bigr)=F_{\mathcal M}\bigl(\rho(u)u\bigr)=1.
	\]
	Thus~\eqref{noflat} is satisfied.
	Notice also that 	$0\in\operatorname{int}D$,
	so the exceptional equality
	$y/F(y)=\widehat W$ cannot occur.
	Proposition~\ref{minimalgauge} therefore yields,
	for every	$y\in\Lambda\setminus\{0\}$,
	\[
		F(y)	=	\min\bigl\{	s>0:G(y-s\widehat W)=s\bigr\}.
	\]
	Using $\widehat W=\eta W$
	and making the change of variable	$t=\eta s$,
	this becomes
	\[
	\eta F(y)	=	\min\bigl\{	t>0:F_{\mathcal M}(y-tW)=t\bigr\}.
\]
Under the assumptions of Proposition~\ref{generalcase},
$\Lambda=T_xM$ and $\mathcal M+W$ is a strongly convex body 	containing $0$ in its interior.
Hence every ray from $0$ meets 	$\partial\mathcal M+W$ exactly once,
and the above positive 	solution is unique.
In particular,
\begin{equation}\label{eqforF}
	F_{\mathcal M}\bigl(y-\eta F(y)W\bigr) 	= 	\eta F(y)
\end{equation}
for every $y\in T_xM\setminus\{0\}$.
Equivalently, $F(y)$ is the	unique positive solution
$s$ of
\[
	F_{\mathcal M}\left(\frac{y}{\eta\, s}-W	\right)=1.
\]
\end{remark}
Let us now consider also the dependence on $x \in M_+$.
Let $\beta$ be the one-form on $M_+$ defined by $\beta_x(y) := - h(Z_x,y)$.
At every $x\in M_+$ with $b_x<1$, using~\eqref{eq:almostalphabeta},
the function $F_{\mathcal M}(x,\cdot)$ is finite and can be written as
an {\em $(\alpha,\beta)$--metric}
(in particular a {\em Matsumoto} one \cite[Def. 5]{Matsumoto1992}),
namely
\begin{equation}\label{Matsumoto}
	F_{\mathcal M}(x,y) = \frac{\alpha_x^2(y)}{\alpha_x(y)-\beta_x(y)},
	\qquad y\in T_xM\setminus\{0_x\},
\end{equation}
where $\alpha_x(y)-\beta_x(y)\geq(1-b_x)\,\alpha_x(y)>0$.
Since $b$ is continuous, the set of such points is open in $M_+$,
and the last formula shows that $F_{\mathcal M}$ is smooth on
$\{(x,y)\in TM_+: b_x<1,\ y\neq0_x\}$.
At the points where $b_x\geq1$, instead, $F_{\mathcal M}(x,\cdot)$
takes the value $+\infty$ in the direction of $-Z_x$, as observed above.
By applying Proposition~\ref{prop:FinslerM-conditions} and Proposition~\ref{generalcase},
we thus have the following two results.
\begin{corollary}
	\label{cor:FM-Finsler}
	The function $F_{\mathcal{M}}\colon TM_+\to [0,+\infty]$
	is a Finsler metric on $M_+$ if and only if 
	$b_x = \alpha_x(Z_x) < 1/2$ for every $x \in M_+$.
\end{corollary}

\begin{corollary}
	\label{cor:F-Finsler}
	The function
	$F\colon TM_+ \to [0,+\infty]$
	is a Finsler metric on $M_+$ if and only if
	\begin{equation}
		\label{eq:FinslerM+conditions}
		b_x = \alpha_x(Z_x) < \frac{1}{2},
		\quad\text{and}\quad
		F_\mathcal{M}(x,-W_x) < 1,
		\qquad \forall x \in M_+.
	\end{equation}
\end{corollary}
\begin{proof}
		Fix $x\in M_+$.
	If $b_x<1/2$, then, by Proposition~\ref{generalcase},
	$F(x,\cdot)$ is a Minkowski norm if and only if $F_{\mathcal M}(x,-W_x)<1$.
	Conversely, assume that $F(x,\cdot)$ is a Minkowski norm,
	so that $\mathcal B_x$ is its unit ball and
	$\partial\mathcal B_x$ is a smooth, strongly convex hypersurface
	enclosing $0_x$.
	We first observe that
	\begin{equation}\label{boundaryinS}
		\partial\mathcal B_x\subseteq\mathcal S_x:=\mathcal Y_x(UT_xM).
	\end{equation}
	Indeed, let $v\in\partial\mathcal B_x$, so that $F(v)=1$.
	Since $\mathcal B_x=\operatorname{rad}\mathcal S_x$, there exist
	$z\in\mathcal S_x$ and $a\in(0,1]$ such that $v=az$;
	as $z\in\mathcal B_x$, we get $1\geq F(z)=F(v)/a=1/a$,
	hence $a=1$ and $v=z\in\mathcal S_x$.

	Let us show that $b_x\leq1/2$. Assume, by contradiction, that $b_x>1/2$,
	and consider the covector $p:=-h(Z_x,\cdot)$.
	For $u\in UT_xM$, setting $t:=h(Z_x,u)\in[-b_x,b_x]$,
	from~\eqref{eq:normalizedForm} we get
	\[
		p\bigl(\mathcal Y_x(u)\bigr)
		=\eta(x)\bigl(-t-t^2+p(W_x)\bigr)
		\leq\eta(x)\Bigl(\frac14+p(W_x)\Bigr)=:m,
	\]
	with equality if and only if $t=-1/2$.
	Since $b_x>1/2$ and $n\geq2$, there exist two distinct vectors
	$u_1,u_2\in UT_xM$ such that $h(Z_x,u_i)=-1/2$, and the points
	$y_i:=\mathcal Y_x(u_i)=\eta(x)\bigl(\tfrac12u_i+W_x\bigr)$,
	$i=1,2$, are distinct.
	The maximum of $p$ on the convex body $\mathcal B_x$ is attained at
	a point of $\partial\mathcal B_x$; hence,
	by~\eqref{boundaryinS} and $\mathcal S_x\subseteq\mathcal B_x$,
	it is equal to $\max_{\mathcal S_x}p=m$.
	Thus, $y_1$ and $y_2$ belong to $\mathcal B_x$ and to its supporting
	hyperplane $\{p=m\}$, and so does the segment joining them,
	which is then contained in $\partial\mathcal B_x$.
	This is impossible, because $\partial\mathcal B_x$ is strongly convex.

	Hence $b_x\leq1/2$, and~\eqref{radialgraph} gives
	$\mathcal S_x=\eta(x)\bigl(\partial\mathcal M_x+W_x\bigr)$.
	By~\eqref{boundaryinS},
	$\partial\mathcal B_x\subseteq
	\eta(x)\bigl(\partial\mathcal M_x+W_x\bigr)$.
	Both sets are homeomorphic to a sphere of dimension $n-1\geq1$;
	hence $\partial\mathcal B_x$, being open by invariance of domain
	and closed by compactness, coincides with $\eta(x)\bigl(\partial\mathcal M_x+W_x\bigr)$.
	Therefore $\partial\mathcal M_x$ is strongly convex,
	and $b_x<1/2$ by Proposition~\ref{prop:FinslerM-conditions}.
	So it remains to show the smoothness of
	$F\colon TM_+ \to [0,+\infty)$
	if~\eqref{eq:FinslerM+conditions} holds.

	For $y\in TM_+\setminus0$, let us set 
	$t(x,y):=\eta(x)F(x,y)$, so that $t(x,y)$
	is the unique positive solution of
	\[
		\Phi(x,y,t):=F_{\mathcal M}(x,y-tW_x)-t=0.
	\]
	(cf. Proposition~\ref{minimalgauge}).
	At such a solution, putting $z=y-tW_x$, convexity and positive	homogeneity of $F_{\mathcal M}$ give
	\[\partial_t\Phi(x,y,t)	=	d_vF_{\mathcal M}(x,z)[-W_x]-1	\leq	F_{\mathcal M}(x,-W_x)-1	<0.\]
	Thus, the implicit function theorem shows that $t$, and therefore	$F=t/\eta$, is smooth on $TM_+\setminus0$.
	It remains to prove continuity along the zero section. Fix  $x_0\in M_+$. Since
	$q(x):=F_{\mathcal M}(x,-W_x)<1$ and both $q$ and $\eta$ are continuous, there exist a neighbourhood
	$U$ of $x_0$ and constants $q_0<1$, $\eta_0>0$ such that
	\[
		q(x)\le q_0,\qquad\eta(x)\ge \eta_0,\qquad x\in U.
	\]
	For $(x,y)\in TM_+\setminus 0$, set $t=\eta(x)F(x,y)$.
	By the defining equation $F_{\mathcal M}(x,y-tW_x)=t$
	and the triangle inequality for the Minkowski norm $F_{\mathcal M}(x,\cdot)$,
	we have
	\[
		t=F_{\mathcal M}(x,y+t(-W_x))\le F_{\mathcal M}(x,y)+tF_{\mathcal M}(x,-W_x).
	\]
	Hence, for $x\in U$, $(1-q_0)t\le F_{\mathcal M}(x,y)$ and therefore
	\[
		0\le F(x,y)
		=\frac{t}{\eta(x)}\le\frac{F_{\mathcal M}(x,y)}{\eta_0(1-q_0)}.
	\]
	Since $F_{\mathcal M}$ is continuous on $TM_+$ and $F_{\mathcal M}(x_0,0)=0$, it follows that $F(x,y)\to 0$ as $(x,y)\to (x_0,0)$.
	Positive homogeneity gives	then the continuous extension $F(x,0)=0$.
\end{proof}

An explicit solution $F$ of \eqref{eqforF} is not always possible;
nevertheless, under some suitable conditions,
a lower and an upper bound for $F$ in terms of $\alpha$ can be given.
This estimate also yields a sufficient condition for backward and forward completeness of $F\colon TM_+ \to [0,+\infty)$,
that will be studied in Subsection~\ref{subsec:FM+completeness}.

\begin{proposition}\label{2bounds}
	Assume that~\eqref{eq:FinslerM+conditions} holds
    and let us set $q_x:=F_{\mathcal M}(-W_x)$.
	Then, for any $(x,y)\in TM_+$,
	\begin{equation}\label{F2bounds}
		\frac{\alpha_x(y)}{\eta(x)(1+b_x)(1+q_x)}
		\leq F(x,y)
		\leq\frac{\alpha_x(y)}{\eta(x)(1-b_x)(1-q_x)}.
	\end{equation}
	In particular,
	\begin{equation}\label{eq:F-LowUpp-bounds}
		\frac{\alpha_x(y)}{3\eta(x)}
		\leq	F(x,y)
		\leq	\frac{2\alpha_x(y)}{\eta(x)\bigl(1-F_{\mathcal M}(x,-W_x)\bigr)}.
	\end{equation}
\end{proposition}
\begin{proof}
	The inequalities are immediate for $y=0_x$.
	Let $y\neq0_x$ and set $t:=\eta(x)F(x,y)$.
	From \eqref{eqforF}, we know that
	$F_{\mathcal M}(x,y-tW_x)=t$.
	Writing $z:=y-tW_x$, we therefore have $F_{\mathcal M}(z)=t$.
	From \eqref{Matsumoto}, we deduce for any $(x,y)\in TM_+$
	\[
		\frac{\alpha_x(y)}{1+b_x}
		\leq F_{\mathcal M}(x,y)
		\leq	\frac{\alpha_x(y)}{1-b_x}.
	\]
	Hence
	\[
		\alpha_x(z)\leq(1+b_x)t,
		\qquad	\alpha_x(W_x)	\leq	(1+b_x)F_{\mathcal M}(x,-W_x)	=	(1+b_x)q_x.
	\]
	Since $y=z+tW_x$, it follows that
	\[
		\alpha_x(y)
		\leq	\alpha_x(z)+t\alpha_x(W_x)
		\leq	(1+b_x)(1+q_x)t,
	\]
	which gives the lower bound in~\eqref{F2bounds}.
	For the upper bound,
	the triangle inequality for $F_{\mathcal M}$ yields
	\[
		t=	F_{\mathcal M}(x,y-tW_x)
		\leq F_{\mathcal M}(x,y)		+tF_{\mathcal M}(x,-W_x)
		\leq\frac{\alpha_x(y)}{1-b_x}+tq_x,
	\]
	i.e. 
	\[	
		t\leq\frac{\alpha_x(y)}{(1-b_x)(1-q_x)}.
	\]
	Finally, as $b_x<1/2$ and $q_x<1$, we obtain
	$(1+b_x)(1+q_x)<3$ and $1/(1-b_x)<2$,
	so~\eqref{eq:F-LowUpp-bounds} follows.
\end{proof}
For future reference, we also record the following result. 
\begin{proposition}\label{Theta}
	Let $x \in M_+$ be such that $b_x <1/2$
	and $F_\mathcal{M}(x,-W_x) < 1$.
	Then the map $\mathcal Y_x\colon UT_xM\to\partial\mathcal B_x$ is a diffeomorphism, and so are
	\[
		\mathcal{U}_x\colon\partial\mathcal B_x\to UT_xM,\quad 	\mathcal{U}_x(v):=\frac{v}{\alpha(v)}, 	
		\qquad\text{and}\qquad 	\Theta_x:=\mathcal{U}_x\circ\mathcal Y_x.
	\]
\end{proposition}
\begin{proof}
	Since $b_x<1$, we have $\rho_x(u)=1+h(Z_x,u)\geq1-b_x>0$ for every
	$u\in UT_xM$, so that $\mathcal X_x(u)=\rho_x(u)u$ maps $UT_xM$
	onto 	$\partial\mathcal M_x$ by~\eqref{radialgraph}, with smooth
	inverse $v\in\partial\mathcal M_x \mapsto v/\alpha(v)$, hence it is
	a diffeomorphism. By Proposition~\ref{generalcase}, 	$\mathcal
	B_x=\eta(x)(\mathcal M_x+W_x)$ is a compact, smooth, strongly
	convex domain  	containing $0_x$ in its interior and
	$\partial\mathcal B_x=\eta(x)(\partial\mathcal M_x+W_x)=\mathcal
	Y_x(UT_xM)$;
	composing $\mathcal X_x$ with the affine diffeomorphism 	$v\mapsto\eta(x)(v+W_x)$ gives the first assertion.
	
Since $0_x\in\operatorname{int}\mathcal B_x$ and $\mathcal B_x$ is convex,
every ray issuing from $0_x$ meets $\partial\mathcal B_x$ in exactly one point,
so $\mathcal{U}_x$ is bijective.
If $n_v$ denotes the outward $h_x$-unit normal to $\partial\mathcal B_x$ at $v$,
convexity gives $h(n_v,z-v)\leq0$ for every 	$z\in\mathcal B_x$;
choosing $\varepsilon>0$ with $\varepsilon n_v\in\mathcal B_x$,
we obtain $h(n_v,v)\geq\varepsilon>0$, so that the radial direction $v$
is transverse to $T_v\partial\mathcal B_x$.
As
\[
	(d\mathcal{U}_x)_v(\xi)=\frac{\xi}{\alpha(v)}-\frac{h(v,\xi)}{\alpha(v)^3}\,v
\]
has kernel $\mathbb Rv$ on $T_xM$,
its restriction to $T_v\partial\mathcal B_x$
is an isomorphism onto $T_{\mathcal{U}_x(v)}UT_xM$.
Hence, $\mathcal{U}_x$, and therefore $\Theta_x$, is a diffeomorphism.
\end{proof}

\subsection{Vanishing self-propulsion}\label{subsec:zero-eta}
We now return to the fundamental law~\eqref{eq:generalForm}
and consider it at a fixed point $x\in M$ where $\eta(x)=0$,
so that
\begin{equation}\label{eq:zero-eta-law}
	\mathcal Y(u)=h(\widehat Z,u)u+\widehat W,
	\qquad u\in UT_xM.
\end{equation}
If $\widehat Z=0$, then $\mathcal Y$ is constant
and $\mathcal B =[0,\widehat W] $.
We therefore assume in the rest of
this subsection that $\widehat Z\neq0$.

Let us first omit the translation by $\widehat W$,
so let us consider the set 
\[
	\mathcal K:=\bigl\{\lambda h(\widehat Z,u)u:
	u \in UT_xM,\ 0\leq\lambda\leq1\bigr\}
	\subset T_xM.
\]
A direct computation gives
\begin{align}
	\mathcal K
&=\bigl\{v\in T_xM:\alpha^2(v)\leq h(\widehat Z,v)\bigr\}\label{K}\\
&=\left\{v\in T_xM:
	\alpha^2\left(v-\frac12\widehat Z\right)
\leq\frac14\alpha^2(\widehat Z)\right\}.\nonumber
\end{align}
Thus $\mathcal K$ is the closed $h$-ball centered at
$\widehat Z/2$ with radius $\alpha(\widehat Z)/2$.
In particular, it is strongly convex,
but the origin belongs to its boundary
and so its gauge function is not finite on the whole $T_xM$.
On the cone
\[
	A:=\bigl\{y\in T_xM:h(\widehat Z,y)>0\bigr\},
\]
one has, for $\lambda>0$, $y\in\lambda\mathcal K$
if and only if
$\alpha^2(y)\leq\lambda h(\widehat Z,y)$.
Hence
\begin{equation}\label{Kropina}
	F_{\mathcal K}(y)
	:=\inf\bigl\{\lambda>0:y\in\lambda\mathcal K\bigr\}
	=\frac{\alpha^2(y)}{h(\widehat Z,y)},
	\qquad y\in A.
\end{equation}
This is a {\em Kropina norm} (see \cite{YosSab12,CaponioJavaloyesSanchez2024}) with indicatrix
\[
\{y\in A:F_{\mathcal K}(y)=1\}=\partial\mathcal K\setminus\{0\}.
\]

Let us now consider the translation and set
$\mathcal C:=\mathcal K+\widehat W$.
By~\eqref{eq:zero-eta-law} and~\eqref{K},
\begin{align}\label{eq:zero-eta-navigation}
	&\mathcal Y(UT_xM)=
	\partial\mathcal K+\widehat W=\partial\mathcal C
	= C + r\,UT_xM,\\
	&\text{where}\quad r:=\tfrac12\alpha(\widehat Z)>0,
	\qquad
	 C:=\tfrac12\widehat Z+\widehat W,\label{RC}
\end{align}
so that $\mathcal C$ is the closed $h$-ball of center $C$
and radius $r$, and $\mathcal B$ is the radial hull
of the sphere $\partial\mathcal C$.
The parametrization $u\mapsto\mathcal Y(u)$ is not injective:
it identifies $u$ with $-u$ and collapses the equator
$\{u\in UT_xM:h(\widehat Z,u)=0\}$ to the single point $\widehat W$;
nevertheless, every point of that sphere is attained.
Let us also define 
\begin{equation}\label{eq:Kscript}
\mathscr K:=\alpha^2(C)-r^2 =\alpha^2(\widehat W)+h(\widehat Z,\widehat W),
\end{equation}
so that the sign of $\mathscr K$ records the position of $0$
relative to the translated Riemannian ball:
interior, boundary or exterior.
These are exactly the mild-, critical- and strong-wind alternatives that appear in Zermelo navigation as identified in~\cite{CaponioJavaloyesSanchez2024}
(see also Subsection~\ref{subsec:zermelo}).
More formally, by applying Lemma~\ref{spheregauge}
we have the following fiberwise classification at
any point where self-propulsion vanishes.

\begin{proposition}\label{zeroeta-classification}
	Let $x\in M$ satisfy $\eta(x)=0$ and $\widehat Z_x\neq0_x$,
	and let $C$, $r$, and $\mathscr K$ be as
	in~\eqref{RC} and~\eqref{eq:Kscript}.
	Then
	\begin{equation*}
		\mathcal C=\bigl\{y\in T_xM:
		\alpha^2(y)\leq h(\widehat Z+2\widehat W,y)-\mathscr K\bigr\},
		\qquad
		\mathcal B=\operatorname{co}\bigl(\{0\}\cup\mathcal C\bigr).
	\end{equation*}
	The gauge of $\mathcal B$ is finite
	on the cone $\Lambda=\operatorname{cone}(\mathcal C)$ and,
	for $y\in\Lambda\setminus\{0\}$, is given by
	\begin{equation}\label{eq:zero-eta-gauge}
		F(y)=\frac{2\alpha^2(y)}{h(\widehat Z+2\widehat W,y)
		+\sqrt{h(\widehat Z+2\widehat W,y)^2-4\mathscr K\,\alpha^2(y)}}.
	\end{equation}
	Moreover:
	\begin{enumerate}
		\item[(i)] if $\mathscr K<0$, then $0\in\operatorname{int}\mathcal B$, $\Lambda=T_xM$, and $F$ is a Randers norm;
		\item[(ii)] if $\mathscr K=0$, then
			\[
				\Lambda=\{y\in T_xM:h(\widehat Z+2\widehat W,y)>0\}\cup\{0\},
			\]
			and on its open half-space of finiteness
			\[
				F(y)=\frac{\alpha^2(y)}{h(\widehat Z+2\widehat W,y)},
			\]
			so $F$ is a Kropina norm;
		\item[(iii)] if $\mathscr K>0$, then $\Lambda$ is a proper closed cone, $0\in\partial\mathcal B$, 
			and $F$ is a conic Minkowski norm on $\operatorname{int}\Lambda$.
	\end{enumerate}
\end{proposition}
\begin{proof}
	By~\eqref{eq:zero-eta-navigation}, $\mathcal C$
	is the closed $h$-ball with center
	$C=\widehat W+\widehat Z/2$
	and radius $r=\alpha(\widehat Z)/2$,
	while $\mathcal B$ is the radial hull of its boundary.
	Lemma~\ref{spheregauge}, with $W=C$, and $\Delta=\mathscr K$,
	gives $\mathcal B=\operatorname{co}(\{0\}\cup\mathcal C)$
	and the formula for its gauge.
	Since $2C=\widehat Z+2\widehat W$,
	the expression in Lemma~\ref{spheregauge}
	is exactly~\eqref{eq:zero-eta-gauge}.
	The three alternatives follow respectively
	from $\Delta<0$, $\Delta=0$, and $\Delta>0$ in that lemma.
	In the critical case the denominator
	in~\eqref{eq:zero-eta-gauge} equals
	$2h(\widehat Z+2\widehat W,y)$
	on the finiteness half-space,
	giving the displayed Kropina formula.
\end{proof}

\begin{remark}\label{rem:gluing}
	If $\eta$ vanishes on an open region $M_0\subset M$ and $\widehat Z\neq 0$
	everywhere on $M_0$, the fiberwise Randers, Kropina, or strong-wind conic
	Minkowski norms given by Proposition~\ref{zeroeta-classification} assemble
	into the Lagrangian~\eqref{eq:zero-eta-gauge} on $M_0$, which is smooth on
	$\bigcup_{x\in M_0}\operatorname{int}\Lambda_x$. However, this structure
	cannot be smoothly glued to the Finsler metric $F$ on $M_+$, because both
	standing hypotheses of the positive self-propulsion regime, namely
	$b_x<1/2$ and~\eqref{translation}, degenerate at the zeros of $\eta$.

	Indeed, the condition $b_x<1/2$ reads
	$\alpha(\widehat Z_x)<\eta(x)/2$. Moreover, since
	$F_{\mathcal M}\geq\alpha/(1+b)$ by~\eqref{Matsumoto},
	condition~\eqref{translation} implies $\alpha(W_x)<1+b_x<3/2$, that is,
	\begin{equation*}
		\alpha(\widehat W_x)<\frac32\,\eta(x),\qquad x\in M_+ .
	\end{equation*}
	Hence, if both hypotheses hold on $M_+$, by continuity $\widehat Z$ and
	$\widehat W$ vanish on $\partial M_+$, and so does the whole velocity law.
	In particular, $\mathcal B_x=\{0_x\}$ at every $x\in\partial M_+$, and
	Proposition~\ref{zeroeta-classification} never applies at such points.
	This is consistent with
	Proposition~\ref{unreachable} below, which
	shows that under the same hypotheses the set $\{\eta=0\}$ is at infinite
	travel time from $M_+$.
\end{remark}

\section{Some classes of examples}\label{examples}
In this section we present several models that fit within
our framework: some well-known and widely studied ones are recovered
as particular cases, while for wildfire propagation and microswimmer
navigation the framework provides alternative kinematic models.

\subsection{Matsumoto slope-of-a-mountain metric}
\label{subsec:Matsumoto}
Let $n\ge2$, let $\Omega\subset\R^n$ be open, and let $\Sigma_f\subset\R^{n+1}$ be the graph of a smooth function $f:\Omega\to\R$. We use Euclidean  coordinates $x=(x^1,\ldots,x^n)$ and write $f_i=\partial f/\partial x^i$, with indices ranging from $1$ to $n$. The Euclidean metric on $\R^{n+1}$ induces on the graph $\Sigma_f$ the Riemannian metric $h$ whose components are
\[
h_{ij}=\delta_{ij}+f_i f_j.
\]
The inverse matrix of $(h_{ij})$ is given by the Sherman--Morrison 
formula \cite{ShermanMorrison1950}:
\[
h^{ij}=\delta^{ij}-\frac{f_i f_j}{1+|Df|^2},
\]
where $|Df|^2=\sum_{i=1}^n f_i^2$. Thus,
\begin{equation}\label{normgradhf}\alpha(\nabla^h f)=\frac{|Df|}{\sqrt{1+|Df|^2}},
\end{equation}
where $\alpha$ is the norm associated with $h$. 

The drift acting on the system is assumed to be proportional to the tangential component of gravity. Accordingly, in this kinematic model we write $V = -g \,\nabla^h f$, where $g>0$ denotes an effective gravity-induced velocity scale (thus $g$ has the physical dimension of speed). 
In~\cite{Matsumoto1989}, M. Matsumoto considered exactly this situation (with $n=2$), neglecting the component of $V$ orthogonal to the direction $u$.
Indeed, he wrote:
\begin{quote}
	\it Actually the earth's gravity $g \sin a$ acts on him;
	the component $(g \sin a) \sin \theta$ perpendicular to $OU$ is regarded to be cancelled by planting his legs on the road $OU$.
\end{quote}
where $g \sin a$ gives the modulus of the projection of the gravitational field on the tangent plane and $\theta$ is the angle between this projection and the  direction $OU$. In our setting, Matsumoto's assumption corresponds to retaining the
entire longitudinal contribution while cancelling the transverse one, that is,
\[
\mu\equiv1,\qquad\sigma\equiv0.
\]
Assume first that $\eta>0$ on $\Omega$. Then, according to~\eqref{eq:Y-oneV}, the velocity law is
\begin{equation}\label{Matsumoto_slope_velocity_law}
    \mathcal{Y}(u) = \big(\eta  - g \,h(\nabla^h f,u)\big)u=\eta\left(1-\frac{g}{\eta}\,h(\nabla^h f,u)\right)u
\end{equation}
and then, by Corollary~\ref{cor:FM-Finsler}, it defines a Finsler metric on $\Sigma_f$
if and only if 
\[
\frac{g}{\eta(x)}\alpha\big(\nabla^h f(x) \big)< \frac 1 2,	\qquad \forall x \in \Omega,
\]
which, by \eqref{normgradhf}, is equivalent to
\[
4g^2|Df(x)|^2<\eta^2(x)(1+|Df(x)|^2).\]
In the case $\eta\equiv 1$ and choosing units such that $g=1$, this becomes $|Df(x)|^2<1/3$  (see  \cite[p. 202]{BaoRobles2004}).
At a point where $\eta(x)=0$ and $\nabla^h f(x)\neq0_x$, one instead has $\widehat Z_x=-g\nabla^h f(x)$ and $\widehat W_x=0_x$. The fiberwise gauge is therefore the Kropina norm~\eqref{Kropina} on the downhill cone $\{y:h_x(\nabla^h f(x),y)<0\}$.

\subsection{Zermelo navigation: Randers, Kropina and strong-wind regimes}
\label{subsec:zermelo}
At a fixed point $x\in M$, Zermelo navigation data consist of a maximal speed $r_x>0$ and a wind vector $C_x\in T_xM$; see, for example, \cite{BaoRoblesShen2004,YosSab12,CaponioJavaloyesSanchez2024}. The
maximal velocity sphere is
\[
\mathscr S_x^{\mathrm{nav}}:=r_xUT_xM + C_x,
\]
that is, the translation by $C_x$ of the indicatrix of the conformal Riemannian metric $h_x/r_x^2$. The corresponding Zermelo velocity ball is
\[
\mathcal C_x:= r_x\overline B_x^h+C_x,\qquad\partial\mathcal C_x=\mathscr S_x^{\mathrm{nav}}.
\]
In accordance with \eqref{admissible}, we associate with the maximal velocity sphere its
radial hull  $\mathcal B_x:=\operatorname{rad}\bigl(\mathscr S_x^{\mathrm{nav}}\bigr)$.
When $0_x\in\mathcal C_x$, corresponding to the mild or critical wind  regimes in the terminology of \cite{CaponioJavaloyesSanchez2024}, one has $\mathcal B_x=\mathcal C_x$.
In the strong-wind regime, however, $0_x\notin\mathcal C_x$ and $\mathcal C_x\subsetneq\mathcal B_x$
so that $\mathcal B_x$ should be regarded as the radial relaxation of the standard Zermelo velocity set.
However, we will show in Theorem~\ref{maximalspeed}
that this distinction does not affect the least travel time branch.

Lemma~\ref{spheregauge}, applied with $W=C_x$ and $r=r_x$,  describes the  gauge of $\mathcal B_x$ in terms of the function $\Delta_x:=\alpha_x^2(C_x)-r_x^2$; the gauge is then given by \eqref{Fwind}
on the cone $\Lambda_x$.  We follow the terminology of~\cite{CaponioJavaloyesSanchez2024}, where $\Sigma_x:=\mathscr S^{\mathrm{nav}}_x$ is a \emph{wind Minkowskian structure} on $T_xM$ and the three cases below are the regions of mild, critical and strong wind:
\begin{itemize}
	\item \emph{mild wind}, $\Delta_x<0$: here $\Lambda_x=T_xM$ and $0_x\in\operatorname{int}\mathcal B_x$, so $F_x$ is a Minkowski norm on the whole of $T_xM$; rationalizing~\eqref{Fwind} gives the Randers metric
	\[
	F_x(y)=\frac{\sqrt{h_x(C_x,y)^2-\Delta_x\alpha_x^2(y)}-h_x(C_x,y)}{-\Delta_x};
	\]
	\item \emph{critical wind}, $\Delta_x=0$: here $\Lambda_x$ is the half-space $\{h_x(C_x,\cdot)>0\}$ together with $0_x$, and~\eqref{Fwind} becomes the Kropina norm $F_x(y)=\alpha_x^2(y)/\bigl(2h_x(C_x,y)\bigr)$; together with the mild case, this is the \emph{Randers--Kropina} class of~\cite[Def.~2.59]{CaponioJavaloyesSanchez2024};
	\item \emph{strong wind}, $\Delta_x>0$: here $\Lambda_x$ is the closed cone of half-aperture $\arcsin\bigl(r_x/\alpha_x(C_x)\bigr)$ about $C_x$, and~\eqref{Fwind} is a conic Minkowski norm on its interior.
\end{itemize}
\begin{remark}[Radial relaxation and least travel time]
In the strong-wind regime our radial admissible set contains velocities which are not 	available in the  interpretation	of strong-wind Zermelo navigation in \cite{CaponioJavaloyesSanchez2024}.
This enlargement does not, however, alter the least-time branch along 	an admissible direction. Indeed, if a ray $\mathbb R_{>0}y$ meets 	$\mathcal C_x$, then it meets the maximal velocity sphere 	$\mathscr S_x^{\mathrm{nav}}=\partial\mathcal C_x$ at a farthest point
$r_+(y)y$. Every point $ry$ of either $\mathcal C_x$ or its radial 	relaxation on that ray satisfies 
\[0\leq r\leq r_+(y).	\]
Hence the largest admissible speed in the direction $y$ is the same 	for both descriptions and is attained on 	$\mathscr S_x^{\mathrm{nav}}$. Equivalently, the least-time gauge is 	determined by the farthest intersection of the ray with the translated 	sphere.
	
In the strong-wind case of \cite{CaponioJavaloyesSanchez2024} the ray generally has two intersections with 	$\mathscr S_x^{\mathrm{nav}}$. Formula~\eqref{Fwind} corresponds to the farthest intersection and therefore to the least-travel-time 	branch. The second intersection gives the Lorentzian branch $F_\ell$ described in \cite[Prop. 2.57-(iii)]{CaponioJavaloyesSanchez2024}.
\end{remark}
There are two ways in which these navigation data arise naturally from our velocity law. First, suppose that $\eta(x)>0$ and $\widehat Z_x=0_x$. Then
\begin{equation}\label{randers}
\mathcal Y_x(u)=\eta(x)u+\widehat W_x,
\end{equation}
so that $r_x=\eta(x)$ and $C_x=\widehat W_x$. 
For a single drift field acting without traction,
$\mu=\sigma=1$ and~\eqref{eq:hat-fields-one} gives exactly this situation: $\widehat Z=0$ and $\widehat W=V$.
Hence~\eqref{randers} becomes
$\mathcal Y(u)=\eta u+V$, and~\eqref{Fwind} gives
\[
	F(y)=\frac{\alpha^2(y)}{h(y,V)+
	\sqrt{h(y,V)^2+\bigl(\eta^2-\alpha^2(V)\bigr)\alpha^2(y)}}.
\]
Second, suppose that $\eta(x)=0$ and $\widehat Z_x\neq0_x$.
Although the parametrization~\eqref{eq:zero-eta-law} is not itself of the usual Zermelo form, its image is exactly the sphere $\mathscr S_x^{\mathrm{nav}}$ with the data displayed in~\eqref{eq:zero-eta-navigation}, and $\Delta_x=\mathscr K$ by~\eqref{eq:Kscript}.
Hence the alternatives $\mathscr K<0$,
$\mathscr K=0$ and $\mathscr K>0$
found in Subsection~\ref{subsec:zero-eta}
are precisely the mild-, critical- and strong-wind regimes,
and~\eqref{eq:zero-eta-gauge} is~\eqref{Fwind} written in terms of $\widehat Z_x$ and $\widehat W_x$.
In particular, a vanishing self-propulsion turns half of the deformation field $\widehat Z_x$ into an effective maximal speed $r_x$ and the other half, together with $\widehat W_x$, into an effective wind $C_x$.

\subsection{Other one-field models}\label{subsec:onevector}
In the case of a single vector field $V$, equation~\eqref{eq:Y-oneV} becomes
\begin{equation}\label{vl_one-vector}
\mathcal Y(u)=\bigl(\eta+(\mu-\sigma)h(V,u)\bigr)u+\sigma V,\qquad\alpha(u)=1.
\end{equation}
As in \eqref{Matsumoto}, let us set 
\[
\beta_x:=-h(Z_x,\cdot),\qquad b_x:=\|\beta_x\|_{\alpha}=\alpha(Z_x).
\]
On the region $M_+$, it is convenient to rewrite the regularity conditions of Proposition~\ref{generalcase} in terms of the parameters appearing in \eqref{vl_one-vector}. Let  $\ell_x:=\alpha(V_x)$. Since 
\[
Z=\frac{\mu-\sigma}{\eta}V,
\]
the condition $b_x<1/2$ is equivalent to
\begin{equation}\label{FMconvex}
\eta(x)>2|\mu(x)-\sigma(x)|\,\ell_x.
\end{equation}
Moreover, under this condition, the translation condition
$F_{\mathcal M}(-W_x)<1$ is equivalent to
\begin{equation}\label{Fconvex}
\eta(x)>|\mu(x)|\,\ell_x.
\end{equation}
If $V_x=0_x$, both conditions are automatic. Hence,  assume $\ell_x>0$ and set $e_x:=V_x/\ell_x$.
The two intersections of
$\mathcal Y_x(UT_xM)$ with the line $\mathbb RV_x$ are
\[\mathcal Y_x(e_x)= \bigl(\eta(x)+\mu(x)\ell_x\bigr)e_x, 
\qquad \mathcal Y_x(-e_x)= \bigl(\mu(x)\ell_x-\eta(x)\bigr)e_x.
\]
The origin lies in the interior precisely
when these points lie on opposite sides of it.
\begin{remark}
Conditions~\eqref{FMconvex} and~\eqref{Fconvex} constrain $\mu$ and 	$\mu-\sigma$, but neither of them involves $\sigma$ alone. Writing 
\[
\frac{\sigma(x)\,\ell_x}{\eta(x)}	=	\frac{\mu(x)\,\ell_x}{\eta(x)}	-	\frac{\bigl(\mu(x)-\sigma(x)\bigr)\ell_x}{\eta(x)},
\]
they  give
\begin{equation}\label{sigmarange}
|\sigma(x)|\,\ell_x<\tfrac32\,\eta(x),
\end{equation}
and \eqref{sigmarange} is sharp: for every $\sigma$ satisfying it, the choice $\mu=\frac23\sigma$ fulfils both~\eqref{FMconvex} and~\eqref{Fconvex}, since 	then $|\mu|\ell_x=\frac23|\sigma|\ell_x<\eta(x)$ and
$|\mu-\sigma|\ell_x=\frac13|\sigma|\ell_x<\frac12\eta(x)$. Thus, the transverse 	response to the drift may exceed the maximal self-speed by up to one half, 	while $F$ is still a Finsler metric. This is not a marginal effect: since
$n\geq2$, we may pick $u\in UT_xM$ with $h(V_x,u)=0$, so that 	$\mathcal Y(u)=\eta(x)u+\sigma(x)V_x$ and
\[
\alpha\bigl(\mathcal Y(u)\bigr)=\sqrt{\eta^2(x)+\sigma^2(x)\ell_x^2},
\]
which, as $|\sigma|\ell_x\to \frac32\eta(x)$, approaches	$\frac{\sqrt{13}}{2}\,\eta(x)\approx1.80\,\eta(x)$, so that steering orthogonally to 	$V_x$ can be almost twice as fast as moving in the absence of any drift.
\end{remark}
We consider now the models introduced in~\cite{AldeaKopacz2025a,AldeaKopacz2026} showing that they fit into our  framework.
Our setting is more general in several respects: the unperturbed speed $\eta\colon M\to[0,+\infty)$ may vary with the base point and may vanish, an arbitrary finite family of smooth vector fields $V_i$ may act simultaneously on the system, and their longitudinal and transverse contributions are
governed by independent smooth coefficients $\mu_i$ and $\sigma_i$.
Consequently, the deformation and translation fields $\widehat Z$ and $\widehat W$ in~\eqref{eq:hat-fields} need not be collinear.
By contrast, the two Aldea--Kopacz models considered below have normalized unperturbed speed $\eta\equiv1$ and involve a single external vector field, so that the corresponding fields $\widehat Z$ and $\widehat W$ (equivalently, $Z$ and $W$) are necessarily collinear.
In~\cite{AldeaKopacz2025a}, this external field is a gravitational field  and the two traction coefficients are constant.
The superwind model of~\cite{AldeaKopacz2026} allows an arbitrary wind and a spatially varying traction coefficient, but retains the one-parameter relation $\mu\equiv1$ and $\sigma=1-\lambda$. We now make these specializations explicit.

We first consider the model of~\cite{AldeaKopacz2025a}. Assume that $\eta\equiv1$ and identify the vector field $V$ in our model with the gravitational field $G_T$.  To avoid a conflict with our notation, let $\gamma,\widetilde\gamma\in[0,1]$ be the constant cross- and along-traction coefficients denoted respectively by $\eta$ and $\widetilde\eta$ in~\cite{AldeaKopacz2025a}. Their resultant velocity \cite[Eqs.~(16)--(18)]{AldeaKopacz2025a} is
\begin{equation}\label{AKvelocity}
v_{\gamma\widetilde\gamma}(u)=u+(1-\gamma)\operatorname{Proj}_{u^\perp}V
	+(1-\widetilde\gamma)\operatorname{Proj}_uV.
\end{equation}
Since $\alpha(u)=1$, we have
\[
\operatorname{Proj}_uV=h(V,u)u,\qquad\operatorname{Proj}_{u^\perp}V=V-h(V,u)u,
\]
and hence
\[
v_{\gamma\widetilde\gamma}(u)=\bigl(1+(\gamma-\widetilde\gamma)h(V,u)\bigr)u+(1-\gamma)V.
\]
Thus~\eqref{AKvelocity} coincides with our velocity law~\eqref{vl_one-vector} under the identification $\mu=1-\widetilde\gamma$ and $\sigma=1-\gamma$. Consequently, the two constructions yield the same admissible velocity sets and hence the same gauge functions. Whenever the common velocity set contains the origin in its interior and is strongly convex, they therefore define the same Finsler metric. Under this identification, \eqref{FMconvex}-\eqref{Fconvex} become 
\[
\ell_x<\frac{1}{2|\gamma-\widetilde\gamma|},\qquad\ell_x<\frac{1}{1-\widetilde\gamma},
\]
with the usual convention that a denominator equal to zero gives the bound $+\infty$. Hence they are equivalent to 
\begin{equation*}
\ell_x<\min\left\{\frac{1}{2|\gamma-\widetilde\gamma|},\;\frac{1}{1-\widetilde\gamma}\right\}.
\end{equation*}
Resolving which of the two terms realizes this minimum gives exactly the four parameter regions and the corresponding strong-convexity alternatives in~\cite[Th.~1 and Lemma~8]{AldeaKopacz2025a}. Thus their case distinction has a simple geometric origin here: one bound is the strong-convexity threshold for the Matsumoto deformation, while the other is the condition that the subsequent translation leave the origin in the interior. On the diagonal $\gamma=\widetilde\gamma$ the deformation disappears, $Z=0$, and only the Zermelo translation condition $(1-\gamma)\ell_x<1$ remains; at $(\gamma,\widetilde\gamma)=(1,1)$ the drift is completely cancelled and no restriction on $\ell_x$ is needed.

At the level of velocity laws, the superwind model of~\cite{AldeaKopacz2026} is obtained from the same algebraic correspondence by replacing the gravitational force  $G_T$ with their arbitrary wind $\overline W$, setting $\widetilde\gamma=0$, and allowing the remaining traction coefficient to depend smoothly on the base point, i.e.  $\gamma=\lambda(x)$, $\lambda\in C^\infty(M,[0,1])$. Combining~\cite[Eqs.~(2)--(3)]{AldeaKopacz2026}, their resultant velocity can be written as 
\[
v_\lambda(u)=u+\lambda\operatorname{Proj}_u\overline W+(1-\lambda)\overline W=
\bigl(1+\lambda h(\overline W,u)\bigr)u+(1-\lambda)\overline W.
\]
Therefore, it coincides with our velocity law~\eqref{eq:Y-oneV} under the identification
\[
\eta\equiv1,\qquad V=\overline W,\qquad\mu\equiv1,\qquad\sigma=1-\lambda.
\]
Equivalently, the normalized fields defined in~\eqref{eq:ZW-formula} are $Z=\lambda\overline W$ and $W=(1-\lambda)\overline W$. Under this identification, conditions~\eqref{FMconvex} and~\eqref{Fconvex} become, respectively,
\[
2\lambda(x)\|\overline W_x\|_\alpha<1,\qquad\|\overline W_x\|_\alpha<1.
\]
The first condition ensures the strong convexity of the underlying Matsumoto metric, whereas the second is the translation condition. Together, they are equivalent to 
\[
\|\overline W_x\|_\alpha<\begin{cases}
	1,	& 0\leq\lambda(x)\leq\frac12,\\[2mm]	\dfrac{1}{2\lambda(x)},
	& \frac12<\lambda(x)\leq1,
\end{cases}
\]
which is precisely the condition obtained in~\cite[Lemma~3.3]{AldeaKopacz2026}.
\subsection{Wildfire models}\label{wildfireJPS}
The model introduced in~\cite{JavaloyesPendasSanchez2023} describes wildfire propagation by means of a time-dependent family of indicatrices incorporating both slope and wind. The Riemannian metric of the terrain is fixed, whereas the remaining parameters may depend on time and position. Therefore, the comparison with our autonomous
model must be understood after fixing the time $t$.

At a more abstract level, Bressan~\cite{Bressan2007} models the burned region at time $t$ as the reachable set, from an initial burned region $R_0$, of a differential inclusion
\[
\dot\gamma(s)\in\mathcal F(\gamma(s)),
\]
where $\mathcal F(x)$ is compact and convex and contains the zero velocity.  Bressan's control problem, however, concerns the progressive construction of barriers that the trajectories cannot cross, whereas here we study the geometry of the underlying propagation law. 

The model considered  in \cite{JavaloyesPendasSanchez2023} prescribes the propagation speed directly as a function of the actual direction of motion. When only the slope is considered (thus neglecting the action of the wind), the authors assume that \emph{``the fire moves faster upwards than downwards, since the (vertical) flame is closer to the ground in the upward direction''} \cite[p.~420]{JavaloyesPendasSanchez2023}. In this case,  directional speed is affine and therefore the velocity law considered is of the type \eqref{eq:generalForm}. When wind is also included, the velocity law  remains radial, but its directional dependence is no longer affine and hence it does not, in general, have the form~\eqref{eq:generalForm}. We now make this distinction precise.

As in Subsection~\ref{subsec:Matsumoto}, represent the terrain as the graph of a smooth function $f\colon\Omega\subset\mathbb R^2\to\mathbb R$, with graph immersion
$\widehat z\colon\Omega\to\mathbb R^3$, $\widehat z(x)=\bigl(x,f(x)\bigr)$.
We identify the graph with $\Omega$ through $\widehat z$ and hence
regard $M=\Omega$ as endowed with the pullback metric $h_x(v,w)=\langle d\widehat z_x(v),d\widehat z_x(w)\rangle_{\mathbb R^3}$. Denoting its Riemannian norm by $\alpha$, we have
$\alpha_x(v)=\sqrt{h_x(v,v)}=\bigl\|d\widehat z_x(v)\bigr\|_{\mathbb R^3}$, where $\|\cdot\|_{\mathbb R^3}$ is the ambient Euclidean norm. We also set $\beta_{f}:=df=h\bigl(\nabla^h f,\cdot\bigr)$.
 To avoid confusing the Riemannian metric $h$ with the scalar function
denoted by $h(t,x)$ in~\cite{JavaloyesPendasSanchez2023}, we denote the latter by $\kappa(t,x)$. 

Let $e_\theta=(\cos\theta,\sin\theta)\in T_xM$ and set $\bar u_\theta:=e_\theta/\alpha(e_\theta)\in UT_xM$.
The vector $d\widehat z_x(\bar u_\theta)$ is precisely the actual unit direction denoted by $u_\theta$ in
\cite[Eq.~(3.3)]{JavaloyesPendasSanchez2023}. Since $d\widehat z_x(v)=\bigl(v,df_x(v)\bigr)$,
equation~\cite[Eq.~(3.2)]{JavaloyesPendasSanchez2023} gives
\[
\cos\delta(x,\theta)=\frac{df_x(e_\theta)}{\alpha(e_\theta)}=df_x(\bar u_\theta)=h\bigl(\nabla^h f(x),\bar u_\theta\bigr).
\]
Consequently, equation~\cite[Eq.~(3.1)]{JavaloyesPendasSanchez2023} becomes
\[
s_0(t,x,\bar u_\theta)=a(t,x)+\kappa(t,x)+\kappa(t,x)h\bigl(\nabla^h f(x),\bar u_\theta\bigr),
\]
where $a$ and $\kappa$ are smooth positive functions.
The  fire velocity, defined  in \cite[p. 421]{JavaloyesPendasSanchez2023}, is
\[
v_{\mathrm{fire}}(t,x,\theta)=s_0(t,x,\bar u_\theta)\,d\widehat z_x(\bar u_\theta)=
d\widehat z_x\left(s_0(t,x,\bar u_\theta)\bar u_\theta\right).
\]
Therefore, under the above identification and writing $u=\bar u_\theta$, their velocity law is
\begin{align}
	\mathcal Y^0_{t,x}(u)&=	s_0(t,x,u)u\nonumber\\
	&=	\Bigl(	a(t,x)+\kappa(t,x)	+\kappa(t,x)h\bigl(\nabla^h f(x),u\bigr)	\Bigr)u,
	\qquad
	u\in UT_xM.\label{wildfireslope}
\end{align}
It therefore coincides, for every fixed $t$, with~\eqref{eq:generalForm} under the identification
\[
\eta_t(x)=a(t,x)+\kappa(t,x),
\qquad
\widehat Z_t(x)=\kappa(t,x)\nabla^h f(x),
\qquad
\widehat W_t(x)=0_x.
\]
Since $a+\kappa>0$, the normalized fields are
\[
Z_t(x)=\frac{\kappa(t,x)}{a(t,x)+\kappa(t,x)}\nabla^h f(x),
\qquad W_t(x)=0_x.
\]
For each $x\in M$, its gauge is
\begin{equation}\label{slopegauge}
F^0(t,y)	=	\frac{\alpha^2(y)}{a(t,x)\alpha(y)+\kappa(t,x)\bigl(\alpha(y)+\beta_{f}(y)\bigr)}, \qquad y\in T_x M\setminus\{0_x\}.
\end{equation}
Equivalently, defining
\begin{equation}\label{betaMt}
(\beta_{\mathcal M,t})_x:=-h_x\bigl(Z_t(x),\cdot\bigr)=-\frac{\kappa(t,x)}{a(t,x)+\kappa(t,x)}\,(\beta_{f})_x,
\end{equation}
we obtain
\[
F^0(t,y)=\frac{1}{a(t,x)+\kappa(t,x)}\frac{\alpha^2(y)}{\alpha(y)-\beta_{\mathcal M,t}(y)}.
\]
Thus~\eqref{slopegauge} is precisely a conformal Matsumoto metric with the sign
convention adopted in~\eqref{Matsumoto}, with factor $1/(a+\kappa)$; equivalently, the isotropic speed $a+\kappa$ scales its indicatrix. The comparison with
Subsection~\ref{subsec:Matsumoto} is now transparent: there $Z=-\frac{g}{\eta}\nabla^h f$ points downhill and
$\beta=\frac{g}{\eta}\beta_f$ is a positive multiple of $\beta_f$, whereas here
$Z_t$ points uphill and, by~\eqref{betaMt}, $\beta_{\mathcal M,t}$ is a negative
multiple of $\beta_f$. Accordingly, in the terminology
of~\cite{JavaloyesPendasSanchez2023}, \eqref{slopegauge} is a \emph{reverse}
Matsumoto metric: fire is assumed to propagate faster uphill, whereas
Matsumoto's original model describes a walker moving faster downhill.

This correspondence also recovers exactly the strong convexity condition in~\cite{JavaloyesPendasSanchez2023}. Indeed, 
\[
b_t(x)=\alpha\bigl(Z_t(x)\bigr)=\frac{\kappa(t,x)}{a(t,x)+\kappa(t,x)}\alpha\bigl(\nabla^h f(x)\bigr).
\]
If $\varsigma(x)$ denotes the inclination angle of the terrain, so that $\cos\varsigma(x)=1/\sqrt{1+|Df(x)|^2}$,
then~\eqref{normgradhf} gives $\alpha\bigl(\nabla^h f(x)\bigr)=\sin\varsigma(x)$. Therefore, $b_t(x)<\frac12$ if and only if 
\[2\sin\varsigma(x)<\frac{a(t,x)+\kappa(t,x)}{\kappa(t,x)},\]
which is precisely~\cite[Eq.~(3.9)]{JavaloyesPendasSanchez2023}.

\begin{remark}[Direct direction--velocity laws]
The velocity law \eqref{wildfireslope} can be described also as follows.  Let 	$\mathscr S_t\colon UTM\to(0,+\infty)$ be smooth and define 
\begin{equation}\label{direct}
\mathcal Y^{\mathrm{dir}}_t(u):=\mathscr S_t(u)u.
\end{equation}
Thus the prescribed direction and the actual direction of motion coincide. 	For fixed $(t,x)$, let
\[
\mathcal B^{\mathrm{dir}}_{t,x}	:=	\operatorname{rad}\bigl(\mathcal Y^{\mathrm{dir}}_{t,x}(UT_xM)\bigr).
\]
By Lemma~\ref{radialgauge}-(i),
\[
\frac{y}{F^{\mathrm{dir}}(t,y)} \in \partial\mathcal B^{\mathrm{dir}}_{t,x} =
\mathcal Y^{\mathrm{dir}}_{t,x}(UT_xM).
\]
Since $\mathscr S_t>0$,  this hypersurface meets the ray through $y$ at the unique point
\[
\mathscr S_t\left(\frac{y}{\alpha(y)}\right) \frac{y}{\alpha(y)},
\]
so that we obtain
\begin{equation}\label{Fdirect}
F^{\mathrm{dir}}(t,y)	=	\frac{\alpha(y)}{\mathscr S_t\bigl(y/\alpha(y)\bigr)}.
\end{equation}

This construction is universal in the sense that  if 	$\widetilde F_t$ is any smooth time-dependent Finsler metric on $M$, then 	the choice
\[
\mathscr S_t(u):=\frac{1}{\widetilde F_t(u)}, 	\qquad u\in UTM,
\]
gives $F^{\mathrm{dir}}(t,\cdot)=\widetilde F_t$, by positive homogeneity.  Conversely, an arbitrary positive smooth function 	$\mathscr S_t$ yields the positively homogeneous function \eqref{Fdirect}, which is a Finsler metric precisely when $\mathcal Y^{\mathrm{dir}}_{t,x}(UT_xM)$ is strongly convex. 
	
At a fixed time, a direct law~\eqref{direct} is also of the  	form~\eqref{eq:generalForm} if and only if
$\widehat W_x=0_x$ and 
\[
\mathscr S_t(u)=\eta_t(x)+h_x(\widehat Z_t(x),u),
\]
with $\alpha(\widehat Z_t)<\eta_t$.

We now apply~\eqref{direct} to the wind-dependent model of~\cite{JavaloyesPendasSanchez2023}.
Motivated by Anderson's empirical double semielliptical description of wildfire propagation with wind \cite{Anderson1983}, Javaloyes, Pendás-Recondo and Sánchez replace the isotropic contribution $a(t,x)$ of the slope model by the radial function of an ellipse of eccentricity $\varepsilon(t,x)\in[0,1)$, whose major axis is aligned with the wind direction and whose ignition point is located at one of its foci. 
Let $\omega_{t}$ be the one-form introduced at \cite[p. 427]{JavaloyesPendasSanchez2023}. For
$u\in UT_xM$, the value $\omega_{t}(u)$ is the cosine of the actual angle between $u$ and the wind direction; in particular,
$\bigl|\omega_{t}(u)\bigr|\leq1$. In this notation, the  speed prescribed in~\cite[Eq.~(3.10)]{JavaloyesPendasSanchez2023} can be written as
\[
\mathscr S_{t}(u)=\frac{a(t,x)\bigl(1-\varepsilon(t,x)^2\bigr)}{1-\varepsilon(t,x)\omega_{t,x}(u)}+
\kappa(t,x)\bigl(1+\beta_{f}(u)\bigr),\qquad u\in UT_xM.
\]
Since $\varepsilon(t,x)<1$, the denominator of the first term is positive. Substitution into \eqref{Fdirect} gives
\begin{equation}\label{windgauge}
F^{\mathrm{dir}}(t,y)=\frac{\alpha^2(y)}{\displaystyle\frac{a(t,x)\bigl(1-\varepsilon(t,x)^2\bigr)\alpha^2(y)}{\alpha(y)-\varepsilon(t,x)\omega_{t}(y)}+ \kappa(t,x)\bigl(\alpha(y)+\beta_{f}(y)\bigr)}.
\end{equation}
This is precisely the function given in \cite[Eq.~(3.12)]{JavaloyesPendasSanchez2023}. When $\varepsilon=0$, formula~\eqref{windgauge} reduces to the slope metric~\eqref{slopegauge}. In general, however, the  dependence on $\omega_{t}(u)$ shows that the wind-dependent velocity law is not of the affine form~\eqref{eq:generalForm}.
Thus,  wind nonlinearly modifies the radial speed profile and hence the shape of the local propagation set. 

A different, Zermelo-type, model  may represent wind by an additive drift $A(t,x)$, thus  considering the velocity law $\mathcal Y^0_{t,x}(u)+A(t,x)$, with $\mathcal Y^0_{t,x}$  as in \eqref{wildfireslope}.  Whenever $A(t,x)$ and $\nabla^hf(x)$ are linearly independent, this provides a genuine instance of the noncollinear two-fields regime, with $\widehat Z_t=\kappa(t,x)\nabla^hf(x)$ and $\widehat W_t=A(t,x)$. We  will discuss further this model  in Subsection~\ref{subsec:additional-drift}.
\end{remark}

\subsection{A free-orientation relaxation of spheroidal microswimmer dynamics}
\label{subsec:microswimmers}
We finally consider a velocity law arising from models of spheroidal microswimmers. The orientational part of these models has a classical lineage. Jeffery derived the rotation law for a passive ellipsoid in a viscous flow~\cite{Jeffery1922}; Kessler introduced the gyrotactic reorientation mechanism for motile algae~\cite{Kessler1985}; and Pedley and Kessler combined these ingredients in a model for
spheroidal microorganisms swimming in a general flow field \cite{PedleyKessler1987}. More recent models incorporate active steering and orientation-dependent settling \cite{QiuEtAl2022}, or specialize the same position--orientation description to surface waves with wind-driven shear~\cite{VentrellaEtAl2023}.

Let $M\subset\mathbb R^3$ be a fluid domain endowed with the Euclidean metric $h$, let $\alpha$ be its associated norm, and let $e_3$ be the upward vertical unit vector. In the reduced particle description, an
axisymmetric spheroidal microswimmer is represented by the position $x\in M$ of its centre and by a unit vector $p\in\mathbb S^2$ directed along the symmetry axis of its body and oriented in the swimming
direction. Thus, its state is $(x,p)\in M\times\mathbb S^2$; the shape of the body is not spatially resolved, but enters the effective translation and orientation laws.

Let $U_t: M \to \mathbb{R}^3$ denote the prescribed background-fluid velocity; in the
setting of~\cite{VentrellaEtAl2023}, it contains the explicitly
time-dependent surface-wave field and,
in the two-dimensional model,
the wind-driven shear contribution.
Set
\[
	S_t:=\frac12\bigl(\nabla U_t+ (\nabla U_t)^{\mathsf T}\bigr)
    \in \mathbb{R}^{3 \times 3},
	\qquad\Omega_t:=\frac12\bigl(\nabla U_t-(\nabla U_t)^{\mathsf T}\bigr)
    \in \mathbb{R}^{3 \times 3}.
\]
The deterministic Jeffery--gyrotaxis equation  augmented, when
appropriate, by an actively imposed angular velocity $\omega_s$,
can be written as
\begin{equation}\label{microswimmer-orientation}
	\dot p	=\Omega_t(x)p	+\lambda_{\mathrm{sh}} 	\bigl(S_t(x)p-h\bigl(p,S_t(x)p\bigr)p\bigr) 	+\frac{1}{2B}\bigl(e_3-h(e_3,p)p\bigr)	+\omega_s\times p .
\end{equation}
Here $B>0$ is the gyrotactic reorientation time and
$\lambda_{\mathrm{sh}}$ is a shape parameter.
In \cite{VentrellaEtAl2023}, $\omega_s=0$ whereas in the
controlled model \cite{QiuEtAl2022} the angular velocity $\omega_s$ is
an action available to the swimmer. In either case, $p$ remains a
state variable subject to a differential constraint and is not an
instantaneously selectable direction.

Let us introduce the dimensionless kinematic parameters in the velocity law of a spheroidal microswimmer: $\nu>0$ is the swimming speed, $v_s>0$ is the settling speed in the orientation of highest drag, and $v_{sr}>1$
is the ratio between the settling speeds in the orientations of lowest and highest drag. The orientation-dependent settling law for a spheroid, standard in microhydrodynamics~\cite{KimKarrila1991} and used
in~\cite[Eq.~(2.1e)]{QiuEtAl2022}, reads in this notation
\begin{equation}\label{microswimmer-position}
	\dot x 	=U_t(x)+\nu p 	-v_s\Bigl(e_3+(v_{sr}-1)h(e_3,p)p\Bigr).
\end{equation}
This is also recalled in \cite[Eq.~(B2)]{VentrellaEtAl2023}, even though in the main body of that paper the authors instead use an orientation-independent vertical settling velocity and restrict the dynamics to a vertical plane.

The coupled system \eqref{microswimmer-orientation}--\eqref{microswimmer-position} is not a particular case of our control system on $M$, because the direction $p$ is dynamically constrained and the actual state space is $M\times\mathbb S^2$.
Nevertheless, it fits our framework if we discard \eqref{microswimmer-orientation} and replace $p$ in \eqref{microswimmer-position} by an arbitrary control direction $u\in UT_xM=\mathbb S^2$.
We call the resulting control system on $M$ the \emph{free-orientation system},
and the corresponding minimum-time problem the \emph{free-orientation problem}. The assumption of a freely steerable direction at fixed swimming speed is used as an idealized model in the literature on point-to-point minimum-time microswimmer navigation; see, for example,~\cite{DaddiMoussaIderEtAl2021}, where  a different hydrodynamic velocity law is used. Thus, let us  rewrite \eqref{microswimmer-position} as
\begin{equation}
\mathcal Y_{t,x}(u)=\bigl(\nu-\delta_s h(e_3,u)\bigr)u 	+U_t(x)-v_s e_3,
	\label{microswimmer-law}
\end{equation}
where $\delta_s:=v_s(v_{sr}-1)>0$.
Consequently, for every fixed $t$, equation~\eqref{microswimmer-law} has  the form~\eqref{eq:generalForm}, with 
\begin{equation*}
\eta_t=\nu, 	\qquad 	\widehat Z_t=-\delta_s e_3,  	\qquad	\widehat W_t=U_t-v_s e_3.
\end{equation*}
In particular, $\widehat Z_t(x)$ and $\widehat W_t(x)$ are linearly independent precisely at the points where $U_t(x)$ has a nonzero component orthogonal to $e_3$. Since the wave velocity in~\cite{VentrellaEtAl2023} depends explicitly on time, this identification belongs to the time-dependent extension of our framework and is to be understood fiberwise at each fixed $t$. 
If $U_t\equiv U$ is autonomous, it is  covered by the autonomous results of the preceding sections.

The normalized fields on $M$ are 
\[Z=-\frac{\delta_s}{\nu}e_3, \qquad W_t=\frac{U_t-v_s e_3}{\nu}\]
and hence $b=\alpha(Z)=\delta_s/\nu$. 
If $\delta_s<\nu$, the  untranslated velocity set is then the boundary of $\nu\mathcal M_x$, where
$\mathcal M_x$ is the unit ball of the Matsumoto gauge
\begin{equation*}
F_{\mathcal M}(y) 	=\frac{\alpha^2(y)}	{\displaystyle\alpha(y)-\frac{\delta_s}{\nu}h(e_3,y)},
\end{equation*}
which  is well defined on the whole tangent space and, by Proposition~\ref{prop:FinslerM-conditions}, it is   strongly convex if and only if $2\delta_s<\nu$. Under~this condition, Proposition~\ref{generalcase} shows that, for each fixed $t$, the radial hull $\mathcal B_{t,x}=\operatorname{rad}\big(\mathcal Y_{t,x}(UT_xM)\big)$
is the unit ball of a Finsler norm $F_t$ wherever $F_{\mathcal M}\bigl(-W_t(x)\bigr)<1$.
\begin{remark}[Velocities and radial relaxation]
Equation~\eqref{microswimmer-law} identifies the physically prescribed 	velocity set $\mathcal S_{t,x}:=\mathcal Y_{t,x}(UT_xM)$.	The passage from $\mathcal S_{t,x}$ to $\mathcal B_{t,x}$ is the radial
relaxation adopted throughout this paper.  In the autonomous case and for distinct endpoints,
Theorem~\ref{maximalspeed} shows that this radial relaxation does not 	change the infimum travel time.
\end{remark}
\begin{remark}[Orientation dynamics and time-optimal navigation]
Every solution $(x,p)$ of the system 	\eqref{microswimmer-position}--\eqref{microswimmer-orientation} 	projects to a trajectory of the free-orientation system by taking $u=p$. Hence the infimum travel time of the free-orientation problem is a lower bound for the infimum over any  class of 	trajectories satisfying the same positional endpoints and any 	additional prescribed orientation boundary conditions.
	
If 	$\delta_s<\nu$, then arguing as in the proof of Proposition~\ref{Theta}, we immediately see that 
 $\mathcal Y_{t,x}:\mathbb S^2\to T_xM$  is injective and, for 	$v\in\mathcal S_{t,x}$,
\[
\mathcal Y_{t,x}^{-1}(v)=\frac{v-\widehat W_t(x)}{\alpha\bigl(v-\widehat W_t(x)\bigr)}.
\]
Accordingly, the orientation along any maximal-velocity trajectory $\gamma$ of the free-orientation system is uniquely recovered as 
\[
u_\gamma(t)	=\frac{\dot\gamma(t)-U_t\bigl(\gamma(t)\bigr)+v_s e_3}
	{\alpha\!\left(		\dot\gamma(t)-U_t\bigl(\gamma(t)\bigr)+v_s e_3\right)}.
\]
If $u_\gamma$ is absolutely continuous and satisfies 	\eqref{microswimmer-orientation} almost everywhere for an admissible 	angular control $\omega_s$ (with $\omega_s=0$ in the model of 	\cite{VentrellaEtAl2023}) as well as the prescribed orientation boundary conditions, then $(\gamma,u_\gamma)$ is a solution  of the 
system \eqref{microswimmer-position}--\eqref{microswimmer-orientation}. If, moreover, $\gamma$ minimizes time for the free-orientation problem, the lower bound is then attained and the two minimum times agree. When $U_t\equiv U$ is autonomous and the Finsler regularity conditions 	above hold, every smooth nonconstant locally time-minimizing trajectory of the free-orientation problem, parametrized by travel time, is a unit-speed 	geodesic of $F$. Such geodesics are therefore candidates for 	time-optimal trajectories of the complete model. Ventrella et al. do not consider  a 	time-optimality issue as their purpose is to study trapping in 	surface waves. The minimum-time problem considered here is an 	additional navigation question, motivated by the optimal microswimmer 	literature. If $U_t$  depends on time, the family $F_t$ 	gives the instantaneous fiberwise gauges, but minimizers are not 	geodesics of a single frozen metric; they are extremals of the 	time-dependent control problem governed by the Hamiltonian of 	Subsection~\ref{subsec:additional-drift}.
\end{remark}

\section{The time-optimal problem}\label{sec:time-travel}

The velocity law \eqref{eq:generalForm} gives rise to a Finsler metric only
in specific regimes (see Propositions~\ref{generalcase}
and~\ref{zeroeta-classification}); in general, the admissible velocity sets
may be nonconvex and may fail to contain the origin in their interior.
For this reason, in this section we study the travel time problem \eqref{TT}
using only tools from the theory of differential inclusions, without relying
on any Finsler structure. Since the arguments depend only on the set-valued
maps $x\mapsto\mathcal S_x$ and $x\mapsto\mathcal B_x$, we do not even
assume that $\mathcal Y$ is of the form \eqref{eq:generalForm}. We compare
the infimum travel times obtained by allowing only the maximal
velocities, their radial hull, or its convex hull, and we then address the
existence of time-optimal trajectories. Finally, in
Subsection~\ref{subsec:FM+completeness}, we return to the regular Finsler
regime on $M_+$.
\subsection{Time-optimal trajectories for locally Lipschitz velocity laws}\label{timeoptimatraj}
Throughout this subsection, we assume that the map
$(x,u)\in UTM\mapsto \mathcal Y_x(u)\in T_xM$
is locally Lipschitz continuous.
For each $x\in M$ let
\[
	\mathcal S_x:=\mathcal Y_x(UT_xM)
\]
be the set of velocities and let us set
$\mathcal B_x=\operatorname{rad}\mathcal S_x$.
We also consider the convex hull
$\widehat{\mathcal B}_x:=\operatorname{co}\mathcal B_x$.
Since $0_x\in\mathcal B_x$ and $\mathcal S_x\subseteq\mathcal B_x$, we have
$\widehat{\mathcal B}_x=\operatorname{co}\bigl(\{0_x\}\cup\mathcal S_x\bigr)$.
Let $\widehat F$ be the gauge function of $\widehat{\mathcal B}_x$. Since
$\mathcal B_x\subseteq\widehat{\mathcal B}_x$, one has $\widehat F\leq F$; in
fact, $\widehat F$ is the largest convex positively homogeneous function
below $F$.
Indeed, $\widehat F$ is convex, being the gauge of a convex set containing
the origin; moreover, if $G\colon T_xM\to[0,+\infty]$ is convex and
positively homogeneous with $G\leq F$, then the convex set $\{G\leq1\}$
contains $\mathcal B_x=\{F\leq1\}$, hence also $\widehat{\mathcal B}_x$, and
Lemma~\ref{radialgauge}-(i), applied to $\widehat{\mathcal B}_x$, gives
$G\leq\widehat F$ by homogeneity.

We first compare~\eqref{TT} with the problem in which only 
velocities in $\mathcal S$ are allowed,
namely with
\[
	\mathcal T_{\mathcal S}(x_1,x_2)
	:= \inf \bigl\{T>0:\exists\gamma\in AC([0,T],M),\
		\gamma(0)=x_1,\ \gamma(T)=x_2,\
		\dot\gamma(t)\in\mathcal S_{\gamma(t)}
	\ \text{a.e.}\bigr\},
\]
Since $\mathcal{B}_x$ is given by arbitrary fractions
of the maximal velocities in $\mathcal{S}_x$,
this enlarges the class of admissible parametrizations,
so that
$\mathcal{T}(x_1,x_2)\le \mathcal{T}_{\mathcal{S}}(x_1,x_2)$
for any $x_1,x_2 \in M$.
We show that the reverse inequality holds for distinct endpoints,
thus obtaining Theorem~\ref{maximalspeed}.

When $\mathcal S=\bigcup_{x\in M}\mathcal S_x$ is a smooth hypersurface
of $TM$ transverse to the vertical distribution and every fiber
$\mathcal S_x$ is strongly convex, $\mathcal S$ is a wind Finslerian
structure in the sense of  \cite[Def.~2.8]{CaponioJavaloyesSanchez2024}; the corresponding
reparametrization statement for piecewise smooth curves is
\cite[Prop.~2.32]{CaponioJavaloyesSanchez2024}; the proof given there proceeds
through the pair of conic pseudo-Finsler metrics $F$ and $F_l$ associated with
$\mathcal S$, and requires in addition suitable continuity of $F\circ\dot\gamma$
and $F_l\circ\dot\gamma$ along the curve.
Here we impose no convexity on $\mathcal S_x$, which
need not be embedded and which a ray issuing from $0_x$ may meet more than
twice; accordingly the subdivision of $\mathcal S_x$ into a
convex and a concave part may fail and the conic pseudo-Finsler metrics  may fail to be  defined. Moreover, we consider  curves which are merely absolutely continuous.
The convexity hypothesis is replaced by compactness and radiality, indeed  Lemma~\ref{radialgauge} yields $y/F(y)\in\mathcal S_x$.

\begin{theorem}[The velocity set  and its radial hull]	\label{maximalspeed}
Let $x_1,x_2\in M$. If $x_1\neq x_2$, then
\begin{equation}\label{Tradial}
\mathcal T_{\mathcal S}(x_1,x_2)=\mathcal T(x_1,x_2),
\end{equation}
with values in $[0,+\infty]$. In particular, $x_2$ is reachable from $x_1$ for the  differential inclusion with the set $\mathcal S$ if and only if it is reachable for the radial inclusion. If $x_1=x_2=:x_0$, the same equality holds provided that
$0_{x_0}\in\mathcal S_{x_0}$.
\end{theorem}
\begin{proof}
Since $\mathcal S_x\subseteq\mathcal B_x$, one always has 	$\mathcal T(x_1,x_2)\leq\mathcal T_{\mathcal S}(x_1,x_2)$. For the 	converse inequality, let us assume first that $x_1\neq x_2$. Let 	$\gamma\in AC([0,T],M)$ satisfy the original differential inclusion 	$\dot\gamma(t)\in\mathcal B_{\gamma(t)}$ a.e., with
$\gamma(0)=x_1$ and $\gamma(T)=x_2$. We first observe that, since 	$\mathcal B_x=\operatorname{rad}\mathcal S_x$, for every $v\in\mathcal B_x\setminus\{0_x\}$ one has $v/F(v)\in\mathcal S_x$. Indeed, by Lemma~\ref{radialgauge}, $1/F(v)$ is the largest $s>0$ 	such that $sv\in\mathcal B_x$, so that  	$v/F(v)\in\mathcal B_x=\operatorname{rad}\mathcal S_x$, and there exist 	$a\in(0,1]$ and $z\in\mathcal S_x$ such that $v/F(v)=az$.  	If $a<1$, then $z=v/(aF(v))\in\mathcal S_x\subseteq\mathcal B_x$,  contradicting the maximality of $1/F(v)$. Hence $a=1$, and the claim follows.
	
Let us now show that the  gauge $F\colon TM\to[0,+\infty]$ is Borel measurable. 	Indeed, since $\mathcal Y$ is continuous, and the sphere bundle over a compact subset of $M$ is compact, the subset $\{(x,v)\in TM:v\in\mathcal B_x\}$ is closed and, by Lemma~\ref{radialgauge}-(ii), for every $c>0$,
\[
\bigl\{(x,v)\in TM:F(v)\leq c\bigr\}	=	\bigl\{(x,v)\in TM:v\in c\mathcal B_x\bigr\}
\]
is closed as well. Moreover, $\{F=0\}$ is the zero section by 	Lemma~\ref{radialgauge}, and for $c<0$ the sublevel $\{F\leq c\}$ is empty. Thus, all the sublevel sets of $F$ are closed, and therefore $F$ is lower semicontinuous. In particular, $F$ is Borel measurable.	

We can therefore set $\lambda(t):=F\bigl(\dot\gamma(t)\bigr)$. Then $\lambda\colon[0,T]\to[0,1]$ is measurable and 
$\lambda(t)=0$ if and only if $\dot\gamma(t)=0_{\gamma(t)}$.
Moreover, whenever $\lambda(t)>0$, the initial observation in the proof   gives
\begin{equation}\label{normalized-velocity}
\frac{\dot\gamma(t)}{\lambda(t)}\in\mathcal S_{\gamma(t)}.
\end{equation}
We now remove the factor $\lambda$, which represents the reduction of 	the assigned velocity.
 Set
\[
s(t):=\int_0^t\lambda(\tau)\,d\tau,
\]
and  $L:=s(T)$. Notice first that $L>0$, since otherwise $\lambda=0$ a.e. on $[0,T]$, hence $\dot\gamma=0$ a.e., contradicting $x_1\neq x_2$. Since $s$ is continuous and nondecreasing, it maps $[0,T]$ onto $[0,L]$.
For $r\in[0,L]$, choose any $t\in[0,T]$ such that $s(t)=r$ and set $\sigma(r):=\gamma(t)$.
This definition does not require $s$ to be injective. Indeed, if $s(t_1)=s(t_2)$ with $t_1<t_2$, then
\[
0=s(t_2)-s(t_1) =\int_{t_1}^{t_2}\lambda(t)dt.
\]
Since $\lambda\geq0$, it follows that $\lambda=0$ a.e. on $[t_1,t_2]$, and therefore $\dot\gamma=0$ a.e. on this interval. Thus $\gamma$ is constant on $[t_1,t_2]$, and in particular $\gamma(t_1)=\gamma(t_2)$. Hence $\sigma$ is well-defined and  $\gamma=\sigma\circ s$.

Since the image of $\gamma$ is compact and $\mathcal Y$ is continuous, 	there exists $C>0$ such that 	$\alpha\bigl(\mathcal Y_x(u)\bigr)\leq C$ for every 	$x\in\gamma([0,T])$ and $u\in UT_xM$. By 	\eqref{normalized-velocity}, it follows that 
\[
\alpha(\dot\gamma(t))\leq C\lambda(t) 	\qquad\text{a.e.}
\]
Therefore, if $r_1\leq r_2$, there exist $t_1\leq t_2$ such that 	$s(t_1)=r_1$, $s(t_2)=r_2$, and
\[
d_h\bigl(\sigma(r_1),\sigma(r_2)\bigr)
=d_h\bigl(\gamma(t_1),\gamma(t_2)\bigr)\leq\int_{t_1}^{t_2}\alpha(\dot\gamma(t))dt
\leq C\int_{t_1}^{t_2}\lambda(t)dt	\leq C(r_2-r_1).
\]
Thus $\sigma$ is Lipschitz and hence differentiable a.e. on $[0,L]$. We now show that its velocity belongs to
$\mathcal S_{\sigma(r)}$ for a.e. $r$. 
	
Let $E\subset[0,T]$ be a measurable set of measure zero such that, for 	every $t\in[0,T]\setminus E$, the functions $s$ and $\gamma$ are differentiable at $t$, $s'(t)=\lambda(t)$, and 	\eqref{normalized-velocity} holds whenever $\lambda(t)>0$. 	We claim that, for every measurable set $A\subset[0,T]$,
\[
\int_A\lambda(t)dt=0 	\quad\Longrightarrow\quad 	\mathcal L^1(s(A))=0.
\]
Indeed, let 
\[
\mu(A):=\int_A\lambda(t)\,dt
\]
and assume that $\mu(A)=0$. Given $\varepsilon>0$, by the outer 	regularity of the measure $\mu$, there exists an open set 	$O\subset[0,T]$ such that $A\subset O$ and $\mu(O)<\varepsilon$. 	Write $O$ as a countable union of pairwise disjoint intervals $I_k$ 	with endpoints $a_k,b_k$. Since $s$ is continuous and nondecreasing, 	$\overline{s(I_k)}=[s(a_k),s(b_k)]$, and hence 
\[
\mathcal L^1(s(I_k)) 	=s(b_k)-s(a_k)	=\int_{a_k}^{b_k}\lambda(t)\,dt	=\mu(I_k).
\]
Denoting by $\mathcal L^{1,*}$ the one-dimensional Lebesgue outer 	measure, we obtain
\[
\mathcal L^{1,*}(s(A))	\leq\sum_k\bigl(s(b_k)-s(a_k)\bigr)	=\sum_k\mu(I_k)	=\mu(O)	<\varepsilon.
\]
Since $\varepsilon$ is arbitrary, $\mathcal L^1(s(A))=0$. 

Now $\mu(E)=0$, because $E$ has Lebesgue measure zero, while
\[
\mu(\{\lambda=0\})	=\int_{\{\lambda=0\}}\lambda(t)\,dt	=0.
\]
It follows that
\[
\mathcal L^1(s(E))=0,	\qquad	\mathcal L^1\bigl(s(\{\lambda=0\})\bigr)=0.
\]
Let $N_\sigma\subset[0,L]$ be the set where the Lipschitz curve 	$\sigma$ is not differentiable. Then 	$\mathcal L^1(N_\sigma)=0$. Set  $N:=s(E)\cup s(\{\lambda=0\})\cup N_\sigma$. By the preceding argument, $\mathcal L^1(N)=0$. Fix 	$r\in[0,L]\setminus N$ and take $t\in[0,T]$ such that $s(t)=r$. Since 	$r\notin s(E)\cup s(\{\lambda=0\})$, necessarily $t\notin E$ and 	$\lambda(t)>0$. Moreover, $r\notin N_\sigma$, so $\sigma$ is
differentiable at $r$. Therefore, using $\gamma=\sigma\circ s$ and 	the chain rule, $\dot\gamma(t)	=\dot\sigma(r)s'(t)	=\lambda(t)\dot\sigma(r)$. Since $\lambda(t)>0$, from~\eqref{normalized-velocity} we obtain 
\[
\dot\sigma(r)	=\frac{\dot\gamma(t)}{\lambda(t)}	\in\mathcal S_{\gamma(t)}	=\mathcal S_{\sigma(r)}.
\]
Thus, $\dot\sigma(r)\in\mathcal S_{\sigma(r)}$ for a.e. 	$r\in[0,L]$, and $\sigma$ joins $x_1$ to $x_2$ in time
\[
L=\int_0^T\lambda(t)\,dt\leq T. 
\]
This proves~\eqref{Tradial} when $x_1\neq x_2$. If 	$x_1=x_2=x_0$ and $0_{x_0}\in\mathcal S_{x_0}$, the constant curve at 	$x_0$ is admissible for the  differential inclusion with $\mathcal S$	on every interval $[0,T]$, $T>0$. Hence,
$\mathcal T_{\mathcal S}(x_0,x_0) =0	=\mathcal T(x_0,x_0)$,
and the proof is complete.
\end{proof}
\begin{remark}[Time-dependent laws]\label{rem:time-dependent}
The reparametrization used in the proof of Theorem~\ref{maximalspeed}
changes the time at which each point of the trajectory is visited.
Therefore, the theorem does not extend to time-dependent velocity laws,
for which waiting may be advantageous.
For instance, let $M=\mathbb R^2$ with the Euclidean metric and
$\mathcal Y_t(u):=u+a(t)e_1$, where $e_1=(1,0)$ and
$a\colon[0,+\infty)\to[-3,0]$ is a smooth function such that
$a\equiv-3$ on $[0,1]$ and $a\equiv0$ on $[2,+\infty)$.
Let $\mathcal S_t:=\{u+a(t)e_1:|u|=1\}$ and
$\mathcal B_t:=\operatorname{rad}\mathcal S_t$, and consider the
trajectories starting from the origin at the time $t=0$.
A trajectory with $\dot\gamma(t)\in\mathcal B_t$ a.e. can rest at the
origin up to the time $t=2$ and then move with velocity $e_1$,
thus reaching the point $(1,0)$ at the time $t=3$.
On the other hand, if $\dot\gamma(t)\in\mathcal S_t$ a.e., then
$\langle\dot\gamma(t),e_1\rangle\leq1+a(t)$, which is equal to $-2$ on
$[0,1]$ and not greater than $1$ afterwards; hence
$\langle\gamma(t),e_1\rangle\leq t-3$ for every $t\geq1$,
and the point $(1,0)$ cannot be reached before the time $t=4$.
\end{remark}
Notice that
Theorem~\ref{maximalspeed} uses only the definition of
$\mathcal B_x$ as a radial hull;
in particular, it does not require convexity of $\mathcal B_x$,
or finiteness of its gauge on all of $T_xM$.
Moreover, only the continuity of $\mathcal Y$ has been used in its proof.
So we compare the original problem~\eqref{TT}
with its convexification:
let $\widehat{\mathcal T}(x_1,x_2)$
denote the infimum travel time for 
\[
	\dot\gamma(t)\in
	\widehat{\mathcal B}_{\gamma(t)}\qquad\text{a.e.}
\]
Since $\mathcal B_x\subseteq\widehat{\mathcal B}_x$
for every $x\in M$,
every admissible curve of the original differential inclusion is also admissible for the convexified one, hence
\begin{equation}\label{Tconvex-ineq}
	\widehat{\mathcal T}(x_1,x_2) 	\leq 	\mathcal T(x_1,x_2), 	\qquad x_1,x_2\in M.
\end{equation}
The reverse inequality follows
under a local controllability condition at the endpoint.

\begin{theorem}[Convexification under endpoint controllability]	\label{convexification}
	Let $x_1,x_2\in M$ and assume that 	
	\begin{equation}\label{endpoint-interiority}
		0_{x_2}\in\operatorname{int}\widehat{\mathcal B}_{x_2}.
	\end{equation}
	Then
	$ \widehat{\mathcal T}(x_1,x_2) 	= 	\mathcal T(x_1,x_2)$.

\end{theorem}

From the point of view of geometric control theory,
the above theorem is not intended as a new relaxation theorem.
Its two ingredients are standard: the Filippov--Wa\.zewski theorem
\cite[Th.~2, p.~124]{AubinCellina} provides uniform approximation of a trajectory of the convexified inclusion by trajectories of the original one with the same initial point, but does not in general preserve the terminal point; a Petrov condition (see, e.g.,~\cite[Definition 8.2.2]{CannarsaSinestrari2004}) allows this terminal error to be corrected in a time tending to zero.
Together with Theorem~\ref{maximalspeed}, this identifies the infimum travel times associated with the maximal velocities, their radial hull, and its convex hull. 

Before proving Theorem~\ref{convexification},
we need the following Hausdorff regularity result
of the sets involved.
\begin{lemma}\label{hausdorff}
In local trivializations of $TM$, the set-valued maps 
\[	x\mapsto\mathcal S_x,\qquad 	x\mapsto\mathcal B_x,\qquad	x\mapsto\operatorname{co}\mathcal S_x,\qquad
	x\mapsto\widehat{\mathcal B}_x
\]
have nonempty compact values and are locally Lipschitz with respect to 	the Hausdorff distance.
\end{lemma}
\begin{proof}
Fix a local smooth $h$-orthonormal frame on a coordinate neighborhood 	$U$.
It identifies $UTM|_U$ with $U\times\mathbb S^{n-1}$ and 	$TM|_U$ with $U\times\mathbb R^n$.
In these trivializations,
$\mathcal Y$ is represented by a map $\widetilde{\mathcal Y}	\colon 
U\times\mathbb S^{n-1}\to \mathbb R^n$,	and $\mathcal S_x=	\widetilde{\mathcal Y}(x,\mathbb S^{n-1})$.
Hence $\mathcal S_x$ is nonempty and compact. After shrinking $U$ to 	a relatively compact coordinate ball, there exists $L>0$ such that 
\[
\left|\widetilde{\mathcal Y}(x,u)	-	\widetilde{\mathcal Y}(x',u)	\right|	\leq L|x-x'|
\]
for every $x,x'\in U$ and $u\in\mathbb S^{n-1}$, where  $L$ is a Lipschitz constant of $\widetilde{\mathcal{Y}}$
on $U\times \mathbb S^{n-1}$.
Recall that, for nonempty compact subsets $A,B\subset\mathbb R^n$,
their Hausdorff distance is
\[
d_H(A,B):=
\max\left\{
\sup_{a\in A}\operatorname{dist}(a,B),\,
\sup_{b\in B}\operatorname{dist}(b,A)
\right\}.
\]
For every $v=\widetilde{\mathcal Y}(x,u)\in\mathcal S_x$, the point
$\widetilde{\mathcal Y}(x',u)$ belongs to $\mathcal S_{x'}$, and hence
\[
\operatorname{dist}(v,\mathcal S_{x'})
\leq
\left|
\widetilde{\mathcal Y}(x,u)
-
\widetilde{\mathcal Y}(x',u)
\right|.
\]
Taking the supremum over $v\in\mathcal S_x$ and then arguing
symmetrically gives
\[
d_H(\mathcal S_x,\mathcal S_{x'})
\leq
\sup_{u\in\mathbb S^{n-1}}
\left|
\widetilde{\mathcal Y}(x,u)
-
\widetilde{\mathcal Y}(x',u)
\right|
\leq L|x-x'|.
\]
It remains to observe that radialization and convexification are 	nonexpansive for the Hausdorff distance. Indeed, let $A,B\subset 	\mathbb R^n$ be nonempty and compact and set $\delta:=d_H(A,B)$.
If $ta\in\operatorname{rad}A$, with $t\in[0,1]$, choose $b\in B$ such that $|a-b|\leq\delta$. Then $tb\in\operatorname{rad}B$ and $|ta-tb|=t|a-b|\leq\delta$.
Arguing symmetrically gives
\[
d_H(\operatorname{rad}A,\operatorname{rad}B)\leq d_H(A,B).
\]
Likewise, if $a=\sum_i\lambda_i a_i\in\operatorname{co}A$, choose 	$b_i\in B$ with $|a_i-b_i|\leq\delta$. Then
$b:=\sum_i\lambda_i b_i\in\operatorname{co}B$ and 
\[
|a-b|\leq	\sum_i\lambda_i|a_i-b_i|	\leq\delta.
\]
The symmetric argument yields 
\[
d_H(\operatorname{co}A,\operatorname{co}B) 	\leq d_H(A,B).
\]
Applying these inequalities to 	$\mathcal B_x=\operatorname{rad}\mathcal S_x$ and 	$\widehat{\mathcal B}_x=\operatorname{co}\mathcal B_x$ proves the 	claim.
\end{proof}
We are now ready to prove Theorem~\ref{convexification}.

\begin{proof}[Proof of Theorem~\ref{convexification}]
In view of~\eqref{Tconvex-ineq}, it remains to prove 	$\mathcal T(x_1,x_2)\leq\widehat{\mathcal T}(x_1,x_2)$. If 	$\widehat{\mathcal T}(x_1,x_2)=+\infty$, there is nothing to prove. 	Let therefore $\gamma\colon[0,T]\longrightarrow M$ be any solution of 	the convexified inclusion such that 	$\dot\gamma(t)\in\widehat{\mathcal B}_{\gamma(t)}$  for a.e. 	$t\in[0,T]$, $\gamma(0)=x_1$, and $\gamma(T)=x_2$. 
	
By Lemma~\ref{hausdorff}, the set-valued maps 	$x\mapsto\mathcal B_x$ and $x\mapsto\widehat{\mathcal B}_x$
have nonempty compact values and are locally Lipschitz with respect to 	the Hausdorff distance in local trivializations. 
Choose then a finite subdivision $0=t_0<t_1<\cdots<t_N=T$ such that, for every $j$, $\gamma([t_{j-1},t_j])$ is contained in a coordinate neighborhood on which $\mathcal B$ and $\widehat{\mathcal B}$ are bounded and
Lipschitz with respect to the Hausdorff distance. We construct the approximating trajectory inductively. On $[t_0,t_1]$, the Filippov--Wa\.zewski relaxation theorem \cite[Th.~2, p.~124]{AubinCellina} gives a trajectory of the original
inclusion with initial point $x_1$ which is uniformly close to $\gamma$.
Suppose that the construction has been carried out up to $t_{j-1}$, and denote by $y_k$ the endpoint obtained there, so that $y_k\to\gamma(t_{j-1})$. Since the relaxation theorem preserves the initial point, it cannot be applied directly to $\gamma|_{[t_{j-1},t_j]}$, whose initial point is $\gamma(t_{j-1})$. By the continuous-dependence theorem for Lipschitz differential inclusions \cite[Th.~1, pp.~120--121]{AubinCellina}, there exists a trajectory $\widetilde\gamma_k$ of the convexified inclusion such that $\widetilde\gamma_k(t_{j-1})=y_k$ and $\widetilde\gamma_k\to\gamma$ uniformly on $[t_{j-1},t_j]$. Applying the Filippov--Wa\.zewski theorem to
$\widetilde\gamma_k$, with the same initial point $y_k$, gives a trajectory of the original inclusion uniformly close to
$\widetilde\gamma_k$, and hence to $\gamma$. Since the subdivision is finite, choosing the successive approximation errors tending to zero and concatenating the resulting trajectories yields solutions $\gamma_k\colon[0,T]\to M$ of
$\dot\gamma_k(t)\in\mathcal B_{\gamma_k(t)}$ for a.e. $t$, with $\gamma_k(0)=x_1$, such that $\gamma_k\to\gamma$ uniformly on $[0,T]$. In particular, setting $x_k:=\gamma_k(T)$, we have $x_k\to x_2$.

We next show that the terminal error can be corrected in a time 	bounded linearly by its distance from $x_2$. Let $\psi$ be a system 	of $h$-normal coordinates centred at $x_2$, associated with an 	$h$-orthonormal basis $e_1(x_2),\ldots,e_n(x_2)$ of $T_{x_2}M$. 	Thus $\psi(x_2)=0$, $d\psi_{x_2}$ is a linear isometry, and 	$|\psi(x)|=d_h(x,x_2)$ in a sufficiently small normal ball around 	$x_2$. 	Extend $e_1(x_2),\ldots,e_n(x_2)$ to a smooth local 	$h$-orthonormal frame and, for 	$\theta=(\theta^1,\ldots,\theta^n)\in\mathbb S^{n-1}$, set 	$u_\theta(x):=\sum_{i=1}^n\theta^i e_i(x)$. In these coordinates 	the original dynamics is represented by
\[
\dot z=\tilde f(z,\lambda,\theta),\qquad 	\tilde f(z,\lambda,\theta) 	:=d\psi_x\bigl(\lambda\mathcal Y_x(u_\theta(x))\bigr), 	\quad x=\psi^{-1}(z), 
\]
with $(\lambda,\theta)\in[0,1]\times\mathbb S^{n-1}$. The control set $[0,1]\times\mathbb S^{n-1}$ is compact and
$\tilde f$ is continuous and locally Lipschitz in $z$, uniformly with respect to $(\lambda,\theta)$.
Multiplying $\tilde f$ by a smooth cut-off which is identically one on a ball $\overline B_r(0)$ contained in the coordinate neighborhood, and  zero outside another ball strictly containing the first one, we obtain a globally defined bounded dynamics which is globally Lipschitz in $z$, uniformly with respect to the controls.
Thus the hypotheses 
in~\cite[Eq. (2.4)]{CannarsaCastelpietra2008} are satisfied. By~\eqref{endpoint-interiority}, there exists $\rho>0$ such that  
\[
	\rho\,\overline{\mathbb B}^{\,n} \subseteq d\psi_{x_2}\bigl(\widehat{\mathcal B}_{x_2}\bigr).
\]
Hence, for every unit vector $q\in\mathbb R^n$, \[
	\min_{v\in\mathcal B_{x_2}}\langle q,d\psi_{x_2}(v)\rangle=
	\min_{v\in\widehat{\mathcal B}_{x_2}}\langle q,d\psi_{x_2}(v)\rangle \leq-\rho,
\]
where the equality follows from $\widehat{\mathcal B}_{x_2}=\operatorname{co}\mathcal B_{x_2}$. 

Since $\tilde f\bigl(0,[0,1]\times\mathbb S^{n-1}\bigr) =d\psi_{x_2}(\mathcal B_{x_2})$, the preceding inequality and the positive homogeneity of both sides with respect to $\nu$ give 
\[
	\min_{(\lambda,\theta)\in[0,1]\times\mathbb S^{n-1}}\langle\tilde f(0,\lambda,\theta),\nu\rangle
	\leq-\rho|\nu|, \qquad \nu\in\mathbb R^n.
\]
For the target $K=\{0\}$,  one has  $\partial K=\{0\}$ and  the normal cone $N_K(0)=\mathbb R^n$. Since the system is autonomous, this is precisely the Petrov condition~\cite[(2.5)]{CannarsaCastelpietra2008}.
We apply \cite[Prop.~2.10]{CannarsaCastelpietra2008}  with $\Omega=\mathbb R^n$, so that the problem is unconstrained. The cut-off dynamics satisfies the assumptions in Eq. (2.4) of that paper. Therefore there exist $\delta>0$ and $C_P>0$ such that its minimum-time function $\widetilde T$ satisfies $\widetilde T(z)\leq C_P|z|$ whenever $|z|<\delta$. For $0<|z|<\delta$, by the definition of the infimum there exists a trajectory from $z$ to $0$ with travel time $\tau\leq\widetilde T(z)+|z|\leq(C_P+1)|z|$. Let $L$ be a bound for the cut-off dynamics. Along such a trajectory,
\[ 
	|y(t)|\leq |z|+Lt \leq\bigl(1+L(C_P+1)\bigr)|z|, \qquad 0\leq t\leq\tau.
\]
Thus, after decreasing $\delta$ so that $\bigl(1+L(C_P+1)\bigr)\delta<r$, and since  the cut-off is identically one on $B_r(0)$, every such trajectory is in fact a trajectory of the original system. Since $|\psi(x)|=d_h(x,x_2)$, we obtain a
neighborhood $V$ of $x_2$ such that 
\begin{equation}\label{Petrov-time-estimate}
\mathcal T(x,x_2)\leq(C_P+1)d_h(x,x_2), 	\qquad x\in V.
\end{equation}
For all sufficiently large $k$, $x_k\in V$. By \eqref{Petrov-time-estimate}, there exists a trajectory of the
original system joining $x_k$ to $x_2$ in a time $\tau_k$ satisfying $\tau_k\leq(C_P+1)d_h(x_k,x_2)+1/k$.
Concatenating this trajectory with 	$\gamma_k$ gives 	$\mathcal T(x_1,x_2)\leq 	T+(C_P+1)d_h(x_k,x_2)+1/k$. Letting $k\to\infty$ yields 	$\mathcal T(x_1,x_2)\leq T$. Since $\gamma$ was an arbitrary 	admissible trajectory for the convexified problem, taking the 	infimum over all such trajectories gives 	$\mathcal T(x_1,x_2)\leq\widehat{\mathcal T}(x_1,x_2)$.  	Together with~\eqref{Tconvex-ineq}, this proves the claim.
\end{proof}
Combining Theorems~\ref{maximalspeed} and~\ref{convexification}, for $x_1\neq x_2$ and $0_{x_2}\in\operatorname{int}\widehat{\mathcal B}_{x_2}$,  the three velocity sets
\[
\mathcal S_x=\mathcal Y_x(UT_xM),
\qquad
\mathcal B_x=\operatorname{rad}\mathcal S_x,
\qquad
\widehat{\mathcal B}_x=\operatorname{co}\mathcal B_x
\]
give the same infimum travel time:
\begin{equation}\label{three-travel-times}
\mathcal T_{\mathcal S}(x_1,x_2)=\mathcal T(x_1,x_2)=\widehat{\mathcal T}(x_1,x_2).
\end{equation}
Let us finally relate these general facts with the positive-self-propulsion
regime considered above. Suppose that $\eta>0$ and $b<1$. Then
$\mathcal X_x(UT_xM)$ is a radial graph and, if $-W_x\in\operatorname{int}\mathcal M_x$,
the compact set $\mathcal J_x:=\eta(x)\bigl(\mathcal M_x+W_x\bigr)$ contains $0_x$ in its interior. Moreover, $\partial\mathcal J_x=\mathcal Y_x(UT_xM)$, and Lemma~\ref{radialhull}-(i) gives $\mathcal J_x\subseteq\operatorname{rad}\partial\mathcal J_x=\mathcal B_x$.
Hence $0_x\in\operatorname{int}\mathcal B_x$, and \eqref{three-travel-times} applies. In particular, this includes the range $1/2\leq b<1$, where the Matsumoto gauge has already ceased to be a Finsler metric (and, for $b>1/2$, is nonconvex).

On the other hand, Theorem~\ref{maximalspeed} remains valid also at   points of $M_+$ where $b_x\geq1$,  and   where $\eta=0$. In such cases the equality with the convexified problem still follows whenever~\eqref{endpoint-interiority} holds.

The preceding results compare the infimum travel times associated with $\mathcal S$, $\mathcal B$, and, under the endpoint controllability condition, $\widehat{\mathcal B}$, but they do not by themselves imply that the relevant infimum is attained. Indeed, when $\mathcal B_x$ is nonconvex, a minimizing sequence may converge only to a trajectory of the convexified inclusion.

A closely related  existence result for autonomous, possibly nonconvex differential inclusions was proved in  \cite[Th.~2.2, p.~121]{CellinaFerrieroMarchini2006}. The assumption $(ii)$ in the statement there reads in our setting:
\[
\text{for every }0_x\neq\xi\in\operatorname{co}\mathcal S_x
\text{ there exists }\lambda\geq1
\text{ such that }\lambda\xi\in\mathcal S_x,
\]
which is then  equivalent to the convexity of $\mathcal B_x=\operatorname{rad}\mathcal S_x$. Indeed, if
$\mathcal B_x$ is convex, then $\operatorname{co}\mathcal S_x\subseteq\mathcal B_x$, and every nonzero
$\xi\in\mathcal B_x$ can be written as $\xi=t z$, with $z\in\mathcal S_x$ and $t\in(0,1]$; hence
$t^{-1}\xi=z\in\mathcal S_x$. Conversely, it  implies $\operatorname{co}\mathcal S_x\subseteq\mathcal B_x$. Since
$\mathcal B_x$ is radial, it follows that  
\[
\operatorname{co}\mathcal B_x
=\operatorname{co}(\operatorname{rad}\mathcal S_x)
=\operatorname{co}(\{0_x\}\cup\mathcal S_x)
=\operatorname{rad}\bigl(\operatorname{co}\mathcal S_x\bigr)
\subseteq\mathcal B_x,
\]
and therefore $\mathcal B_x$ is convex.

Under the hypotheses below, however, existence can be obtained directly
from the classical compactness theorem for convex-valued differential
inclusions, applied to $x\mapsto\mathcal B_x$. The reparametrization
constructed in the proof of Theorem~\ref{maximalspeed} then promotes an
optimal $\mathcal B$-trajectory to an optimal
$\mathcal S$-trajectory.

\begin{proposition}[Existence of optimal trajectories]\label{nonfinsler-existence}
Let $(M,h)$ be complete and assume that 	
\begin{equation}\label{bounded-maximal-speed}
C_{\mathcal Y}:=\sup\bigl\{	\alpha\bigl(\mathcal Y_x(u)\bigr):	x\in M,\ u\in UT_xM\bigr\}	<+\infty.
\end{equation}
Suppose moreover that $\mathcal B_x$ is convex for every $x\in M$. 	Then, for every $x_1\neq x_2$ such that
$\mathcal T(x_1,x_2)<+\infty$, there exists $\sigma\in	AC\bigl([0,\mathcal T(x_1,x_2)],M\bigr)$
such that
\[
\sigma(0)=x_1,\qquad	\sigma\bigl(\mathcal T(x_1,x_2)\bigr)=x_2,\qquad	\dot\sigma(t)\in\mathcal S_{\sigma(t)}
	\quad\text{a.e.}
\]
\end{proposition}
\begin{proof}
Set $T_*:=\mathcal T(x_1,x_2)$.	Since every $v\in\mathcal B_x$ satisfies	$\alpha(v)\leq C_{\mathcal Y}$, the assumptions	$x_1\neq x_2$ and $T_*<+\infty$ imply $T_*\geq d_h(x_1,x_2)/C_{\mathcal Y}>0$.
For every $k\geq1$, choose a $\mathcal B$-admissible trajectory
$\gamma_k\colon[0,T_k]\to M$ such that
\[
\gamma_k(0)=x_1,\qquad	\gamma_k(T_k)=x_2,\qquad	T_*\leq T_k<T_*+\frac1k.
\]
Set $T_0:=T_*+1$ and extend $\gamma_k$ to $[0,T_0]$ by putting 	$\gamma_k(t):=x_2$ for $t\in[T_k,T_0]$. This extension is still 	$\mathcal B$-admissible because $0_x\in\mathcal B_x$ for every $x$. 	The curves $\gamma_k$ are $C_{\mathcal Y}$-Lipschitz and their images 	are contained in
\[
K:=\overline{B_h\bigl(x_1,C_{\mathcal Y}T_0\bigr)},
\]
which is compact by  completeness of $h$. 	Moreover, Lemma~\ref{hausdorff} shows that $x\mapsto\mathcal B_x$ is locally Lipschitz with respect to the 	Hausdorff distance. Its values are nonempty and compact, and they are
convex by hypothesis.
Consequently, the compactness theorem for convex-valued differential inclusions
\cite[Th.~1, p.~104]{AubinCellina} (localized to $K$	using a finite  covering by charts of $M$) yields a subsequence, still denoted by $\gamma_k$,
and a curve $\gamma\in AC([0,T_0],M)$ such that
$\gamma_k\to \gamma$ uniformly on $[0,T_0]$, and $\dot\gamma(t)\in\mathcal B_{\gamma(t)}$, a.e.
Since $T_*\leq T_k$ and $\gamma_k(T_k)=x_2$, we have
\[
d_h\bigl(\gamma_k(T_*),x_2\bigr)\leq	\int_{T_*}^{T_k}\alpha\bigl(\dot\gamma_k(t)\bigr)\,dt	\leq
	C_{\mathcal Y}(T_k-T_*)	\longrightarrow0.
\]
Uniform convergence therefore gives $\gamma(T_*)=x_2$. Thus the 	restriction of $\gamma$ to $[0,T_*]$ is an optimal 	$\mathcal B$-trajectory from $x_1$ to $x_2$. 	Applying the reparametrization constructed in the proof of
Theorem~\ref{maximalspeed} to this trajectory gives an 	$\mathcal S$-trajectory $\sigma\colon[0,L]\to M$ with the same 	endpoints and with $L\leq T_*$. Since 	$\mathcal S_x\subseteq\mathcal B_x$, the definition of $T_*$ also
implies $T_*\leq L$. Hence $L=T_*$, and $\sigma$ has all the stated properties.
\end{proof}
Condition~\eqref{bounded-maximal-speed} is automatic when $M$ is compact. It is also implied by the boundedness of $\eta$, $\widehat Z$, and $\widehat W$, since
\[
\alpha\bigl(\mathcal Y_x(u)\bigr)\leq\eta(x)+\alpha(\widehat Z_x)+\alpha(\widehat W_x).
\]
For the velocity law~\eqref{eq:generalForm}, the convexity of
$\mathcal B_x$ can be read off from the data, and it holds well beyond
the regime in which $F$ is a Finsler metric.
\begin{corollary}[Existence for the two-field law]\label{cor:model-existence}
Let $\mathcal Y$ be given by~\eqref{eq:generalForm}.
Assume that $(M,h)$ is complete, that $\eta$, $\widehat Z$ and $\widehat W$ are bounded, and that
\begin{equation}\label{halfcondition}
\alpha(\widehat Z_x)\leq\frac12\,\eta(x)\qquad\text{for every }x\in M_+.
\end{equation}
Then every reachable pair $x_1\neq x_2$ is joined by a time-optimal trajectory having  maximal velocity almost everywhere.
In particular, this holds when $\eta\equiv0$.
\end{corollary}
\begin{proof}
Let us show that $\mathcal B_x$ is convex for every $x\in M$.
Let first $\eta(x)=0$.
If $\widehat Z_x\neq0_x$, Lemma~\ref{spheregauge} shows that 	$\mathcal B_x$ is convex. If $\widehat Z_x=0_x$, then 	$\mathcal S_x=\{\widehat W_x\}$ and 
$\mathcal B_x=\operatorname{rad}\{\widehat W_x\}$ is a segment, hence 	is again convex.
Let now $x\in M_+$, so that $b_x\leq1/2$ by~\eqref{halfcondition}.
Then $\mathcal M_x$ is convex: for $b_x<1/2$ this follows from
Proposition~\ref{prop:FinslerM-conditions}, while for $b_x=1/2$ the set
$\mathcal M_x$ is the limit, in the Hausdorff distance, of the convex sets
obtained by replacing $Z_x$ with $sZ_x$, $s\nearrow1$.
Since $b_x<1$, \eqref{radialgraph} gives
$\mathcal S_x=\eta(x)\bigl(\partial\mathcal M_x+W_x\bigr)$ and then,
by Lemma~\ref{radialhull}-(i),
\[
\mathcal B_x
=\eta(x)\operatorname{rad}\bigl(\partial(\mathcal M_x+W_x)\bigr)
=\eta(x)\operatorname{rad}\bigl(\mathcal M_x+W_x\bigr)
=\eta(x)\operatorname{co}\bigl(\{0_x\}\cup(\mathcal M_x+W_x)\bigr),
\]
where the last equality holds because $\mathcal M_x+W_x$ is convex.
Hence $\mathcal B_x$ is convex also in this case.
Moreover,
$\alpha\bigl(\mathcal Y_x(u)\bigr) 	\leq	\eta(x)+\alpha(\widehat Z_x)+\alpha(\widehat W_x)$ is bounded, so that Proposition~\ref{nonfinsler-existence} applies.
\end{proof}
Notice that no restriction is imposed on $\widehat W$, apart from
its boundedness: Corollary~\ref{cor:model-existence} applies also when
$F_{\mathcal M}(-W_x)\geq1$ at some point $x\in M_+$,
where $F$ is not a Finsler metric.
\begin{remark}[Comparison with wind Finslerian structures]
Let $\Sigma$ be a wind Finslerian structure in the sense of 	\cite[Def.~2.8]{CaponioJavaloyesSanchez2024}, and let $\mathcal K_x$ 	be the compact strongly convex body bounded by $\Sigma_x$. Choosing 	$\mathcal Y$ so that $\mathcal S_x=\Sigma_x$, one has 
\[
\mathcal B_x=\operatorname{rad}\Sigma_x	=\operatorname{co}\bigl(\{0_x\}\cup\mathcal K_x\bigr),
\]
and hence $\mathcal B_x$ is convex, independently of the position of 	$0_x$ with respect to $\mathcal K_x$. Consequently, 	Proposition~\ref{nonfinsler-existence} applies to every such structure 	using any  auxiliary complete  Riemannian  $h$ such that the 	velocities satisfy \eqref{bounded-maximal-speed} with respect to the norm $\alpha$ associated with $h$. Notice that this also 	includes the strong-wind case, where $0_x$ need not belong to 	$\mathcal K_x$.
	
In particular, for a wind Riemannian structure whose indicatrix is the 	$h_x$-sphere of radius $r_x$ centred at $C_x$ (see Subsection~\ref{subsec:zermelo}), consider 	$\widetilde h:=h/r^2$. Then 	$\widetilde\alpha(\mathcal Y_x(u)) 	\leq1+\alpha_x(C_x)/r_x$, so that  	Proposition~\ref{nonfinsler-existence} applies whenever 	$\widetilde h$ is complete and $\alpha(C)/r$ is bounded. In the mild 	and critical regimes $\alpha(C)\leq r$, the latter condition is 	automatic. For a Randers--Kropina structure associated with 	$(g_0,\omega,\Lambda)$, $\Lambda\geq0$, this gives
\[
\widetilde h 	=\frac{g_0}{\Lambda+\alpha^2(\omega^\sharp)}.
\]
Thus the completeness of this metric is a simple sufficient condition 	for the existence of a time-optimal trajectory between any two distinct  	reachable points.  For related existence results for minimizing wind and
Randers--Kropina geodesics using Lorentzian methods  see \cite[Th.~6.15]{CaponioJavaloyesSanchez2024} and
\cite[Rem.~5.4]{cjmRK} (as corrected in \cite[p. 835]{cjmRKcorr}); in the Kropina case see also 
\cite[Prop.~4.1, Cor.~4.3, Th. 4.4]{CaponioGiannoniMasielloSuhr2021}. 
The multiplicity result of \cite[Th.~5.1]{cjmRK} (with assumption $(c)$  corrected in \cite[p.834]{cjmRKcorr}) is  of a different nature and is not addressed by Proposition~\ref{nonfinsler-existence}.
\end{remark}

\subsection{Completeness of \texorpdfstring{$F$}{F} within $M_+$ and time-optimal trajectories}
\label{subsec:FM+completeness}

In this subsection, we study the completeness of $M_+ = \{x \in M: \eta(x) > 0\}$ with respect to $F:TM_+ \to [0,+\infty)$,
assuming that~\eqref{eq:FinslerM+conditions} holds, so that $F$ is a Finsler metric on $M_+$.
Some care is needed here, because the set
on which $F$ is a Finsler metric need not be all of $M$.
More precisely,  if $\partial M_+\neq\emptyset$
then $(M_+,h)$ is never complete, whether or not $(M,h)$ is.
Thus,  the completeness of the 	background metric $h$ must be assumed on $M$, and the behaviour of $F$ near the zero set $\{x \in M: \eta(x)=0\}$ of $\eta$ must be controlled 	separately.
That control is provided by the following estimate, whose only ingredient is the local Lipschitz estimate $\eta(x)\leq C d_h(x,x_0)$ near a zero $x_0$ of the non-negative smooth function $\eta$.
This upper bound is sufficient to force infinite travel time to the  set $\{\eta=0\}$.
\begin{lemma}\label{loglemma}
	Assume that~\eqref{eq:FinslerM+conditions} holds,
	and let $x_0\in M$ be such that $\eta(x_0)=0$. Let $r>0$ be such that $\overline{B_h(x_0,r)}$ is compact and contained in a convex neighborhood (w.r.t. $h$) of $x_0$, and assume $B_h(x_0,r)\cap M_+\neq\emptyset$. Then
\[
C:=\sup_{B_h(x_0,r)}|d\eta|_h\in(0,+\infty),
\]
and, for every $\gamma\in AC\bigl([t_0,t_1],M_+\cap B_h(x_0,r)\bigr)$,
\begin{equation}\label{logestimate}
\ell_F(\gamma)\ \geq\ \frac{1}{3C}\left|\,\log\frac{d_h\bigl(\gamma(t_1),x_0\bigr)}{d_h\bigl(\gamma(t_0),x_0\bigr)}\,\right|.
\end{equation}
\end{lemma}
\begin{proof}
The quantity $C$ is finite because $\overline{B_h(x_0,r)}$ is compact and $\eta$ is smooth. If $C=0$, then $\eta$ would be constant on the connected set $B_h(x_0,r)$, hence identically $0$ there, contradicting $B_h(x_0,r)\cap M_+\neq\emptyset$; thus $C>0$.

Let $x\in B_h(x_0,r)$ and let $c$ be the minimizing $h$-geodesic from $x_0$ to $x$, which is contained in $B_h(x_0,r)$ because every point of $c$ is at $h$-distance at most $d_h(x,x_0)<r$ from $x_0$. Applying the mean value theorem to $\eta\circ c$ and using $\eta(x_0)=0$ we obtain
\begin{equation}\label{linearvanishing}
\eta(x)\leq C\,d_h(x,x_0),\qquad x\in B_h(x_0,r).
\end{equation}
Set now $\ell(t):=d_h\bigl(\gamma(t),x_0\bigr)$. Since $\gamma(t)\in M_+$ we have $\gamma(t)\neq x_0$,
so $\ell>0$. Moreover, for $t_0\leq s\leq t\leq t_1$,
\[
\bigl|\ell(t)-\ell(s)\bigr|\leq d_h\bigl(\gamma(s),\gamma(t)\bigr)\leq\int_s^t\alpha_{\gamma}(\dot\gamma)\,d\tau.
\]
As $\alpha(\dot\gamma)\in L^1(t_0,t_1)$, the function $\ell$ is
absolutely continuous and $|\ell'|\leq\alpha(\dot\gamma)$ a.e.; hence, by the lower bound in~\eqref{eq:F-LowUpp-bounds} and by~\eqref{linearvanishing},
\[
\ell_F(\gamma)=\int_{t_0}^{t_1}F(\gamma,\dot\gamma)\,dt
\geq\int_{t_0}^{t_1}\frac{\alpha_{\gamma}(\dot\gamma)}{3\eta(\gamma)}\,dt
\geq\frac{1}{3C}\int_{t_0}^{t_1}\frac{|\ell'|}{\ell}\,dt
\geq\frac{1}{3C}\left|\int_{t_0}^{t_1}\frac{\ell'}{\ell}\,dt\right|,
\]
which is~\eqref{logestimate}. 
\end{proof}
\begin{theorem}\label{bfcompleteness}
	Assume that~\eqref{eq:FinslerM+conditions} holds.
Suppose that $(M,h)$ is complete and that there exists a constant $\delta>0$ such that $\eta(x)\leq\delta$ for all $x\in M$. Then $(M_+,F)$ is both forward and backward complete.
\end{theorem}
\begin{proof}
By~\eqref{eq:F-LowUpp-bounds} and the upper bound on $\eta$,
\begin{equation}\label{globallowerbound}
	F(x,y)\geq \frac{\alpha_x(y)}{3\eta(x)}\geq \frac{\alpha_x(y)}{3\delta},
	\qquad \forall\, (x,y)\in TM_+. 
\end{equation}
Both~\eqref{globallowerbound} and~\eqref{logestimate} bound from below the $F$-length of \emph{any} curve joining two given points of $M_+$, irrespective of its orientation. Consequently, each of the lower bounds derived below applies simultaneously to $d_F(p,q)$ and to $d_F(q,p)$, and it suffices to argue for a sequence $(x_k)\subset M_+$ that is forward  Cauchy for $d_F$ since the other case can be proved analogously. 

By~\eqref{globallowerbound}, every forward Cauchy sequence $(x_k)$ for $d_F$ is 
 Cauchy also for $d_h$ and, $(M,h)$ being complete, $x_k\to x_\infty$ for some $x_\infty\in M$.
Let us see that  the limit lies in $M_+$. Suppose, by contradiction, that $\eta(x_\infty)=0$. Choose $r>0$ as in Lemma~\ref{loglemma} for $x_0=x_\infty$; since $x_k\in B_h(x_\infty,r)\cap M_+$ for $k$ large, the lemma applies and yields a constant $C>0$. Set $r_k:=d_h(x_k,x_\infty)>0$, so that $r_k\to0$, and put
\[
\varepsilon_0:=\frac12\min\Bigl\{\frac{r}{6\delta},\,1\Bigr\}>0 .
\]
Let $k$ be large enough that $r_k<r/2$  and let $j\geq k$ with $r_j<r/2$. Let $\gamma$ be any curve in $M_+$ joining $x_k$ and $x_j$. Then one of the following occurs:
\begin{itemize}
\item[(i)] $\gamma$ is contained in $B_h(x_\infty,r)$. Then~\eqref{logestimate} gives
$\ell_F(\gamma)\geq\frac{1}{3C}\bigl|\log(r_j/r_k)\bigr|$.
\item[(ii)] $\gamma$ leaves $B_h(x_\infty,r)$. Then its $h$-length is at least $r-\max\{r_k,r_j\}>r/2$, so that  $\ell_F(\gamma)\geq r/(6\delta)$ by~\eqref{globallowerbound}.
\end{itemize}
Taking the infimum over all such $\gamma$ we obtain
\[
d_F(x_k,x_j)\geq\min\left\{\frac{1}{3C}\left|\log\frac{r_j}{r_k}\right|,\ \frac{r}{6\delta}\right\}.
\]
Since $r_k\to0$, we may choose $j\geq k$ with $r_j<r_k\,e^{-6C\varepsilon_0}$, so that the first term is at least $2\varepsilon_0$; the second is at least $2\varepsilon_0$ by the definition of $\varepsilon_0$. Therefore,
$d_F(x_k,x_j)>\varepsilon_0$, contradicting the assumption that $(x_k)$ is a forward Cauchy sequence. Hence $x_\infty\in M_+$.
On $M_+$ the metric $F$ is a Finsler metric, so the $h$-balls and the forward and backward balls of $d_F$ all generate the manifold topology (see \cite[\S 6.2.C]{bcs}). Thus,  $x_k\to x_\infty\in M_+$ in that topology, and a forward (resp. backward) Cauchy sequence converging in the topology converges with respect to  the topology generated by forward balls of $d_F$, i.e.   $(M_+,F)$ is forward complete.
\end{proof}
\begin{remark}
Only the local Lipschitz continuity of $\eta$ has been used in the proof, through \eqref{linearvanishing}.  If one weakens this to continuity, the conclusion can fail. For example, let $M=\R^n$ with the Euclidean metric, $\widehat Z=\widehat W=0$ and
\[
\eta(x)=\min\bigl\{\sqrt{|x^1|},\,1\bigr\},
\]
which is continuous but not Lipschitz across $\{x^1=0\}$. Here $b\equiv0$, $F_{\mathcal M}(-W)\equiv0$ and $F=\alpha/\eta$, so all the remaining hypotheses of Theorem~\ref{bfcompleteness} hold with $\delta=1$. Nevertheless, the segment $t\mapsto(1-t,0,\ldots,0)$, $t\in[0,1)$, has finite $F$-length $\int_0^1s^{-1/2}\,ds=2$ and leaves every compact subset of $M_+$, so $(M_+,F)$ is not complete. 
\end{remark}
The estimate~\eqref{logestimate} is local: it requires neither the completeness of $(M,h)$ nor any bound on $\eta$. Hence, the following invariance property holds without the global assumptions of Theorem~\ref{bfcompleteness}.
\begin{proposition}\label{unreachable}
Assume that~\eqref{eq:FinslerM+conditions} holds. Then the  set $\{\eta=0\}$ is at infinite travel time from $M_+$. More precisely, every curve $\gamma\in AC([0,T],M)$ which is admissible for~\eqref{TT} with $\gamma(0)\in M_+$ is contained in $M_+$.
\end{proposition}
\begin{proof}
Let $\gamma\in AC([0,T],M)$ be admissible for~\eqref{TT} with $\gamma(0)\in M_+$, assume by contradiction that $\gamma$ is not contained in $M_+$, and let $\tau:=\inf\{t\in[0,T]:\gamma(t)\notin M_+\}$. Then $\tau>0$ and $\gamma(\tau)\in\partial M_+$. Let $r$ and $C$ be as in Lemma~\ref{loglemma} for $x_0=\gamma(\tau)$. Choose $t<\tau$ sufficiently close to $\tau$ so that $\gamma([t,\tau])\subset B_h(\gamma(\tau),r)$. Then, for every $t'\in(t,\tau)$, \eqref{logestimate} gives 
\[
\ell_F(\gamma|_{[t,t']})\geq\frac{1}{3C}\left|\log\frac{d_h(\gamma(t'),\gamma(\tau))}{d_h(\gamma(t),\gamma(\tau))}
\right|.
\]
Since $d_h(\gamma(t'),\gamma(\tau))\to0$ as $t'\nearrow\tau$, it follows that $\ell_F(\gamma|_{[t,\tau)})=+\infty$.
On the other hand, admissibility gives $F(\dot\gamma)\leq1$ a.e. on $[0,\tau)$, and hence 
\[
\ell_F(\gamma|_{[0,\tau)})\leq\tau<+\infty,
\]
a contradiction. Hence every admissible trajectory issuing from $M_+$ remains in $M_+$ for all time. 
\end{proof}
\begin{remark}[Rest points]\label{rem:restpoints}
Proposition~\ref{unreachable} is an instance of a general fact, which does not depend on the Finsler structure. Let $\mathcal Y$ be a locally Lipschitz velocity law and let $x_0\in M$ be a rest point, i.e. $\mathcal S_{x_0}=\{0_{x_0}\}$. Arguing as in the proof of Lemma~\ref{hausdorff}, there exist a neighbourhood $U$ of $x_0$ and $L>0$ such that $\alpha(v)\leq L\,d_h(x,x_0)$ for every $x\in U$ and $v\in\mathcal B_x$. Hence, along every admissible curve $\gamma$ contained in $U$, the function $\ell(t):=d_h(\gamma(t),x_0)$ satisfies $|\ell'|\leq\alpha(\dot\gamma)\leq L\ell$ a.e. and, by Gronwall's inequality, it either vanishes identically or never vanishes. Thus, a rest point can be neither reached nor left in finite time. Under~\eqref{eq:FinslerM+conditions}, every point of $\partial M_+$ is a rest point of~\eqref{eq:generalForm} (see Remark~\ref{rem:gluing}); therefore, no admissible trajectory crosses $\partial M_+$, in either direction.
\end{remark}
 If $F$ is forward or backward complete, the Finslerian Hopf--Rinow theorem ensures that the distance $d_F(x_1,x_2)$ is attained by a minimizing geodesic
from $x_1$ to $x_2$ (see, e.g.,~\cite[Th.~6.6.1]{bcs}); parametrized with unit $F$-speed, such a geodesic
is then a smooth solution of  problem \eqref{TT} as stated in the following corollary.
\begin{corollary}\label{cor:geodesics}
Assume that~\eqref{eq:FinslerM+conditions} holds. Then,  for every $x_1,x_2$ in the same connected component of $M_+$,  $\mathcal T(x_1,x_2)$ in \eqref{TT} is  equal to $d_F(x_1,x_2)$. If, moreover, the assumptions 	of Theorem~\ref{bfcompleteness} hold and $x_1\neq x_2$, then
$\mathcal T(x_1,x_2)>0$ is attained by an $F$-unit-speed  geodesic from $x_1$ to $x_2$.
\end{corollary}
\begin{proof}
By Proposition~\ref{unreachable}, every curve admissible for~\eqref{TT} and issuing from $x_1\in M_+$ stays in $M_+$, so we may argue entirely on $M_+$. There $\mathcal B_x$ is the closed unit ball of $F(x,\cdot)$, so a curve 	$\gamma\in AC([0,T],M)$ is admissible for \eqref{TT} if and only if 	$F(\dot\gamma)\leq1$ a.e.; hence if $\gamma\in AC([0,T],M)$ is admissible and connects $x_1$ to $x_2$, we get $T\geq\int_0^TF(\dot\gamma)\,dt=\ell_F(\gamma) \geq d_F(x_1,x_2)$, so that $\mathcal T(x_1,x_2)\geq d_F(x_1,x_2)$. Conversely, the equality being trivial for $x_1=x_2$, let $x_1\neq x_2$ and let $c$ be a piecewise smooth curve in $M_+$ from $x_1$ to $x_2$. After an affine change of parameter, which does not affect its $F$-length, we can assume that $F(\dot c)\leq1$, i.e. that $c$ is admissible for~\eqref{TT}. The reparametrization constructed in the proof of Theorem~\ref{maximalspeed} then yields an admissible curve connecting $x_1$ to $x_2$ in the time $\int F(\dot c)\,dt=\ell_F(c)$; taking the infimum over all such curves, we get $\mathcal T(x_1,x_2)\leq d_F(x_1,x_2)$, hence equality 	holds. Finally, under the assumptions of Theorem~\ref{bfcompleteness}, 	by the same result and the Hopf--Rinow theorem there exists a 	minimizing $F$-geodesic in $M_+$ from $x_1$ to $x_2$ which, parametrized with unit 	$F$-speed, is defined on $[0,d_F(x_1,x_2)]$ and is admissible, hence the infimum is attained.
\end{proof}

\section{Support Hamiltonian}\label{hamiltonian}

The time-optimal problem \eqref{TT} admits a dual description in terms of the
\emph{support function} of the admissible velocity sets, that is, the function
$H\colon T^*M\to\mathbb R$ defined by
\[
H(p):=\max_{z\in\mathcal B_x}p(z), \qquad p\in T_x^*M.
\]
Indeed, for a fixed $x_0\in M$, the minimum-time function $\tau_{x_0}:=\mathcal T(x_0,\cdot)$
increases at most at unit rate along admissible trajectories, and exactly at
unit rate along time-minimizing ones. Hence $\tau_{x_0}$ satisfies the eikonal equation
$H(d\tau_{x_0})=1$ wherever it is differentiable and, under suitable controllability
assumptions, in the viscosity sense (see, e.g.,~\cite{CannarsaSinestrari2004}).

For every $x\in M$, the restriction of $H$ to $T_x^*M$ is convex and positively
homogeneous, regardless of whether $\mathcal B_x$ is convex. Moreover, $H\geq0$
because $0_x\in\mathcal B_x$. Since a linear functional has the same maximum on
a set and on its convex hull, we also have
\[
H(p)=\max_{z\in\operatorname{co}\mathcal B_x}p(z).
\]
Thus, when $\mathcal B_x$ is nonconvex, $H$ does not distinguish the original
differential inclusion from its convexification.
This is consistent with
Theorem~\ref{convexification}: under the hypotheses
of Section~\ref{sec:time-travel}, this replacement does not affect the
infimum travel time.
Finally, unlike the gauge
$F$, which is a Finsler metric only in the regular regime, $H$ can be computed directly
from the velocity law $\mathcal Y_x$ given by \eqref{eq:generalForm} for every $\eta\geq0$, as we now show.

We denote by $\alpha^*$ the norm associated with the dual Riemannian metric $h^*$. Whenever $\mathcal B_x$ is the unit ball of a Finsler norm $F$, its dual norm is denoted by $F^*$.

\begin{proposition}\label{support-general}
For every $x\in M$ and $p\in T_x^*M$,
\begin{equation}\label{Hgeneral-hat}
H(p)=\max\left\{0,\, p(\widehat W_x)+\max_{u\in UT_xM}\bigl(\eta(x)+h(\widehat Z_x,u)\bigr)p(u)	\right\}.
\end{equation}
In particular, at a point where $\eta(x)=0$, 
\begin{equation}\label{H-zero-eta}
H(p)=\max\left\{0,\,p(C_x)+r_x\alpha^*(p)\right\}=\max\left\{0,\,p(\widehat W_x)
		+\frac{	p(\widehat Z_x)+\alpha(\widehat Z_x)\alpha^*(p)}{2}	\right\},
\end{equation}
where $C_x$ and $r_x$ are defined in~\eqref{RC}, here extended to the case $\widehat Z_x=0_x$ (so that $r_x=0$).
If instead $\eta(x)>0$, $b_x<1/2$ and 	$F_{\mathcal M}(-W_x)<1$, then $\mathcal B_x$ is the unit ball of a
Finsler norm $F$, and
\begin{equation}\label{Hgen}
H(p)=F^*(p)	=\eta(x)\bigl(p(W_x)+H_Z(p)\bigr),
\end{equation}
where
\begin{equation}\label{HZ}
H_Z(p):=\max_{u\in UT_xM}\bigl(1+h(Z_x,u)\bigr)p(u)	=F_{\mathcal M}^*(p).
\end{equation}
\end{proposition}
\begin{proof}
From~\eqref{admissible}, we get
\[
H(p)=\max_{\substack{u\in UT_xM\\0\leq\lambda\leq1}}\lambda p\bigl(\mathcal Y_x(u)\bigr)=
		\max\left\{	0,\,\max_{u\in UT_xM}p\bigl(\mathcal Y_x(u)\bigr)\right\},
\]
which gives~\eqref{Hgeneral-hat} after substituting	\eqref{eq:generalForm}. The maximum with zero takes into account exactly 	the radial hull in~\eqref{admissible}.
	
Suppose now that $\eta(x)=0$, so that $\mathcal Y_x$ is given by \eqref{eq:zero-eta-law}.  By~\eqref{eq:zero-eta-navigation},   the maximal velocity set is the possibly degenerate sphere
$\mathcal Y_x(UT_xM)=r_xUT_xM+C_x$, whose  support function is then $p\in T^*_xM\mapsto p(C_x)+r_x\alpha^*(p)$, and the second equality in \eqref{H-zero-eta} follows from~\eqref{RC}.
If $\widehat Z_x=0_x$, then $r_x=0$ and the same formula reduces to $H(p)=\max\{0,p(\widehat W_x)\}$,	consistently with $\mathcal B_x=[0_x,\widehat W_x]$.
	
Finally, assume that $\eta(x)>0$, $b_x<1/2$ and 	$F_{\mathcal M}(-W_x)<1$. By Proposition~\ref{generalcase}, 	$0_x\in\operatorname{int}\mathcal B_x$ and $\mathcal B_x$ is the unit ball of the Finsler metric $F$. Hence the  maximum in~\eqref{Hgeneral-hat} is redundant 	and $H=F^*$ on $T_x^*M$. Substituting 	in \eqref{Hgeneral-hat} $\widehat Z=\eta Z$ and $\widehat W=\eta W$ gives \eqref{Hgen}. Since $\mathcal M_x$ is the unit ball of 	$F_{\mathcal M}$ and $\partial\mathcal M_x=\left\{\bigl(1+h(Z_x,u)\bigr)u:u\in UT_xM\right\}$,
its support function is precisely~\eqref{HZ}.
\end{proof}

By Lemma~\ref{separation}, the  maximum in
\eqref{Hgeneral-hat} is redundant for every covector if and only if $0_x\in\operatorname{co}\bigl(\mathcal Y_x(UT_xM)\bigr)$. In the navigation regimes of Subsections~\ref{subsec:zero-eta} and~\ref{subsec:zermelo}, this includes both the mild and the critical cases. 
\subsection{Scalar representation in the regular regime}
We now reduce the maximization in~\eqref{Hgen} to a one-dimensional problem. Let us introduce some notation.
For $p\in T_x^*M\setminus\{0\}$, let
\[
v_p:=p^{\sharp_h},	\qquad	r:=\alpha(v_p)=\alpha^*(p),	\qquad	e_p:=\frac{v_p}{r},
\]
and let us set
\[
\zeta_p:=h(Z_x,e_p)=\frac{p(Z_x)}{r},	\qquad	\chi_p:=	\alpha\bigl(Z_x-\zeta_pe_p\bigr).
\]
\begin{proposition}
\label{prop:Hscalar}
Fix $x\in M_+$ and assume that  $b_x<\frac12$ and $F_{\mathcal M}(-W_x)<1$.
Then, for $p\neq 0$,
\begin{equation}\label{Hscalar}
H(p)=\eta(x)\left[p(W_x)+r\cos\theta_p
\bigl(1+\zeta_p\cos\theta_p+\chi_p\sin\theta_p\bigr)\right],
\end{equation}
where $\theta_p=0$ if $\chi_p=0$, while if $\chi_p>0$, it is   
the unique solution in $(0,\pi/4)$ of 
\begin{equation}\label{Hcritical}
	-\sin\theta	-\zeta_p\sin(2\theta)+\chi_p\cos(2\theta)=0.
\end{equation}
\end{proposition}
\begin{proof}
We suppress the base point $x$ from the notation. Assume first that 	$\chi_p=0$. Then $Z=\zeta_pe_p$. Writing 	$t=h(e_p,u)\in[-1,1]$, we obtain 
\[
\bigl(1+h(Z,u)\bigr)p(u)=rt(1+\zeta_pt).
\]
Since $|\zeta_p|<1/2$,
\[
\frac{d}{dt}\bigl[rt(1+\zeta_pt)\bigr]	=	r(1+2\zeta_pt)>0
\]
on $[-1,1]$. The maximum is therefore attained at $t=1$, and
\[
H_Z(p)	=	r(1+\zeta_p)	=	\alpha^*(p)+p(Z).
\]
This agrees with~\eqref{Hscalar} when $\theta_p=0$. In particular, if 	$Z=0$, then $H_Z(p)=\alpha^*(p)$. 

Assume now that $\chi_p>0$ and define
\[e_p^\perp	:=	\frac{Z-\zeta_pe_p}{\chi_p}.\]
Every $u\in UT_xM$ can be written as
\[
u=t_1e_p+t_2e_p^\perp+\xi,\qquad	\xi\perp\operatorname{span}\{e_p,e_p^\perp\},
\]
with
\[t_1^2+t_2^2+\alpha^2(\xi)=1.\]
Since
\[p(u)=rt_1,\qquad h(Z,u)=\zeta_pt_1+\chi_pt_2,
\]
the function to be maximized is 
\begin{equation}\label{t1t2}
rt_1\bigl(1+\zeta_pt_1+\chi_pt_2\bigr).	
\end{equation} 
The value at $u=e_p$ is
\[
r(1+\zeta_p)\geq r(1-b_x)>0.
\]
Moreover, $1+h(Z,u)>0$ because $b_x<1/2$. Hence, no vector with 	$t_1\leq0$ can be a maximizer.
For fixed $t_1>0$, \eqref{t1t2} is strictly increasing with respect to 	$t_2$. Under the unit-sphere constraint, its maximum is therefore attained 	when
\[\xi=0, \qquad	t_2=\sqrt{1-t_1^2}.	\]
Writing
\[t_1=\cos\theta,	\qquad	t_2=\sin\theta,	\qquad	0\leq\theta\leq\frac{\pi}{2},
\]
we obtain
\[H_Z(p)=r\max_{0\leq\theta\leq\pi/2}g_p(\theta), \]
where
\[
g_p(\theta):=\cos\theta	\bigl(1+\zeta_p\cos\theta+\chi_p\sin\theta\bigr).
\]
Its derivative is
\[
g_p'(\theta)=-\sin\theta-\zeta_p\sin(2\theta)+\chi_p\cos(2\theta).
\]
On $(0,\pi/4)$, the equation of the critical points is equivalent to 
\[\chi_p=R_{\zeta_p}(\theta):=	\frac{\sin\theta+\zeta_p\sin(2\theta)}{\cos(2\theta)}.
\]
A direct computation gives
\[
R_{\zeta_p}'(\theta)=\frac{\cos\theta\bigl(\cos^2\theta+3\sin^2\theta\bigr)+2\zeta_p}
	{\cos^2(2\theta)}.
\]
For $0<\theta<\pi/4$, $\cos\theta\bigl(\cos^2\theta+3\sin^2\theta\bigr)>1$, and therefore
\[
R_{\zeta_p}'(\theta)>	\frac{1+2\zeta_p}{\cos^2(2\theta)}	>0.
\]
Furthermore,
\[
R_{\zeta_p}(0)=0,	\qquad	\lim_{\theta\to\pi/4^-}R_{\zeta_p}(\theta)=+\infty.
\]
Thus~\eqref{Hcritical} has a unique solution	$\theta_p\in(0,\pi/4)$.  For $\theta\in[\pi/4,\pi/2]$,
\[
g_p'(\theta)=-\sin\theta-\zeta_p\sin(2\theta)+\chi_p\cos(2\theta)\leq		-\sin\theta+|\zeta_p|\sin(2\theta)=	-\sin\theta	\bigl(1-2|\zeta_p|\cos\theta\bigr)<0.
\]
Consequently, $\theta_p$ is the unique global maximizer of $g_p$ on 	$[0,\pi/2]$.
Formula~\eqref{Hscalar} then follows from~\eqref{Hgen}.

\end{proof}
\subsection{Time-dependent laws and wind as a translational drift}
\label{subsec:additional-drift}
The preceding support Hamiltonian construction is fiberwise and extends without change to smoothly time-dependent data. Indeed, suppose that $\eta_t$, $\widehat Z_t$ and $\widehat W_t$ depend smoothly on time and consider
\[
\mathcal Y_t(u)= \bigl(\eta_t(x)+h(\widehat Z_t(x),u)\bigr)u+\widehat W_t(x),\qquad u\in UT_xM.
\]
If
\[
\mathcal B_{t,x}:=\operatorname{rad}\bigl(\mathcal Y_t(UT_xM)\bigr),
\]
then the same argument as in Proposition~\ref{support-general} gives
\begin{equation*}
H_t(p)=\max\left\{	0,p\bigl(\widehat W_t(x)\bigr)+\max_{u\in UT_xM}
\bigl(\eta_t(x)+h(\widehat Z_t(x),u)\bigr)p(u)\right\}.
\end{equation*}
In  the time-dependent model of Subsection~\ref{wildfireJPS}  wind enters through $\mathscr S_t$ in \eqref{direct} and therefore changes the shape of the local propagation set.

A different kinematic model consists in representing wind as an additive translational drift. Let $\mathcal Y^0_t$ be a time-dependent maximal velocity law such that
\[
\mathcal Y^0_t(UT_xM)=\partial\mathcal B^0_{t,x},
\]
where $\mathcal B^0_{t,x}$ is the unit ball of a Finsler norm $F^0_t$. Let $A=A(t,x)$ be a smooth time-dependent vector field and set
\begin{equation*}
\mathcal Y^A_t(u) 	:=	\mathcal Y^0_t(u)+A(t,x),	\qquad u\in UT_xM.
\end{equation*}
In accordance with \eqref{admissible}, the associated admissible set is
\[
\mathcal B^A_{t,x}:=\operatorname{rad}\bigl(\mathcal Y^A_t(UT_xM)\bigr)=
\operatorname{rad}\bigl(\partial\mathcal B^0_{t,x}+A(t,x)\bigr).
\]
Its support Hamiltonian is
\begin{equation}\label{Htime-drift}
H^A_t(p)=\max\left\{0,\,\bigl(F^0_t\bigr)^*(p)	+p\bigl(A(t,x)\bigr)\right\}.
\end{equation}
Thus, an additive drift contributes a linear term to the Hamiltonian, while the  maximum in \eqref{Htime-drift} takes into account the possibility, built into  \eqref{admissible}, of reducing the  velocity to zero.
At the kinematic level, $A$ may model a persistent mean wind, an airflow
induced by the terrain, or another systematic transport, and therefore it
should not necessarily be interpreted as the physical wind velocity.
Translations of local propagation sets also occur in classical wildfire
models, as for instance in Richards' model~\cite{Richards1990}.

It is important to distinguish the radial-hull convention above from the usual Zermelo-type translated differential inclusion
\begin{equation}\label{translated-ball}
\widetilde{\mathcal B}^{A}_{t,x}:=\mathcal B^0_{t,x}+A(t,x).
\end{equation}
In~\eqref{translated-ball}, the drift remains fully active when the controllable velocity inside $\mathcal B^0_{t,x}$ is reduced. Its support Hamiltonian is therefore
\[
\widetilde H^A_t(p)=\bigl(F^0_t\bigr)^*(p)+p\bigl(A(t,x)\bigr),
\]
without the maximum with zero. 
In the case where  $F^0_t\bigl(-A(t,x)\bigr)>1$, the translated set~\eqref{translated-ball} does not contain the origin, whereas the radial hull does. Hence, for the radial hull, every admissible trajectory reaching a point $x$ at time $s$ can be extended by the constant trajectory on any interval $[s,t]$.
Consequently, the corresponding reachable sets are nondecreasing by inclusion in time. 
However,  monotonicity of the reachable sets does not imply global reachability. The support Hamiltonian provides the following general obstruction, which applies to any of the time-dependent admissible velocity sets considered above.

\begin{proposition}[An obstruction to reachability]\label{prop:obstruction}
Let $\mathcal B_{t,x}\subset T_xM$ be a time-dependent family of admissible 	velocity sets and let
\[
H_t(p):=\max_{z\in\mathcal B_{t,x}}p(z), 	\qquad p\in T_x^*M,
\]
be its support Hamiltonian. Suppose that there exist a smooth function 	$f\colon M\to\mathbb R$ and a non-negative function 	$m\in L^1([0,+\infty))$ such that
\begin{equation*}
H_t(df_x)\leq m(t)
\end{equation*}
for every $(t,x)\in[0,+\infty)\times M$. Then every admissible trajectory 	$\gamma\colon[0,T]\to M$ satisfies, 
for almost every $t\in[0,T]$, 
\[
\frac{d}{dt}f(\gamma(t))=df_{\gamma(t)}(\dot\gamma(t))\leq	\max_{z\in\mathcal B_{t,\gamma(t)}}
df_{\gamma(t)}(z)=	H_t\bigl(df_{\gamma(t)}\bigr)\leq m(t).
\]
Consequently, a point $q\in M$ is not reachable from $x_0$ if
\[
f(q)-f(x_0)	>	\int_0^{+\infty}m(t)\,dt.
\]
\end{proposition}
\begin{proof}
The first equality is the chain rule, and the first inequality holds because $\dot\gamma(t)\in\mathcal B_{t,\gamma(t)}$ for almost every $t$. The last claim follows by integrating on $[0,T]$.
\end{proof}
For instance, let $M=\mathbb R^n$ with the Euclidean metric and consider
the autonomous law $\mathcal Y(u)=ru+A$, where $r>0$ and
$A\in\mathbb R^n$ is a constant vector such that $|A|\geq r$
(critical or strong wind). By~\eqref{Hgeneral-hat},
$H(p)=\max\{0,\,r|p|+p(A)\}$, and the linear function
$f(x):=-\langle A,x\rangle$ satisfies
$H(df)=\max\{0,\,|A|(r-|A|)\}=0$.
Hence, no point $q$ such that $\langle A,q-x_0\rangle<0$ is reachable
from $x_0$.

\section*{Acknowledgments}
\noindent E. Caponio is partially supported by MUR under the Program ``Department of Excellence'' Legge 232/2016 (Grant No. CUP - D93C23000100001).\\
E. Caponio and D. Corona are both partially supported by ``INdAM - GNAMPA Project'' CUP E53C25002010001.\\
All the authors gratefully acknowledge support from the Abdus Salam International Centre for Theoretical Physics, Trieste (Italy), and in particular the framework ``ICTP-INdAM Collaborative Grants and Research in Pairs Programme 2026'' during which the research for this paper was performed.

\appendix
\section{Fiberwise representation formulae}
\label{appendix}
While the gauge of the velocity set 
is generally determined only implicitly,
we show in this appendix that on $M_+$
the gauge formula admits a semi-explicit representation.
In particular, when $Z$ and $W$ are linearly dependent,
the gauge $F(y)$ depends on the norm of $y$
and on a scalar value through an inverse function.
In the non-collinear case, the inverse function
is instead defined on a closed elliptic disk of $\mathbb R^2$
(on its boundary, when $\dim M=2$).
The collinear and non-collinear cases are addressed
separately in the following subsections.

\subsection{Single vector field}\label{subsec:single} 

When $Z$ and $W$ are collinear,
for every $x\in M_+$ such that $Z_x\neq0$
the gauge function admits a fiberwise representation as an $(\alpha,\beta)$--metric, i.e.
\begin{equation}\label{Falpha-beta}
	F(y)=\alpha(y)\phi_x\left(\frac{\beta(y)}{\alpha(y)}\right),
	\qquad y\in T_xM\setminus\{0_x\}, 
\end{equation}
where $\phi_x$ is defined in terms of
the inverse of a scalar function
(see Theorem~\ref{thm:single} and Remark~\ref{beta}).
Whenever $Z_x = 0$,
then the velocity law reduces to the pure-translation case,
thus $F$ is explicitly computable
and is a Randers metric if $\alpha(W)<1$.
Therefore, $Z_x\neq0$ is the difficult case
when seeking a non-trivial representation formula.

It is convenient to use the notation of a single vector field.
In the one-field model~\eqref{eq:Y-oneV}, we have
$\widehat Z=(\mu-\sigma)V$ and $\widehat W=\sigma V$
by~\eqref{eq:hat-fields-one}; more generally, if $Z_x=a_1V_x$ and
$W_x=a_2V_x$ for a suitable $V_x\neq0_x$, the same holds at $x$ with
$\mu(x):=\eta(x)(a_1+a_2)$ and $\sigma(x):=\eta(x)a_2$.
In both cases, the
velocity law~\eqref{eq:generalForm} reads
\[
	\mathcal{Y}(u)
	= \Big(\eta(x) + \big(\mu(x) - \sigma(x)\big)h(V,u)\Big)u
	+ \sigma(x)V,
\]
and we set $\ell_x:=\alpha_x(V_x)$.
Although the representation~\eqref{Falpha-beta} is written in terms of
the one-form $\beta=-h(Z,\cdot)$ of~\eqref{Matsumoto}, it is convenient to
perform the scalar inversion using the one-form $\vartheta:=h(V,\cdot)$.
For $u\in UT_xM$, set $t:=\vartheta(u)=h(V_x,u)\in[-\ell_x,\ell_x]$.
The key observation is that the two quantities
\begin{align}
	\varrho_x(t)&:=	\alpha\bigl(\mathcal Y(u)\bigr)=	\sqrt{\bigl(\eta(x)+\mu(x)t\bigr)^2+\sigma(x)^2\bigl(\ell_x^2-t^2\bigr)},\nonumber\\
	\psi_x(t)&:=\frac{\vartheta\bigl(\mathcal Y(u)\bigr)}{\alpha\bigl(\mathcal Y(u)\bigr)}=
	\frac{\bigl(\mu(x)-\sigma(x)\bigr)t^2+\eta(x)t+\sigma(x)\ell_x^2}{\varrho_x(t)}\label{psi}
\end{align}
depend on $u$ only through $t$. On the other hand, for every
$y\neq0_x$, Lemma~\ref{radialgauge}-(i) gives $y/F(y)\in\partial\mathcal B_x
=\mathcal Y(UT_xM)$, that is, $y=F(y)\,\mathcal Y(u)$ for some
$u\in UT_xM$. Therefore, $F(y)=\alpha(y)/\varrho_x(t)$ and
$\psi_x(t)=\vartheta(y)/\alpha(y)$. Hence, inverting $\psi_x$ allows us to
recover $t$ from $\vartheta(y)/\alpha(y)$, and then to compute $F(y)$.

\begin{theorem}\label{thm:single}
	Assume that~\eqref{eq:FinslerM+conditions} holds,
	hence $F$ is a Finsler metric on $M_+$;
	in the present notation, this amounts to~\eqref{FMconvex}
	and~\eqref{Fconvex}.
	Then, for every $x\in M_+$ such that $V_x\neq 0$,
	the function
	$\psi_x\colon[-\ell_x,\ell_x]\to[-\ell_x,\ell_x]$, defined in \eqref{psi},	is an increasing diffeomorphism and,  for all $y\in T_xM\setminus \{0_x\}$, 
	\begin{equation}\label{Fvartheta}
		F(y)=\alpha(y)\widehat\phi_x\left(\frac{\vartheta(y)}{\alpha(y)}\right), 
	\end{equation}
	where $\widehat\phi_x(s):=1/\varrho_x\bigl(\psi_x^{-1}(s)\bigr)$, $s\in[-\ell_x,\ell_x]$.
\end{theorem}
\begin{proof}
	Let us fix $x\in M_+$ with $V_x\neq0$, and let us omit the dependence on $x$ from the 	notation.
	By~\eqref{Fconvex}, $\eta+\mu t>0$, for all  $t\in[-\ell,\ell]$,	and therefore
	$\varrho^2(t)	=	(\eta+\mu t)^2+\sigma^2(\ell^2-t^2)>0$. Moreover, from \eqref{psi}, 
	\[\psi(-\ell)=-\ell,	\qquad	\psi(\ell)=\ell.	\]
	Thus, it is enough to prove that $\psi'(t)>0$ on $[-\ell,\ell]$.
	Let us set
	\[
	\varphi(t):=\varrho^2(t),	\qquad	\zeta(t):=\psi(t)\varrho(t)	=	(\mu-\sigma)t^2+\eta t+\sigma\ell^2.
	\]
	Then
	\[
	\psi'(t)	=	\frac{2\zeta'(t)\varphi(t)-\zeta(t)\varphi'(t)}	{2\varrho^3(t)}.
	\]
	A direct computation gives
	\[
	2\zeta'(t)\varphi(t)-\zeta(t)\varphi'(t)	=	2\bigl(\eta+(\mu-\sigma)t\bigr)Q(t),
	\]
	where
	\begin{align*}
		Q(t)&=(\mu^2-\sigma^2)t^2+\eta(2\mu-\sigma)t+\eta^2+(\sigma-\mu)\sigma\ell^2\\
		&=(\eta+\mu t)\bigl(\eta+(\mu-\sigma)t\bigr)+\sigma(\mu-\sigma)(t^2-\ell^2).
	\end{align*}
	By~\eqref{FMconvex}, $\eta+(\mu-\sigma)t>0$, for all $t\in[-\ell,\ell]$, so it remains to show that $Q>0$ on $[-\ell,\ell]$. At the endpoints,
	\[	Q(\pm\ell)	=	(\eta\pm\mu\ell)	\bigl(\eta\pm(\mu-\sigma)\ell\bigr)>0\]
	by~\eqref{FMconvex} and~\eqref{Fconvex}. 
	
	If $\mu^2-\sigma^2\leq0$, then $Q$ is  concave, hence its 	minimum on $[-\ell,\ell]$ is attained at an endpoint, and therefore
	$Q>0$.
	
	Suppose now that $\mu^2-\sigma^2>0$. The only additional case to consider is when the vertex
	\[
	t^*	=	-\frac{\eta(2\mu-\sigma)}	{2(\mu^2-\sigma^2)}	\]
	belongs to $(-\ell,\ell)$. Then
	\[
	\frac{\eta^2(2\mu-\sigma)^2}	{4(\mu^2-\sigma^2)}	<	(\mu^2-\sigma^2)\ell^2,
	\]
	and hence
	\begin{align*}
		Q(t^*)&=\eta^2+\sigma(\sigma-\mu)\ell^2	-\frac{\eta^2(2\mu-\sigma)^2}{4(\mu^2-\sigma^2)}>	\eta^2-	\bigl[\mu^2-\sigma(2\sigma-\mu)\bigr]\ell^2.
	\end{align*}
	Thus it suffices to prove $\eta^2/\ell^2>\mu^2-\sigma(2\sigma-\mu)$.
	If $\sigma(2\sigma-\mu)\geq0$, this follows immediately from 	\eqref{Fconvex}.
	If $\sigma(2\sigma-\mu)<0$, then either $\sigma>0$ and $\mu>2\sigma>\sigma$, or $\sigma<0$ and $\mu<2\sigma<\sigma$; hence $(\mu-\sigma)(\mu-2\sigma)>0$,
	and consequently
	\[
	4(\mu-\sigma)^2	-	\bigl[\mu^2-\sigma(2\sigma-\mu)\bigr]	=	3(\mu-\sigma)(\mu-2\sigma)>0.
	\]
	The desired inequality now follows from~\eqref{FMconvex}.
	Hence $Q(t)>0$ also in this case and $\psi'(t)>0$ for $t\in[-\ell,\ell]$ and therefore 	$\psi$ is an increasing diffeomorphism.
	
	Let now $y\in T_xM\setminus\{0_x\}$. By Proposition~\ref{minimalgauge} we have $y/F(y)\in\partial\mathcal B_x$ and since,  from Proposition~\ref{Theta}, $\mathcal Y$ is a diffeomorphism  onto $\partial\mathcal B_x$, there is a unique $u\in UT_xM$ with $\mathcal Y(u)=y/F(y)$.
	Since $\vartheta$ and $\alpha$ are positively homogeneous, \eqref{psi} gives
	\[
	\psi_x\bigl(\vartheta(u)\bigr)=\frac{\vartheta\bigl(\mathcal Y(u)\bigr)}{\alpha\bigl(\mathcal Y(u)\bigr)}=
	\frac{\vartheta(y)}{\alpha(y)},
	\]
	and therefore $\vartheta(u)=\psi_x^{-1}\bigl(\vartheta(y)/\alpha(y)\bigr)$. Taking the $\alpha$-norms in $y=F(y)\,\mathcal Y(u)$ and using $\alpha(\mathcal Y(u))=\varrho_x(\vartheta(u))$, we conclude
	\[
	F(y)=\frac{\alpha(y)}{\varrho_x\Bigl(\psi_x^{-1}\bigl(\vartheta(y)/\alpha(y)\bigr)\Bigr)},
	\]
	which is \eqref{Fvartheta}.
\end{proof}
\begin{remark}\label{beta}
	Let us assume that $Z_x\neq 0$ and let us set 
	\[
	a_x:=\frac{\mu(x)-\sigma(x)}{\eta(x)}.
	\]
	Then from \eqref{Fvartheta}  $F$ equivalently admits the representation \eqref{Falpha-beta}
	where 
	\begin{equation}\label{phibeta}
		\phi_x(s):=	\widehat\phi_x\left(-\frac{s}{a_x}\right),\qquad s\in[-b_x,b_x].
	\end{equation}
	Indeed, we have
	\[Z_x=a_xV_x,\qquad\beta_x=-a_x\vartheta_x,\qquad b_x=|a_x|\ell_x.
	\]
	Since  $Z_x\neq0$, we have $a_x\neq0$ and 
	$\vartheta(y)/\alpha(y)=- \beta(y)/\big(a_x\,  \alpha(y)\big)$, so that \eqref{phibeta} and \eqref{Falpha-beta}
	follow.
\end{remark}	
\subsection{Two vector fields}
We now consider the case in which $Z_x$ and $W_x$ are linearly
independent, the collinear case being covered by
Subsection~\ref{subsec:single}. The idea is the same, but now both
$\alpha\bigl(\mathcal Y(u)\bigr)$ and the direction of $\mathcal Y(u)$
depend on $u$ through the two quantities $h(Z_x,u)$ and $h(W_x,u)$, so
that, instead of a scalar function, we have to invert a map of two
variables. Throughout this subsection, we fix $x\in M_+$ such that
\begin{equation}\label{twoconditions}
	b_x<\frac12,
	\qquad
	F_{\mathcal M}(-W_x)<1,
\end{equation}
and $Z_x$, $W_x$ are linearly independent.

Let us define
\[
\beta_1:=-h(Z,\cdot),
\qquad
\beta_2:=h(W,\cdot),
\]
and 
\[
P_x\colon T_xM\to\mathbb R^2,
\qquad
P_x(y):=\bigl(\beta_1(y),\beta_2(y)\bigr),
\qquad
\mathcal D_x:=P_x(UT_xM).
\]
Let 
\[
E_x:=\operatorname{span}\{Z_x,W_x\},
\qquad
L_x:=P_x|_{E_x}.
\]
Since $Z_x$ and $W_x$ are linearly independent, $L_x\colon E_x\to\mathbb R^2$ is a linear isomorphism. 
Writing $u=v+z$, with $v\in E_x$ and $z\in E_x^\perp$, we have
\[
P_x(u)=L_x(v), \qquad \alpha^2(z)=1-\alpha^2(v), \qquad u\in UT_xM. \]
Consequently, if $\dim M\geq3$, $\mathcal D_x	=	L_x\bigl(\{v\in E_x:\alpha(v)\leq1\}\bigr)$
is a closed elliptic disk, whereas, if $\dim M=2$, then $\mathcal D_x$ is the boundary of such a disk.

If $u\in UT_xM$ and $\mathbf s=P_x(u)$, then $h(Z_x,u)=-s_1$, $h(W_x,u)=s_2$,
and~\eqref{eq:normalizedForm} gives
\[
\alpha^2\bigl(\mathcal Y(u)\bigr)=\eta^2(x)\Bigl((1-s_1)^2+2(1-s_1)s_2+\alpha^2(W_x)\Bigr).
\]
So, for $\mathbf s=(s_1,s_2)\in\mathbb R^2$, let us set
\[
\chi_x(\mathbf s):=\sqrt{	(1-s_1)^2	+2(1-s_1)s_2	+\alpha^2(W_x)}.
\]
Wherever the radicand is positive, we define
\begin{equation}\label{Psirhoexplicit}
	\varrho_x(\mathbf s)	:=	\eta(x)\chi_x(\mathbf s),	\qquad	\Psi_x(\mathbf s)	:=	\frac{(1-s_1)\mathbf s+P_x(W_x)}
	{\chi_x(\mathbf s)},
\end{equation}
where, by definition of $P_x$, $P_x(W_x)=\bigl(-h(Z_x,W_x),\alpha^2(W_x)\bigr)$.
By~\eqref{twoconditions} and Proposition~\ref{generalcase}, $\mathcal Y(u)\neq0_x$ for every $u\in UT_xM$. Hence
$\chi_x>0$ on $\mathcal D_x$ and, by compactness, on an open neighborhood of $\mathcal D_x$. Thus $\varrho_x$ and $\Psi_x$ are smooth there. Moreover, by the linearity of $P_x$,
\[
P_x\bigl(\mathcal Y(u)\bigr) = \eta(x)\bigl((1-s_1)\mathbf s+P_x(W_x)\bigr),
\]
and therefore
\begin{equation}\label{Psirho}
	\varrho_x\bigl(P_x(u)\bigr)	=	\alpha\bigl(\mathcal Y(u)\bigr),	\qquad	\Psi_x\bigl(P_x(u)\bigr)
	=	P_x\bigl(\Theta_x(u)\bigr),
\end{equation}
where $\Theta_x$ is the map defined in Proposition~\ref{Theta}.
In particular, $\Psi_x(\mathcal D_x)\subset\mathcal D_x$.
To study $\Psi_x$, let us consider
\[
G_x^{Z,W}:=\bigl\{A\in O(T_xM,h_x):A(Z_x)=Z_x,\ A(W_x)=W_x\bigr\}.
\]
The $G_x^{Z,W}$-orbits in $UT_xM$ are precisely the fibers of $P_x|_{UT_xM}$. Indeed, $P_x(u)$ determines the $E_x$-component of $u$ through the isomorphism $L_x$, while two vectors in $E_x^\perp$ having the same $h_x$-norm are related by an element of $O(E_x^\perp,h_x)$. Conversely, every element of $G_x^{Z,W}$ leaves both $\beta_1$ and $\beta_2$ invariant. Moreover, $\Theta_x$ is $G_x^{Z,W}$-equivariant. In fact, for $A\in G_x^{Z,W}$,
\[
\mathcal Y(Au)=\eta(x)\Bigl(\bigl(1+h(Z_x,u)\bigr)Au+A(W_x)\Bigr)=A\mathcal Y(u),
\]
and therefore $\Theta_x(Au)=A\Theta_x(u)$. We are  ready to state the following theorem.
\begin{theorem}\label{thm:two}
	Assume~\eqref{twoconditions} and let $Z_x,W_x$ be linearly independent. 	Then
	$\Psi_x\colon\mathcal D_x\to\mathcal D_x$ 	is a diffeomorphism and
	$\varrho_x\colon\mathcal D_x\to(0,+\infty)$ is smooth. Moreover, for every $y\in T_xM\setminus\{0_x\}$,
	\begin{equation}\label{Ftwoforms}
		F(y)=\alpha(y)\,\Phi_x\left(\frac{\beta_1(y)}{\alpha(y)},\frac{\beta_2(y)}{\alpha(y)}\right),
	\end{equation}
	where
	\[
	\Phi_x(\mathbf s):=	\frac{1}{\varrho_x\bigl(\Psi_x^{-1}(\mathbf s)\bigr)},	\qquad	\mathbf s\in\mathcal D_x.
	\]
	In particular, $\Phi_x$ is smooth on $\mathcal D_x$ (up to its 	boundary when $\dim M\geq3$).
\end{theorem}
\begin{proof}
	By Proposition~\ref{Theta}, $\Theta_x$ is a diffeomorphism of 	$UT_xM$. Since both $\Theta_x$ and $\Theta_x^{-1}$ are
	$G_x^{Z,W}$-equivariant and the $G_x^{Z,W}$-orbits are precisely 	the fibers of $P_x|_{UT_xM}$, they induce mutually inverse maps on 	the orbit space. In view of~\eqref{Psirho}, the map induced by
	$\Theta_x$ is precisely $\Psi_x$. Hence $\Psi_x$ is bijective.
	
	It remains to prove the smoothness of its inverse. If $\dim M=2$, 	then $E_x=T_xM$ and 
	$P_x|_{UT_xM}\colon UT_xM\to\mathcal D_x$	is a diffeomorphism. Thus, 
	\[
	\Psi_x^{-1}	=	P_x\circ\Theta_x^{-1}	\circ\bigl(P_x|_{UT_xM}\bigr)^{-1},
	\]
	and the assertion follows. 
	
	Assume now $\dim M\geq3$. We show that any  smooth  function $f\colon UT_xM\to\R$ which is $G_x^{Z,W}$-invariant
	induces a smooth function on $\mathcal D_x$, up to its boundary. Choose $e\in E_x^\perp$ with $\alpha(e)=1$. Since
	$L_x=P_x|_{E_x}$ is an isomorphism, every $\mathbf s\in\mathcal D_x$ determines
	\[
	v=L_x^{-1}(\mathbf s)\in E_x, \qquad \alpha(v)\leq1.
	\]
	The fiber of $P_x|_{UT_xM}$ over $\mathbf s$ consists of the vectors
	$v+z$ with
	\[
	z\in E_x^\perp, \qquad\alpha^2(z)=1-\alpha^2(v),
	\]
	and is a single $G_x^{Z,W}$-orbit. Hence,  the function induced by $f$ on $\mathcal D_x$ is
	\[\bar f(\mathbf s)	=	f\left(	v+\sqrt{1-\alpha^2(v)}\,e\right),	\qquad 	v=L_x^{-1}(\mathbf s).
	\]
	The only point requiring attention is the smoothness at $\partial\mathcal D_x$, where $\alpha(v)=1$. Since $UT_xM$ is the
	$\alpha$-unit sphere, a smooth $G_x^{Z,W}$-invariant function $f\colon UT_xM\to\mathbb R$ admits the explicit smooth extension $\widetilde f(y):=f\left(y/\alpha(y)\right)$ to an annular neighborhood of $UT_xM$ in $T_xM$. Moreover,
	$\widetilde f$ is still $G_x^{Z,W}$-invariant, since every $A\in G_x^{Z,W}$ is an $h_x$-isometry and hence
	$\alpha(Ay)=\alpha(y)$. Fix $e\in E_x^\perp$ with $\alpha(e)=1$, and let $R\in G_x^{Z,W}$ be the orthogonal transformation which is the identity on $E_x$ and sends $e$ to $-e$. Then 
	\[\widetilde f(v+te)=\widetilde f\bigl(R(v+te)\bigr)=\widetilde f(v-te),\]
	so $(v,t)\mapsto\widetilde f(v+te)$ is smooth and even in $t$. Hence, by Whitney's theorem on differentiable even functions
	\cite{Whitney1943},  there exists locally a smooth function $g(v,r)$ such that
	\[
	\widetilde f(v+te)=g(v,t^2).
	\]
	Therefore, for $\mathbf s\in\mathcal D_x$, with $v=L_x^{-1}(\mathbf s)$,
	\[
	\bar f(\mathbf s) = g\left( L_x^{-1}(\mathbf s),1-\alpha^2\bigl(L_x^{-1}(\mathbf s)\bigr) \right),
	\]
	which is smooth also at $\partial\mathcal D_x$.

	We apply this observation to the two components of $P_x\circ\Theta_x^{-1}\colon UT_xM\to\R^2$.
	They are smooth and $G_x^{Z,W}$-invariant, and therefore induce a smooth map $\widehat\Psi_x\colon\mathcal D_x\to\mathcal D_x$. If $\mathbf s=P_x(w)$ and $u=\Theta_x^{-1}(w)$, then
	\[
	\Psi_x\bigl(\widehat\Psi_x(\mathbf s)\bigr)	=	\Psi_x\bigl(P_x(u)\bigr)	=	P_x\bigl(\Theta_x(u)\bigr)	=
	P_x(w)	=	\mathbf s.
	\]
	Hence, since $\Psi_x$ is bijective, $\widehat\Psi_x=\Psi_x^{-1}$, which implies that  $\Psi_x^{-1}$ is smooth and  
	$\Psi_x\colon\mathcal D_x\to\mathcal D_x$ is a diffeomorphism. Finally, the smoothness and positivity of
	$\varrho_x$ follow directly from~\eqref{Psirhoexplicit}.

	Let now $y\in T_xM\setminus\{0_x\}$. By 	Proposition~\ref{minimalgauge}, $y/F(y)\in\partial\mathcal B_x$,
	and Proposition~\ref{Theta} gives a unique $u\in UT_xM$ such that
	$\mathcal Y(u)=y/F(y)$.
	Therefore
	\[
	\Theta_x(u)	=	\frac{\mathcal Y(u)}{\alpha(\mathcal Y(u))}	=	\frac{y}{\alpha(y)}.
	\]
	Using~\eqref{Psirho}, we obtain
	\[
	\left(	\frac{\beta_1(y)}{\alpha(y)},	\frac{\beta_2(y)}{\alpha(y)}	\right)	=	P_x\bigl(\Theta_x(u)\bigr)
	=\Psi_x\bigl(P_x(u)\bigr),
	\]
	and hence
	\[
	P_x(u)	=	\Psi_x^{-1}\left(	\frac{\beta_1(y)}{\alpha(y)},	\frac{\beta_2(y)}{\alpha(y)}	\right).
	\]
	Taking $\alpha$-norms in 	$y=F(y)\mathcal Y(u)$ gives
	\[
	F(y)	=	\frac{\alpha(y)}	{\alpha\bigl(\mathcal Y(u)\bigr)}	=	\frac{\alpha(y)}{\varrho_x\bigl(P_x(u)\bigr)},
	\]
	which yields~\eqref{Ftwoforms}. Since $\Psi_x^{-1}$ and	$\varrho_x$ are smooth and $\varrho_x>0$, the same holds for
	$\Phi_x$.
\end{proof}
Thus $F$ admits a fiberwise $(\alpha,\beta_1,\beta_2)$--representation. In contrast with the collinear case of Subsection~\ref{subsec:single}, where the representation requires the inversion of the scalar function $\psi_x$, the non-collinear case requires the inversion of the two-dimensional map $\Psi_x$.

\end{document}